\documentclass{tm-amsart-latex}

\begin{document}

\title{
  Equivariant aspects of de-completing cyclic homology
}
  \tmnote{{\citetexmacs{TeXmacs:website}}}

\author{
  Zhouhang Mao
}
  \tmaffiliation{BCAM -- Basque Center for Applied Mathematics\\
  Alameda de Mazarredo 14, 48009 Bilbao, Spain}
  \tmemail{zhouhang.mao@gmail.com}

\subjclass[2020]{14F30, 14F40, 19D55}
\keywords{cyclic homology, de-completion, K-theory, de Rham cohomology, crystalline cohomology}

\date{\tmdate}

\maketitle

\begin{abstract}
  Derived de Rham cohomology turns out to be important in $p$-adic geometry,
  following Bhatt's discovery {\cite{Bhatt2012a}} of conjugate filtration in
  char $p$, de-Hodge-completing results in {\cite{Beilinson2012}}. In
  {\cite{Kaledin2018a}}, Kaledin introduced an analogous de-completion of the
  periodic cyclic homology, called the polynomial periodic cyclic homology,
  equipped with a conjugate filtration in char $p$, and expected to be related
  to derived de Rham cohomology.
  
  In this article, using genuine equivariant homotopy structure on Hochschild
  homology as in {\cite{Angeltveit2018,Blumberg2022}}, we give an equivariant
  description of Kaledin's polynomial periodic cyclic homology. This leads to
  Morita invariance without any Noetherianness assumption as in
  {\cite{Kaledin2018a}}, and the comparison to derived de Rham cohomology
  becomes transparent. Moreover, this description adapts directly to
  ``topological'' analogues, which gives rise to a de-Nygaard-completion of
  the topological periodic cyclic homology, which admits an extension to
  linear categories over truncated Brown--Peterson spectra.
  
  As an application, we establish a noncommutative crystalline--de Rham
  comparison, which de-completes the result in {\cite{Petrov2019}}, and
  extends it to prime $p = 2$. We also compare polynomial periodic cyclic
  homology to topological Hochschild homology over $\mathbb{F}_p$, and produce
  a conjugate filtration in char $p$ from our description.
\end{abstract}

\section{Introduction}

The main goal of this article is to give a genuine equivariant homotopy
theoretic approach to (noncommutative) de-(Hodge-)completion of negative and
periodic cyclic homology, which conceptually explains Kaledin's
{\tmdfn{co-periodic cyclic homology}} in {\cite{Kaledin2018a}}, giving rise to
motivic filtrations. This approach also adapts to topological versions, which
extends Devalapurkar--Hahn--Raksit--Yuan's $\tmop{THH}_{\Prism}$ to
DG-categories (or even linear categories over truncated Brown--Peterson
spectra). As an application, we give a de-completion of noncommutative
crystalline--de Rham comparison in {\cite{Petrov2019}}, which works for the
even prime $p = 2$ as well. We start with backgrounds of this de-completion in
algebraic geometry.

Grothendieck's algebraic de Rham cohomology, introduced in
{\cite{Grothendieck1966}}, turns out to be an important tool to study the
cohomology of smooth schemes. However, it does not behave well beyond the
smooth case. Illusie, following ideas of Quillen, introduced {\tmdfn{derived
de Rham cohomology}}, along with its Hodge-completion, in
{\cite[Ch.~VIII]{Illusie1972}}.

Hodge-completion makes derived de Rham cohomology easier to control, since the
associated graded pieces are given by shifts of (derived) exterior powers of
the cotangent complex. In particular, it coincides with algebraic de Rham
cohomology for smooth schemes.

On the other hand, (non-Hodge-completed) derived de Rham cohomology was
intractable until Bhatt's discovery in {\cite{Bhatt2012a}} of
{\tmdfn{conjugate filtration}} on it in char $p$, whose associated graded
pieces are equivalent to shifts of Frobenius twists of algebraic differential
forms. He also observed the triviality of derived de Rham cohomology after
rationalization. Using this new tool, he identified derived de Rham cohomology
with crystalline cohomology for lci maps between $\mathbb{Z}/ p^r$-schemes.
Later on, this non-Hodge-complete version becomes useful in $p$-adic geometry.
For example, Fontaine's period rings $A_{\tmop{cris}}$ and $C_{\tmop{st}}$ are
equipped with non-complete Hodge-filtration, and Bhatt applied this
non-Hodge-complete version to prove some of Beilinson's conjectures in
{\cite{Beilinson2012}}.

Periodic cyclic homology is a noncommutative counterpart of Hodge-completed
derived de Rham cohomology, defined for general DG-categories. For morphisms
of $\mathbb{Q}$-schemes, the relation is particularly simple: periodic cyclic
homology is a product of shifts of derived de Rham cohomology, as recently
shown by Konrad {\tmname{Bals}} in full generality in {\cite{Bals2024}}. This
relation was firstly discovered by Loday--Quillen for smooth morphisms of
$\mathbb{Q}$-schemes, cf.~{\cite[§5.1.12]{Loday1998}}, and studied in
{\cite{Toen2011}}. For schemes beyond char $0$, it was studied in
{\cite{Majadas1996}}. Breakthroughs were made in
{\cite{Bhatt2018,Antieau2019}}, which proved that, there is a complete
filtration on periodic cyclic homology, whose associated graded pieces are
shifts of Hodge-completed derived de Rham cohomology.

In view of usefulness of non-Hodge-completed derived de Rham cohomology in
$p$-adic geometry, it is natural to ask whether there is a
``non-Hodge-completion'' of periodic cyclic homology for DG-categories, which
carries a filtration with associated graded pieces being shifts of
non-Hodge-completed derived de Rham cohomology? In {\cite{Kaledin2018a}},
following Kontsevich's suggestion {\cite[2.32]{Kontsevich2008}}, Kaledin
defined {\tmdfn{polynomial periodic cyclic homology}}, equipped with a
conjugate filtration in char $p$ (when $p = 2$, it is later constructed in
{\cite{Kaledin2017}}), whose associated graded pieces are equivalent to shifts
of Frobenius twisted Hochschild homology. Using this, he showed that a certain
completion of polynomial periodic cyclic homology, called {\tmdfn{co-periodic
cyclic homology}}, is a derived Morita invariant when the base is Noetherian.
Moreover, he expected that polynomial periodic cyclic homology is closely
related to derived de Rham cohomology.

As explained in {\cite{Kaledin2018a}}, Kaledin's defines polynomial periodic
cyclic homology and deals with it by explicit manipulations of chain
complexes, which makes the arguments technical and difficult, and the
homotopy-theoretic functoriality of this construction becomes opaque. A main
goal of this article is to give a ``more invariant treatment'' as he wished in
the introduction, which overcomes these difficulties.

The key to our description is the genuine equivariant homotopy structure on
the usual $\mathbb{Z}$-linear Hochschild homology. More precisely, let
$\mathcal{C}$ be a DG-category. Then the Hochschild homology of $\mathcal{C}$,
being a Borel $\mathbb{T}$-equivariant $\mathbb{Z}$-module spectrum, has the
formula
\[ \tmop{HH} (\mathcal{C}/\mathbb{Z}) = \tmop{THH} (\mathcal{C})
   \otimes_{\tmop{THH} (\mathbb{Z})}^{\mathbb{L}} \mathbb{Z}. \]
This admits an obvious cyclonic (à la {\cite{Barwick2016}}) structure: the
$\mathbb{T}$-equivariant ring $\mathbb{Z}$ is the underlying object of the
constant Tambara functor $\underline{\mathbb{Z}}$, and the universal property
of $\tmop{THH}$ in {\cite{Angeltveit2018}} gives rise to a map $\tmop{THH}
(\mathbb{Z}) \rightarrow \underline{\mathbb{Z}}$ of
$\mathbb{T}$-$\mathbb{E}_{\infty}$-rings. This gives rise to an enhancement
\[ \tmop{HH} (\mathcal{C}/ \underline{\mathbb{Z}}) = \tmop{THH} (\mathcal{C})
   \otimes_{\tmop{THH} (\mathbb{Z})}^{\mathbb{L}} \underline{\mathbb{Z}} \]
as a $\underline{\mathbb{Z}}$-module in cyclonic spectra. This was generalized
to a genuine version of factorization homology in {\cite{Blumberg2022}}.
However, up to our knowledge, such a genuine equivariant homotopy structure on
$\mathbb{Z}$-linear Hochschild homology does not seem to be studied in the
literature. By definition, such an enhancement is a derived Morita invariant.

\begin{remark}
  This genuine equivariant homotopy structure has other applications. In a
  companion paper {\cite{Mao2024}}, we use a ``thickening'' of it to define
  prismatic Hochschild homology. In our paper {\cite{Mao-Witt-norm}}, we use a
  similar genuine equivariant structure to streamline Kaledin's
  Hochschild--Witt homology, a noncommutative counterpart of de Rham--Witt
  complex.
\end{remark}

It turns out that such a structure contains enough information to recover
polynomial periodic cyclic homology. Recall that, in terms of explicit chain
complexes, the usual Tate construction involves a product totalization in one
direction, and Kaledin's polynomial periodic cyclic homology is taking the
direct sum totalization instead, so that it has good colimit-preserving
properties. Inspired by this, we introduce the following definitions for all
primes:

\begin{definition}[\Cref{def:decpl-Tate,def:HP-poly}]
  The {\tmdfn{($\underline{\mathbb{Z}}$-)de-completed $\mathbb{T}$-Tate
  construction}} $(-)^{\theta_{\underline{\mathbb{Z}}} \mathbb{T}}$ is the
  filtered-colimit-preserving approximation of the composite functor
  \[ \tmop{Mod}_{\underline{\mathbb{Z}}} (\tmop{Sp}^{g^{<} \mathbb{T}})
     \longrightarrow \tmop{Sp}^{B\mathbb{T}} \xrightarrow{(-)^{t\mathbb{T}}}
     \tmop{Sp} . \]
  The {\tmdfn{polynomial periodic cyclic homology}} $\tmop{HP}^{\tmop{poly}}
  (\mathcal{C}/\mathbb{Z})$ of a DG-category $\mathcal{C}$ is defined to be
  $\tmop{HH} (\mathcal{C}/
  \underline{\mathbb{Z}})^{\theta_{\underline{\mathbb{Z}}} \mathbb{T}}$,
  applying the de-completed $\mathbb{T}$-Tate construction to Hochschild
  homology $\tmop{HH} (\mathcal{C}/ \underline{\mathbb{Z}})$.
\end{definition}

This is a point where we could explain a similar strategy does not work if we
only look at the usual non-genuine equivariant Hochschild homology: the
colimit-preserving approximation of the functor $(-)^{t\mathbb{T}} \of D
(\mathbb{Z})^{B\mathbb{T}} \rightarrow D (\mathbb{Z})$ is zero, as its value
on the compact generator $\mathbb{Z} [\mathbb{T}]$ is, more or less by
definition, zero.

The same construction works for any $t$-bounded (i.e. bounded with respect to
the $t$-structure) animated ring as base in place of $\mathbb{Z}$, as done in
the main text.

\begin{remark}
  The de-completed Tate construction depends on the choice of base. However,
  in some cases, it does not quite depend on that. We will prove relevant
  results in \Cref{sec:dR-prism-comp}.
\end{remark}

From this description, it is immediate that polynomial periodic cyclic
homology is rationally zero (\Cref{rem:decpl-Tate-rat}), since the functor
$\tmop{Mod}_{\underline{\mathbb{Z}}} (\tmop{Sp}^{g^{<} \mathbb{T}})
\rightarrow \tmop{Sp}^{B\mathbb{T}}$ becomes an equivalence after
rationalization. With slightly more efforts, we show that

\begin{theorem}[\Cref{prop:asm-eq-smooth}]
  Let $\mathcal{C}$ be a smooth and bounded DG-category\footnote{Smooth
  DG-categories are necessarily modules categories over smooth
  $\mathbb{E}_1$-rings, by {\cite[Prop~11.3.2.4]{Lurie2018}}. A DG-category is
  {\tmdfn{bounded}} if, for every pair $(x, y)$ of compact objects, the
  mapping $\mathbb{Z}$-module spectrum $\tmop{Hom} (x, y) \in D (\mathbb{Z})$
  has bounded $\tmop{Tor}$-amplitude (or equivalently, $t$-bounded, since
  $\mathbb{Z}$ is of finite flat dimension).}. Then the {\tmdfn{assembly map}}
  $\tmop{HP}^{\tmop{poly}} (\mathcal{C}/\mathbb{Z}) \rightarrow \tmop{HP}
  (\mathcal{C}/\mathbb{Z})$ is an equivalence after profinite completion.
\end{theorem}

Note that, for every quasicompact quasiseparated scheme $X$, its derived
category $D (X)$ is bounded, by {\cite[Cor~3.1.8]{Bondal2003}}. When $X$ is in
addition smooth, then its derived category $D (X)$ is also smooth. We refer to
{\cite{Orlov2016}} for general discussions. This proposition, along with
colimit-preservation of de-completed Tate construction, implies that, on
animated rings, the polynomial periodic cyclic homology is left Kan extended
from polynomial rings, and thus it coincides with various ad-hoc constructions
in the literature\footnote{Also compare with {\cite[§1]{Manam2024}}.}, such
as in {\cite{Bhatt2018,Antieau2020}}. Consequently, we address Kaledin's
expectation in the following.

\begin{theorem}[Construction~\ref{cons:HKR-fil}]
  Let $R$ be a commutative ring. Then there exists a functorial filtration
  $\tmop{Fil}_{\tmop{HKR}}^{\ast}$ on the profinite completion
  $\tmop{HP}^{\tmop{poly}} (R /\mathbb{Z})^{\wedge}$ of polynomial periodic
  cyclic homology with associated graded pieces equivalent to shifts of
  derived de Rham cohomology $\tmop{dR}_{R /\mathbb{Z}}$ after profinite
  completion.
\end{theorem}

Our description also suggests a ``topological'' analogue.

\begin{definition}[\Cref{def:TP-poly}]
  Let $S$ be a perfectoid ring, and $\mathcal{C}$ a DG-category over $S$. Then
  {\tmdfn{topological polynomial periodic cyclic homology}}
  $\tmop{TP}^{\tmop{poly} / S} (\mathcal{C})$ is defined to be $\tmop{THH}
  (\mathcal{C})^{\theta_{\tmop{THH} (S)} \mathbb{T}}$, defined by applying
  $\tmop{THH} (S)$-de-completed $\mathbb{T}$-Tate construction to topological
  Hochschild homology $\tmop{THH} (\mathcal{C})$.
\end{definition}

Note that this is an arena where explicit chain complex manipulations cannot
arrive. Previous results for polynomial periodic cyclic homology adapts to its
topological analogue as well:

\begin{theorem}[\Cref{cor:asm-T-eq-sm}]
  Let $S$ be a perfectoid ring, and $R$ a $p$-completely smooth and bounded
  DG-category over $S$. Then the assembly map $\tmop{TP}^{\tmop{poly} / S} (R)
  \rightarrow \tmop{TP} (R)$ is an equivalence after $(p, \ker
  (\theta))$-completion, where $\theta \of A_{\inf} (S) \twoheadrightarrow S$
  is Fontaine's map.
\end{theorem}

\begin{theorem}[Construction~\ref{cons:mot-fil}]
  Let $S$ be a perfectoid ring, and $R$ a commutative $S$-algebra. Then there
  exists a functorial filtration $\tmop{Fil}_M^{\ast}$ on the $(p, \ker
  (\theta))$-completed topological polynomial periodic cyclic homology
  $\tmop{TP}^{\tmop{poly} / S} (R)_{(p, \ker (\theta))}^{\wedge}$ with
  associated graded pieces equivalent to shifts of Frobenius twisted prismatic
  cohomology $\varphi_A^{\ast}  \Prism_{R / A}$ after $(p, \ker
  (\theta))$-completion, where $A \assign A_{\inf} (S)$.
\end{theorem}

In a forthcoming work of Devalapurkar--Hahn--Raksit--Yuan, they produce
de-completions of $\tmop{TC}^-$ and $\tmop{TP}$ which carry filtrations with
associated graded pieces given by {\tmdfn{absolute prismatic cohomology}} and
its Nygaard filtered pieces. Our techniques adapt to this situation as well,
giving rise to {\tmdfn{de-completed Borel completion}} $(-)^{\eta_{\tmop{THH}
(\mathbb{Z})}} \of \tmop{Mod}_{\tmop{THH} (\mathbb{Z})} (\tmop{Sp}^{g^{<}
\mathbb{T}}) \rightarrow \tmop{Mod}_{\tmop{THH} (\mathbb{Z})^h}
(\tmop{Sp}^{g^{<} \mathbb{T}})$.

\begin{theorem}[\Cref{cor:THH-prism-noncomm}]
  The functor
  \[ \tmop{THH}_{\Prism} \of \tmop{CAlg}_{\mathbb{Z}}^{\tmop{an}} \rightarrow
     \tmop{RMod}_{\tmop{THH} (\mathbb{Z})^h} (\tmop{Sp}^{g^{<}
     \mathbb{T}})_p^{\wedge} \]
  in Devalapurkar--Hahn--Raksit--Yuan extends to a localizing invariant over
  $\mathbb{Z}$, which is equivalent to $(\tmop{THH} (-)^{\eta_{\tmop{THH}
  (\mathbb{Z})}})_p^{\wedge}$. We can replace $\mathbb{Z}$ by truncated
  Brown--Peterson spectra $\tmop{BP} \langle n \rangle$ for $n \in \mathbb{N}$
  as well (\Cref{rem:THH-Prism-BP-n}).
\end{theorem}

Thus this extended $\tmop{THH}_{\Prism}$ is indeed a noncommutative
counterpart of (Nygaard-filtered) absolute prismatic cohomology. Actually, we
can recover $\tmop{TP}^{\tmop{poly} / S}$ from this new construction:

\begin{proposition}[\Cref{ex:TP-poly-de-cpl-Borel-cpl,cor:de-cpl-THH-S-THH-Z}]
  Let $S$ be a perfectoid ring, and $\mathcal{C}$ a DG-category over $S$. Then
  we have equivalences
  \begin{eqnarray*}
    \left( (\tmop{THH} (\mathcal{C})^{\eta_{\tmop{THH}
    (\mathbb{Z})}})^{C_{p^{\infty}}} \right)_{(p, v_1)}^{\wedge} & \simeq &
    \tmop{TC}^{-, \tmop{poly} / S} (\mathcal{C})_{(p, \xi)}^{\wedge},\\
    \left( ((\tmop{THH} (\mathcal{C})^{\eta_{\tmop{THH} (\mathbb{Z})}})^{\Phi
    C_p})^{C_{p^{\infty}} / C_p} \right)_{(p, v_1)}^{\wedge} & \simeq &
    \tmop{TP}^{\tmop{poly} / S} (\mathcal{C})_{(p, \xi)}^{\wedge} .
  \end{eqnarray*}
\end{proposition}

It seems slightly surprising that cyclotomic $\tmop{THH}$ (with $\tmop{THH}
(\mathbb{Z})$-module structure) is already enough to de-Nygaard-complete
topological periodic cyclic homology, but this phenomenon is demystified by
Efimov's rigidity of localizing motives (\Cref{rem:rigid-loc-mot}).

As for noncommutative geometry on its own, we first give a noncommutative
crystalline--de Rham comparison. Let $R$ be an animated ring. Then the
crystalline--de Rham comparison theorem tells us that the crystalline
cohomology of $R \otimes_{\mathbb{Z}}^{\mathbb{L}} \mathbb{F}_p$ over
$\mathbb{Z}_p$ is equivalent to the $p$-completion of the
(non-Hodge-completed) derived de Rham of $R$ over $\mathbb{Z}$. In
{\cite{Petrov2019}}, when $p$ is an odd prime, for DG-categories
$\mathcal{C}$, they give a comparison between $\tmop{TP} (\mathcal{C}
\otimes_{\mathbb{Z}}^{\mathbb{L}} \mathbb{F}_p)$ with $\tmop{HP}
(\mathcal{C}/\mathbb{Z})$, which corresponds to such a comparison under
certain completions. This result is improved by
{\cite[Thm~0.5.1]{Devalapurkar2025}} as well for odd primes $p$. Equipped with
previous constructions, adapting main ideas in {\cite{Petrov2019}} along with
an observation by A.~{\tmname{Raksit}}, we prove a comparison of de-completed
constructions which corresponds to crystalline--de Rham comparison on the
nose:

\begin{theorem}[\Cref{cor:Z-lin-cat-crys-comp}]
  Let $\mathcal{C}$ be a DG-category. Then there exists a lax symmetric
  monoidal (in $\mathcal{C}$) equivalence
  \[ \tmop{TP}^{\tmop{poly} /\mathbb{F}_p} (\mathcal{C}
     \otimes_{\mathbb{Z}}^{\mathbb{L}} \mathbb{F}_p)_p^{\wedge} \simeq
     \tmop{HP}^{\tmop{poly}} (\mathcal{C}/\mathbb{Z})_p^{\wedge} . \]
\end{theorem}

We then give two proofs for the following comparison, due to Kaledin in
{\cite[Cor~11.15]{Kaledin2020}}, but our proof is much simpler.

\begin{proposition}[\Cref{cor:poly-HP-THH}]
  Let $\mathcal{C}$ be a DG-category over $\mathbb{F}_p$. Then the polynomial
  periodic cyclic homology $\tmop{HP}^{\tmop{poly}}
  (\mathcal{C}/\mathbb{F}_p)$ is equivalent to $\tmop{THH} (\mathcal{C})
  [\sigma^{- 1}]$ as $\mathbb{Z}^{t\mathbb{T}}$-module spectra.
\end{proposition}

We also produce a conjugate filtration on polynomial periodic cyclic homology
in char $p$ in \Cref{sec:conj-fil}, and prove that

\begin{proposition}[\Cref{cor:conj-fil-complete-bnd-below-cats}]
  Let $k$ be a commutative $\mathbb{F}_p$-algebra, and $\mathcal{C}$ a
  DG-category over $k$. Then the conjugate filtration on
  $\tmop{HP}^{\tmop{poly}} (\mathcal{C}/ k)$ is complete in the following two
  cases:
  \begin{enumerate}
    \item $\mathcal{C}= D (R)$ for some $(- 1)$-connective
    $\mathbb{E}_1$-$k$-algebra $R$ (which includes all associative
    $k$-algebras $R$);
    
    \item $\mathcal{C}= D (X)$ for a quasicompact quasiseparated $k$-scheme
    $X$.
  \end{enumerate}
\end{proposition}

\begin{remark}
  In comparison to Efimov's refined negative cyclic homology and its
  continuation in Scholze's ongoing work\footnote{An abstract, along with
  recordings, of the talk can be found at
  \url{https://www.mpim-bonn.mpg.de/node/13359}.} on refined $\tmop{TC}^-$, as
  mentioned above, in \Cref{rem:rigid-loc-mot}, we explain that topological
  Hochschild homology as a cyclotomic spectrum already sees ``all'' $p$-adic
  formal information. However, their versions capture rigid analytic
  information. For example, for smooth $\mathbb{F}_p$-schemes, their versions
  see rigid cohomology. It might be worth understanding whether equivariant
  homotopy theory could say something for their versions as well.
\end{remark}

\begin{notation}
  Let $G$ be a finite group. We denote by $\tmop{Sp}^{g G}$ the symmetric
  monoidal $\infty$-category of genuine $G$-spectra, by
  $\tmop{Mack}_G^{\tmop{coh}} (k)$ the abelian category of $k$-linear
  cohomological $G$-Mackey functors (we omit $k$ when $k =\mathbb{Z}$), and by
  $\tmop{Sp}^{g^{<} \mathbb{T}}$ (resp. $\tmop{Sp}^{g_p \mathbb{T}}$) the
  symmetric monoidal $\infty$-category of cyclonic (resp. $p$-cyclonic)
  spectra as in {\cite{Barwick2016}}.
\end{notation}

\begin{acknowledgments*}
  We would like to thank Lukas {\tmname{Brantner}}, Christian
  {\tmname{Carrick}}, Kaif {\tmname{Hilman}}, Dmitry {\tmname{Kaledin}},
  Kirill {\tmname{Magidson}}, Alexander {\tmname{Petrov}}, Arpon
  {\tmname{Raksit}}, Maxime {\tmname{Ramzi}}, Georg {\tmname{Tamme}}, and
  Dmitry {\tmname{Vaintrob}}. We also thank Robert {\tmname{Burklund}} and
  Jeremy {\tmname{Hahn}}. We thank the anonymous referee for questions and
  comments. This project has received funding from the European Research
  Council (ERC) under the European Union's Horizon 2020 research and
  innovation programme (grant agreement No. 864145).
\end{acknowledgments*}

\section{A de-completion of Tate construction}

Let $k$ be a commutative ring, and $G$ a finite group or $\mathbb{T}$. Recall
that the $G$-Tate construction $(-)^{t G} \of D (k)^{B G} \rightarrow D (k)$
does not preserve filtered colimits. In this section, we try to
``de-complete'' it when the input is further equipped with a genuine
equivariant structure. When $G = C_p$, we will show that it can be expressed
in terms of the geometric fixed points. We will also establish a de-completed
version of the Tate orbit lemma.

\begin{definition}
  \label{def:decpl-Tate}Let $G$ be a finite group (resp. $\mathbb{T}$), and
  $A$ an $\mathbb{E}_1$-algebra in the symmetric monoidal $\infty$-category
  $\tmop{Sp}^{g G}$ (resp. $\tmop{Sp}^{g^{<} \mathbb{T}}$) of $G$-spectra
  (resp. cyclonic spectra). Then
  \begin{itemize}
    \item The {\tmdfn{($A$-)de-completed homotopy $G$-fixed points}}
    $(-)^{\eta_A G} \of \tmop{RMod}_A (\tmop{Sp}^{g G}) \rightarrow \tmop{Sp}$
    (resp. $\tmop{RMod}_A (\tmop{Sp}^{g^{<} \mathbb{T}}) \rightarrow
    \tmop{Sp}$) is the filtered-colimit-preserving approximation of the
    homotopy $G$-fixed points $(-)^{h G} \of \tmop{RMod}_A (\tmop{Sp}^{g G})
    \rightarrow \tmop{Sp}$ (resp. $\tmop{RMod}_A (\tmop{Sp}^{g^{<}
    \mathbb{T}}) \rightarrow \tmop{Sp}$), equipped with an {\tmdfn{assembly
    map}} $(-)^{\eta_A G} \rightarrow (-)^{h G}$. We omit ``$A$-'' when the
    context is clear.
    
    \item The {\tmdfn{($A$-)de-completed $G$-Tate construction}}
    $(-)^{\theta_A G} \of \tmop{RMod}_A (\tmop{Sp}^{g G}) \rightarrow
    \tmop{Sp}$ (resp. $\tmop{RMod}_A (\tmop{Sp}^{g^{<} \mathbb{T}})
    \rightarrow \tmop{Sp}$) is the filtered-colimit-preserving approximation
    of the $G$-Tate construction $(-)^{t G} \of \tmop{RMod}_A (\tmop{Sp}^{g
    G}) \rightarrow \tmop{Sp}$ (resp. $\tmop{RMod}_A (\tmop{Sp}^{g^{<}
    \mathbb{T}}) \rightarrow \tmop{Sp}$), equipped with an {\tmdfn{assembly
    map}} $(-)^{\theta_A G} \rightarrow (-)^{t G}$ which canonically fits into
    a commutative diagram
    \begin{equation}
      \begin{array}{ccc}
        (-)^{\eta_A G} & \longrightarrow & (-)^{h G}\\
        \longdownarrow &  & \longdownarrow\\
        (-)^{\theta_A G} & \longrightarrow & (-)^{t G}
      \end{array} . \label{eq:asm-fixed-pts-Tate}
    \end{equation}
  \end{itemize}
\end{definition}

\begin{remark}
  \label{rem:asm-htpy-fixed-pts-Tate-cart}Let $G$ be a finite cyclic group
  (resp. $\mathbb{T}$), and $A$ a commutative algebra in $\mathcal{E} \assign
  \tmop{Sp}^{g G}$ (resp. $\mathcal{E} \assign \tmop{Sp}^{g^{<} \mathbb{T}}$).
  Then the diagram \eqref{eq:asm-fixed-pts-Tate} is Cartesian in the
  $\infty$-category $\tmop{Fun}^{\tmop{Ex}} (\mathcal{E}, \tmop{Sp})$: the
  fiber of the canonical natural transformation $(-)^{h G} \rightarrow (-)^{t
  G}$ is $(-)_{h G} [d] \in \tmop{Fun}^{\tmop{Ex}} (\mathcal{E}, \tmop{Sp})$
  (where $d = 0$ if $G$ is a finite cyclic group, and $d = 1$ when $G
  =\mathbb{T}$), which preserves filtered colimits as well, and consequently,
  the diagram \eqref{eq:asm-fixed-pts-Tate} induces an equivalence on fibers
  of vertical arrows.
\end{remark}

\begin{remark}
  Let $G$ be a finite group (resp. $\mathbb{T}$), and $A \rightarrow B$ a map
  of $\mathbb{E}_1$-algebras in the symmetric monoidal $\infty$-category
  $\mathcal{E} \assign \tmop{Sp}^{g G}$ (resp. $\mathcal{E} \assign
  \tmop{Sp}^{g^{<} \mathbb{T}}$) of $G$-spectra (resp. cyclonic spectra). Then
  the assembly maps on compact objects of $\tmop{RMod}_B (\mathcal{E})$ induce
  ``relative'' assembly maps $(-)^{\eta_A G} \rightarrow (-)^{\eta_B G}$ and
  $(-)^{\theta_A G} \rightarrow (-)^{\theta_B G}$, which fit into a
  commutative diagram
  \[ \begin{array}{ccc}
       (-)^{\eta_A G} & \longrightarrow & (-)^{\eta_B G}\\
       \longdownarrow &  & \longdownarrow\\
       (-)^{\theta_A G} & \longrightarrow & (-)^{\theta_B G}
     \end{array} \]
  which is Cartesian by \Cref{rem:asm-htpy-fixed-pts-Tate-cart}.
\end{remark}

\begin{remark}
  \label{rem:borel-cpl-asm-eq}Let $G$ be a finite group (resp. $\mathbb{T}$),
  and $A$ an $\mathbb{E}_1$-algebra in the symmetric monoidal
  $\infty$-category $\mathcal{E} \assign \tmop{Sp}^{g G}$ (resp. $\mathcal{E}
  \assign \tmop{Sp}^{g^{<} \mathbb{T}}$) of $G$-spectra (resp. cyclonic
  spectra). Let $A \rightarrow A^h$ denote the Borel completion of $A$. Then
  it follows immediately from the definitions that the de-completed homotopy
  $G$-fixed points $(-)^{\eta_A G}$ factors as
  \[ \tmop{RMod}_A (\mathcal{E}) \xrightarrow{(-) \otimes_A^{\mathbb{L}} A^h}
     \tmop{RMod}_{A^h} (\mathcal{E}) \xrightarrow{(-)^{\eta_{A^h} G}}
     \tmop{Sp} \]
  and similarly for the de-completed $G$-Tate construction. Roughly speaking,
  it does not hurt to replace all genuine equivariant bases by their Borel
  completions. However, sometimes it seems to be convenient to consider
  genuine equivariant bases.
\end{remark}

\begin{remark}
  \label{rem:de-cpl-Tate-symm-mon-struct}Let $G$ be a finite group, and $A$ an
  $\mathbb{E}_{\infty}$-algebra in $\tmop{Sp}^{g G}$. Then the lax symmetric
  monoidal structure on the homotopy fixed points $(-)^{h G}$ (resp. the Tate
  construction $(-)^{t G}$) gives rise\footnote{Here the same argument does
  not work for $G =\mathbb{T}$, as in this case, the object $A \in
  \tmop{Mod}_A (\tmop{Sp}^{g G})$ is not compact.} to a lax symmetric monoidal
  structure on the de-completed homotopy fixed points $(-)^{\eta_A G}$ (resp.
  the de-completed Tate construction $(-)^{\theta_A G}$). The assembly maps
  are equipped with a lax symmetric monoidal structure as well. In particular,
  the objects $M^{\eta_A G}$ (resp. $M^{\theta_A G}$) carries a canonical
  $A^{h G}$-(resp. $A^{t G}$-)module structure, which is functorial in $M \in
  \tmop{Mod}_A (\tmop{Sp}^{g G})$.
\end{remark}

\begin{remark}
  In desirable situations, the $A$-de-completed homotopy $G$-fixed points
  (resp. $G$-Tate construction) does not quite depend on $A$. We discuss some
  independences of this form in \Cref{sec:dR-prism-comp}.
\end{remark}

\begin{remark}
  \label{rem:borel-rigid}Let $G$ be a finite group, and $A$ an
  $\mathbb{E}_{\infty}$-algebra in $\tmop{Sp}^{g G}$. Recall that the
  symmetric monoidal $\infty$-category $\tmop{Sp}^{g G}$ is rigid, thus the
  forgetful functor $\tmop{Mod}_A (\tmop{Sp}^{g G}) \rightarrow \tmop{Mod}_A
  (\tmop{Sp}^{B G})$ factors through the rigidification $(\tmop{Mod}_A
  (\tmop{Sp}^{B G}))^{\tmop{rig}}$ of the target. Then the $A$-de-completed
  $G$-Tate construction $(-)^{\theta_A G}$ coincides with the composite
  \[ \tmop{Mod}_A (\tmop{Sp}^{g G}) \longrightarrow (\tmop{Mod}_A
     (\tmop{Sp}^{B G}))^{\tmop{rig}} \xrightarrow{(-)^{t G, \tmop{rig}}}
     \tmop{Sp}, \]
  where $(-)^{t G, \tmop{rig}}$ is the composite functor
  \[ (\tmop{Mod}_A (\tmop{Sp}^{B G}))^{\tmop{rig}} \longhookrightarrow
     \tmop{Ind} (\tmop{Mod}_A (\tmop{Sp}^{B G})) \xrightarrow{\tmop{Ind}
     ((-)^{t G})} \tmop{Ind} (\tmop{Sp}) \xrightarrow{\tmop{colim}} \tmop{Sp}
     . \]
  This can be checked by restricting to compact objects of $\tmop{Mod}_A
  (\tmop{Sp}^{g G})$. The same holds for de-completed homotopy $G$-fixed
  points.
\end{remark}

\begin{remark}
  \label{rem:Tate-Cn-mod-n}Let $A$ be an $\mathbb{E}_{\infty}$-algebra in $D
  (\mathbb{Z}) \otimes \tmop{Sp}^{g^{<} \mathbb{T}}$. Recall that, for every
  $A$-module $M$ in $\tmop{Sp}^{g^{<} \mathbb{T}}$ and every positive integer
  $n \in \mathbb{N}_{> 0}$, the canonical maps
  \[ M^{t\mathbb{T}} / n \longleftarrow M^{t\mathbb{T}}
     \otimes_{A^{t\mathbb{T}}} A^{t C_n} \longrightarrow M^{t C_n} \]
  are equivalences, {\cite[Lem~IV.4.12]{Nikolaus2018}}. It follows that, the
  canonical maps
  \[ M^{\theta_A \mathbb{T}} / n \longleftarrow M^{\theta_A \mathbb{T}}
     \otimes_{A^{t\mathbb{T}}} A^{t C_n} \longrightarrow M^{\theta_A C_n} \]
  are equivalences as well. Consequently, up to profinite completion, the
  $A$-de-completed $\mathbb{T}$-Tate construction $M^{\theta_A \mathbb{T}}$
  can be recovered from de-completed $C_n$-Tate constructions, where $n$ runs
  through all positive integers. This gives rise to a lax symmetric monoidal
  structure on the profinitely completed de-completed $\mathbb{T}$-Tate
  construction $(M^{\theta_A \mathbb{T}})^{\wedge}$.
\end{remark}

\begin{remark}
  \label{rem:decpl-Tate-rat}Let $k$ be a commutative ring. Then the constant
  Green functor $\underline{k}$ can be viewed as an object of $D (\mathbb{Z})
  \otimes \tmop{Sp}^{g^{<} \mathbb{T}}$. Then for every $n \in \mathbb{N}_{>
  0}$, we have
  \[ (\underline{k} \otimes \Sigma_{\mathbb{T}}^{\infty}  [\mathbb{T}/
     C_n]_+)^{t\mathbb{T}} \otimes_{\mathbb{Z}}^{\mathbb{L}} \mathbb{Q} \simeq
     0. \]
  Consequently, the rationalized de-completed $\mathbb{T}$-Tate construction
  $(-)^{\theta_{\underline{k}} \mathbb{T}}$ vanishes. Combining with
  \Cref{rem:Tate-Cn-mod-n}, we see that the de-completed $\mathbb{T}$-Tate
  construction $(-)^{\theta_{\underline{k}} \mathbb{T}}$ can be completely
  recovered from de-completed $C_n$-Tate constructions, where $n$ runs through
  all positive integers.
\end{remark}

The Tate orbit lemma admits a de-completion. We first observe that, since
$C_p$ is a simple group\footnote{In this article, we only consider cyclic
group actions, but the argument works for any simple group.}, the de-complete
$C_p$-Tate construction has a fairly simple formula:

\begin{lemma}
  \label{lem:decpl-Cp-Tate}Let $A$ be an $\mathbb{E}_1$-algebra in the
  symmetric monoidal category $\tmop{Sp}^{g C_p}$ of $C_p$-spectra. Then the
  $A$-de-completed $C_p$-Tate construction $(-)^{\theta_A C_p}$ coincides with
  the composite functor
  \[ \tmop{RMod}_A (\tmop{Sp}^{g C_p}) \xrightarrow{(-)^{\Phi C_p}} D (A^{\Phi
     C_p}) \longrightarrow D (A^{t C_p}) \]
  where the second functor is the base change along the map $A^{\Phi C_p}
  \rightarrow A^{t C_p}$ of $\mathbb{E}_1$-rings. Moreover, when $A$ is an
  $\mathbb{E}_{\infty}$-algebra in $\tmop{Sp}^{g C_p}$, this identification is
  (lax) symmetric monoidal.
\end{lemma}

\begin{proof}
  It suffices to restrict to compact objects of $\tmop{RMod}_A (\tmop{Sp}^{g
  C_p})$. For compact objects, there are many ways to see this. For example,
  the functor $(-)^{t C_p}$ vanishes on induced $C_p$-spectra, thus it
  canonically factors through $D (A^{\Phi C_p})$ in $\Pr^L$ (here we use the
  fact that $C_p$ is a simple group), cf.~{\cite[§5]{Ayala2021}}, and then to
  see that the resulting functor $D (A^{\Phi C_p}) \rightarrow D (A^{t C_p})$
  in $\Pr^L$ coincides with the base change, it suffices to check on the
  generator, which is straightforward. The symmetric monoidal structure
  follows from a similar argument.
\end{proof}

\begin{corollary}
  \label{cor:decpl-Cp-Tate-str-cts}Let $A$ be an $\mathbb{E}_1$-algebra in the
  symmetric monoidal category $\tmop{Sp}^{g C_p}$ of $C_p$-spectra. Then the
  $A$-de-completed $C_p$-Tate construction $(-)^{\theta_A C_p} \of
  \tmop{RMod}_A (\tmop{Sp}^{g C_p}) \rightarrow D (A^{t C_p})$ is strongly
  continuous\footnote{A functor $F \of \mathcal{C} \rightarrow \mathcal{D}$ in
  $\Pr_{\tmop{St}}^L$ is {\tmdfn{strongly continuous}} if its right adjoint
  $F^R \of \mathcal{D} \rightarrow \mathcal{C}$ preserves filtered colimits.
  When $\mathcal{C}$ is compactly generated, it is equivalent to $F$
  preserving compact objects.}.
\end{corollary}

\begin{corollary}
  Let $A$ be an $\mathbb{E}_{\infty}$-algebra in the symmetric monoidal
  category $\tmop{Sp}^{g C_p}$ of $C_p$-spectra. Then the $A$-de-completed
  $C_p$-Tate construction $(-)^{\theta_A C_p}$ is symmetric monoidal.
\end{corollary}

\begin{lemma}
  \label{lem:decpl-Tate-orb-lem-T}Let $A$ be a bounded below
  $\mathbb{E}_1$-algebra in cyclonic spectra. Then the $p$-completed
  $A$-de-completed $\mathbb{T}$-Tate construction\footnote{Note that the
  forgetful functor $\tmop{Sp}^{g^{<} \mathbb{T}} \rightarrow \tmop{Sp}^{g_p
  \mathbb{T}}$ induces an equivalence on $p$-complete objects, thus we could
  work $p$-typically throughout.}
  \[ ((-)^{\theta_A \mathbb{T}})_p^{\wedge} \of \tmop{RMod}_A
     (\tmop{Sp}^{g^{<} \mathbb{T}}) \longrightarrow D
     (A^{t\mathbb{T}})_p^{\wedge}, \]
  factors as\footnote{Thanks to \Cref{rem:borel-cpl-asm-eq}, it does not
  matter what cyclonic structure we put on $A^{t C_p}$.}
  \[ \tmop{RMod}_A (\tmop{Sp}^{g^{<} \mathbb{T}}) \xrightarrow{(-)^{\theta_A
     C_p}} \tmop{RMod}_{A^{t C_p}} (\tmop{Sp}^{g_p (\mathbb{T}/ C_p)})
     \xrightarrow{\left( (-)^{\eta_{A^{t C_p}} (\mathbb{T}/ C_p)}
     \right)_p^{\wedge}} D (A^{t\mathbb{T}})_p^{\wedge} \]
  in $\Pr_{\tmop{St}}^L$. 
\end{lemma}

\begin{proof}
  Since every functor preserves filtered colimits, it suffices to check on
  compact objects in $\tmop{RMod}_A (\tmop{Sp}^{g^{<} \mathbb{T}})$, and this
  follows from \Cref{cor:decpl-Cp-Tate-str-cts} and the Tate orbit lemma
  (since $A$ is bounded below, so is $A \otimes \Sigma_{\mathbb{T}}^{\infty} 
  [\mathbb{T}/ C_m]_+$ for every $m \in \mathbb{N}_{> 0}$).
\end{proof}

The same argument works for finite cyclic groups (for which we can further
keep track of the lax symmetric monoidal structure on the de-completed Tate
construction).

\begin{lemma}
  \label{lem:decpl-Tate-orb-lem-fin-cyc}Let $r \in \mathbb{N}_{> 0}$, and $A$
  a bounded below $\mathbb{E}_1$-(resp. $\mathbb{E}_{\infty}$-)algebra in
  $C_{p^r}$-spectra. Then the $A$-de-completed $C_{p^r}$-Tate construction
  \[ (-)^{\theta_A C_p} \of \tmop{Mod}_A \left( \tmop{Sp}^{C_{p^r}} \right)
     \longrightarrow D \left( A^{t C_{p^r}} \right), \]
  as a presentable (resp. lax symmetric monoidal) functor, factors as
  \[ \tmop{Mod}_A \left( \tmop{Sp}^{g C_{p^r}} \right)
     \xrightarrow{(-)^{\theta_A C_p}} \tmop{Mod}_{A^{t C_p}} \left(
     \tmop{Sp}^{g_p (C_{p^r} / C_p)} \right) \xrightarrow{(-)^{\eta_{A^{t
     C_p}} (C_{p^r} / C_p)}} D \left( A^{t C_{p^r}} \right) . \]
\end{lemma}

\begin{remark}
  Let $k$ be a commutative ring. By
  \Cref{rem:Tate-Cn-mod-n,rem:decpl-Tate-rat,lem:decpl-Tate-orb-lem-fin-cyc},
  we can also keep track of the lax symmetric monoidal structure on the
  $p$-completed de-completed $\mathbb{T}$-Tate construction
  $((-)^{\theta_{\underline{k}} \mathbb{T}})_p^{\wedge}$, showing that, as an
  exact lax symmetric monoidal functor, it factors through the composite exact
  symmetric monoidal functor
  \[ \tmop{Mod}_{\underline{k}} (\tmop{Sp}^{g^{<} \mathbb{T}})
     \xrightarrow{(-)^{\Phi C_p}} \tmop{Mod}_{\underline{k}^{\Phi C_p}}
     (\tmop{Sp}^{g^{<} (\mathbb{T}/ C_p)}) \longrightarrow \tmop{Mod}_{k^{t
     C_p}} (\tmop{Sp}^{g^{<} (\mathbb{T}/ C_p)}), \]
  where the remaining functor $\tmop{Mod}_{k^{t C_p}} (\tmop{Sp}^{g^{<}
  (\mathbb{T}/ C_p)}) \rightarrow D (k^{t\mathbb{T}})_p^{\wedge}$ is the limit
  of de-completed homotopy $C_{p^r}$-fixed points along $r \in \mathbb{N}$
  (which at least a priori does not necessarily preserve filtered colimits).
\end{remark}

\section{Polynomial periodic and negative cyclic homology}

We briefly review Kaledin's {\tmdfn{polynomial periodic cyclic homology}} of
cyclic objects, and then describe it in terms of the de-completed
$\mathbb{T}$-Tate construction (\Cref{def:HP-poly}), informally explaining why
it coincides with Kaledin's original construction.

Let $k$ be a commutative ring, and $X_{\bullet} \of \Lambda^{\tmop{op}}
\rightarrow \tmop{Mod}_k$ a cyclic objects in $k$-modules. Recall that, for
every $[n]_{\Lambda} \in \Lambda^{\tmop{op}}$, the $k$-module $X_n \assign
X_{\bullet} ([n]_{\Lambda})$ carries an $k$-linear $C_n$ action, which gives
rise to a $2$-periodic complex
\[ \begin{array}{cccccccccc}
     & \cdots & \xleftarrow{1 - \sigma_n} & X_n & \xleftarrow{N_n} & X_n &
     \xleftarrow{1 - \sigma_n} & X_n & \xleftarrow{N_n} & \cdots\\
     \text{weight} &  &  & - 1 &  & 0 &  & 1 &  & 
   \end{array} \]
of $k$-modules, where $\sigma_n \of X_n \rightarrow X_n$ is the generator of
$C_n$, and $N_n \assign 1 + \sigma_n + \cdots + \sigma_n^{n - 1}$ is the
$C_n$-norm. This complex represents the shifted Tate construction $X_n^{t C_n}
[- 1]$. These complexes compile into a double complex
\[ \begin{array}{ccccccccc}
     &  & \vdots &  & \vdots &  & \vdots &  & \\
     &  & \longdownarrow &  & \longdownarrow &  & \longdownarrow &  & \\
     \cdots & \xleftarrow{1 - \sigma} & X_1 & \xleftarrow{N} & X_1 &
     \xleftarrow{1 - \sigma} & X_1 & \xleftarrow{N} & \cdots\\
     &  & \longdownarrow &  & \longdownarrow &  & \longdownarrow &  & \\
     \cdots & \xleftarrow{1 - \sigma} & X_0 & \xleftarrow{N} & X_0 &
     \xleftarrow{1 - \sigma} & X_0 & \xleftarrow{N} & \cdots
   \end{array} \]
where we suppress the subscripts of $N$ and $\sigma$, and vertical
differentials are appropriately $2$-periodically given (for details, see
{\cite[§3.1]{Kaledin2018a}}). Then
\begin{itemize}
  \item The periodic cyclic homology $\tmop{HP} (X_{\bullet} / k)$ of the
  cyclic object $X_{\bullet}$, is the object in the derived category $D (k)$
  of $k$-modules represented by the product totalization of this double
  complex.
  
  \item The polynomial periodic cyclic homology $\tmop{HP}^{\tmop{poly}}
  (X_{\bullet} / k)$ of the cyclic object $X_{\bullet}$ is the object in the
  derived category $D (k)$ of $k$-modules represented by the direct sum
  totalization of this double complex.
\end{itemize}
We now give an alternative, more conceptual description of this double complex
and the polynomial periodic cyclic homology. Recall that the geometric
realization $| X_{\bullet} |_{\Lambda} \in D (k)^{B\mathbb{T}}$ can be
rewritten as a geometric realization
\[ \underset{[n] \in \tmmathbf{\Delta}^{\tmop{op}}}{\tmop{colim}}
   \tmop{Ind}_{C_n}^{\mathbb{T}} (X_n) \in D (k)^{B\mathbb{T}} . \]
If we apply the $\mathbb{T}$-Tate construction to it, and {\tmem{incorrectly}}
interchange $(-)^{t\mathbb{T}}$ with $\tmop{colim}_{[n]}$, we get
\[ \underset{[n] \in \tmmathbf{\Delta}^{\tmop{op}}}{\tmop{colim}} 
   (\tmop{Ind}_{C_n}^{\mathbb{T}} (X_n))^{t\mathbb{T}} \in D (k) \]
where $(\tmop{Ind}_{C_n}^{\mathbb{T}} (X_n))^{t\mathbb{T}} \simeq X_n^{t C_n}
[- 1]$ {\cite[Prop~3]{Hesselholt2020}}, thus we see that this computes the
polynomial periodic cyclic homology. Note that every $G$-module $M$ gives rise
to a cohomological $G$-Mackey functor $\underline{M} \of [G / H] \mapsto M^H$,
we can endow $C_n$-module $X_n$ with a cohomological $C_n$-Mackey functor
structure, and thus realize the geometric realization $| X_{\bullet}
|_{\Lambda}$ as an object of $\tmop{Mod}_{\underline{k}}
(\tmop{Sp}^{\mathbb{T}})$, and apply the de-completed $\mathbb{T}$-Tate
construction $(-)^{\theta_{\underline{k}} \mathbb{T}}$ to it, obtaining
polynomial periodic cyclic homology.

\begin{remark}
  This procedure can be made more rigorous by considering cohomological Mackey
  functors over Connes' cyclic category $\Lambda^{\tmop{op}}$. Since we do not
  depend on Kaledin's original construction, we skip such a development.
  However, a toy version of this is explained in \Cref{app:Tate-cohom-cplx}.
\end{remark}

\begin{question}
  This comparison does not compare the lax symmetric monoidal structure of
  Kaledin's polynomial periodic cyclic homology and ours. How do we compare
  the ring structure on the two?
\end{question}

Now we give our formal definitions, and explain why it corresponds to the
construction above for associative flat $k$-algebras.

{\construction{\label{cons:cyclonic-HH}Let $k$ be a $t$-bounded\footnote{This
means that it is bounded with respect to the canonical $t$-structure. We add
the prefix ``$t$-'' to avoid confusion with boundedness of $p$-power torsion
(which we do not use in this article anyways).} animated ring. The universal
property of $\tmop{THH}$ as in {\cite{Angeltveit2018}} gives rise to a map
$\tmop{THH} (k) \rightarrow \underline{k}$ of
$\mathbb{T}$-$\mathbb{E}_{\infty}$-rings. Thus for every $\tmop{THH}
(k)$-module $M$ in cyclonic spectra, we get an object $M \otimes_{\tmop{THH}
(k)}^{\mathbb{L}} \underline{k} \in \tmop{Mod}_{\underline{k}}
(\tmop{Sp}^{g^{<} \mathbb{T}})$. In particular, let $\mathcal{C}$ be a
dualizable presentable stable $k$-linear $\infty$-category, we have a
canonical genuine equivariant enhancement of the $k$-linear Hochschild
homology of $\mathcal{C}$, denoted by $\tmop{HH} (\mathcal{C}/ \underline{k})
\in \tmop{Mod}_{\underline{k}} (\tmop{Sp}^{g^{<} \mathbb{T}})$.}}

\begin{definition}
  \label{def:HP-poly}Let $k$ be a $t$-bounded animated ring, and $\mathcal{C}$
  a dualizable presentable stable $k$-linear $\infty$-category. Then the
  {\tmdfn{polynomial periodic cyclic homology}} $\tmop{HP}^{\tmop{poly}}
  (\mathcal{C}/ k)$ (resp. the {\tmdfn{polynomial negative cyclic homology}}
  $\tmop{HC}^{-, \tmop{poly}} (\mathcal{C}/ k)$) is defined to be the
  $\underline{k}$-de-completed $\mathbb{T}$-Tate construction $\tmop{HH}
  (\mathcal{C}/ \underline{k})^{\theta_{\underline{k}} \mathbb{T}}$ (resp. the
  $\underline{k}$-de-completed homotopy $\mathbb{T}$-fixed points $\tmop{HH}
  (\mathcal{C}/ \underline{k})^{\eta_{\underline{k}} \mathbb{T}}$).
\end{definition}

Let $k$ be a commutative ring, and $R$ an associative flat $k$-algebra. Recall
that, as in {\cite[§6]{Angeltveit2018}}, the Hochschild homology $\tmop{HH}
(R / \underline{k})$ relative to the constant Tambara functor $\underline{k}$
is informally given by the geometric realization of relative norms $[n]
\mapsto R^{\otimes_{\underline{k}}^{\mathbb{L}} C_n}$. In
{\cite{Mao-Witt-norm}}, we formally identify these norms with the naive tensor
powers equipped with the obvious cyclic action. Thus the previous discussion
basically identifies the polynomial periodic cyclic homology
$\tmop{HP}^{\tmop{poly}} (R / k)$ with Kaledin's original one.

\section{de Rham and prismatic comparison}\label{sec:dR-prism-comp}

In this section, we will first establish a comparison result
(\Cref{prop:asm-eq-smooth}), which says that, for smooth commutative algebras,
polynomial periodic cyclic homology is the same as periodic cyclic homology.
As a consequence, polynomial periodic cyclic homology of commutative algebras
acquires a filtration whose associated graded pieces are equivalent to shifts
of (non-Hodge-completed) derived de Rham cohomology. We then define a
``topological'' analogue (\Cref{def:TP-poly}), and for smooth algebras over
perfectoid rings, it acquires a {\tmdfn{motivic filtration}}
(Construction~\ref{cons:mot-fil}) whose associated graded pieces are
equivalent to shifts of (non-Nygaard-completed) prismatic cohomology.

We first give some sufficient conditions for assembly maps being equivalences.

\begin{lemma}
  \label{lem:asm-map-fin-gp-eq-projs}Let $A$ is a commutative algebra in
  cyclonic spectra (resp. $G$-spectra for a finite group $G$) whose underlying
  Borel equivariant spectrum, denoted by $A^h$, is $t$-bounded. Then the
  assembly map $(-)^{\theta_A G} \rightarrow (-)^{t G}$ (resp. $(-)^{\eta_A G}
  \rightarrow (-)^{h G}$) is an equivalence on the idempotent-complete stable
  subcategory of $\tmop{Mod}_A (\tmop{Sp}^{g^{<} \mathbb{T}})$ (resp.
  $\tmop{Mod}_A (\tmop{Sp}^{g G})$) generated by $A$-modules of the form
  $\bigoplus_{i \in I} A \otimes [\mathbb{T}/ H_i]$ (resp. $\bigoplus_{i \in
  I} A \otimes [G / H_i]$) for an indexed family $(H_i)_{i \in I}$ of finite
  cyclic groups $H_i$ (resp. finite groups $H_i \subseteq G$).
\end{lemma}

\begin{proof}
  We write the argument for the finite group case. The cyclonic case is
  similar.
  \begin{itemize}
    \item By construction, the assembly map is an equivalence on $A \otimes [G
    / H]$ for finite groups $H \subseteq G$.
    
    \item The family $\left\{ A^h \otimes [G / H] \in \tmop{Sp}^{B G}
    \barsuchthat \left( \text{$H \subseteq G$ is a finite subgroup} \right)
    \right\}$ is uniformly $t$-bounded, thus the canonical map
    \[ \bigoplus_i (A \otimes [G / H_i])^{h G} \longrightarrow \left(
       \bigoplus_i A \otimes [G / H_i] \right)^{h G} \]
    is an equivalence for an indexed family $(H_i \subseteq G)_i$ of finite
    subgroups of $G$. It follows that the assembly map is an equivalence on
    $\bigoplus_i A \otimes [G / H_i]$.
    
    \item The result follows from the fact that the functors $(-)^{\theta_A
    G}$, $(-)^{\eta G}$, $(-)^{t G}$ and $(-)^{h G}$ are exact.
  \end{itemize}
\end{proof}

When $G$ is a finite cyclic group, the situation is particularly simple.

\begin{corollary}
  \label{cor:fin-cyc-asm-eq-bnded}Let $G$ be a finite cyclic group. Then the
  assembly map $(-)^{\theta_{\underline{\mathbb{Z}}} G} \rightarrow (-)^{t G}$
  (resp. $(-)^{\eta_{\underline{\mathbb{Z}}} G} \rightarrow (-)^{h G}$) is an
  equivalence on $t$-bounded objects in $\tmop{Mod}_{\underline{\mathbb{Z}}}
  (\tmop{Sp}^{g G})$.
\end{corollary}

\begin{proof}
  Since the abelian category $\tmop{Mack}_G^{\tmop{coh}}$ has finite
  projective dimension when $G$ is finite cyclic {\cite[Cor~7.2]{Bouc2017}},
  every $t$-bounded object is represented by a finite complex of projective
  objects. Thus the assembly map $(-)^{\theta G} \rightarrow (-)^{t G}$ is an
  equivalence on these objects by \Cref{lem:asm-map-fin-gp-eq-projs}.
\end{proof}

\begin{corollary}
  \label{cor:fin-cyc-bnd-base-chg}Let $G$ be a finite cyclic group, and $A$ a
  commutative algebra in $D^b \tmop{Mack}_G^{\tmop{coh}} \subseteq
  \tmop{Mod}_{\underline{\mathbb{Z}}} (\tmop{Sp}^{g G})$. Then the assembly
  map $(-)^{\theta_{\underline{\mathbb{Z}}} G} \rightarrow (-)^{\theta_A G}$
  (resp. $(-)^{\eta_{\underline{\mathbb{Z}}} G} \rightarrow (-)^{\eta_A G}$)
  is an equivalence. Consequently, the $A$-de-completed $G$-Tate construction
  $(-)^{\theta_A G}$ (resp. the $A$-de-completed homotopy $G$-fixed points
  $(-)^{\eta_A G}$) coincides with the composite functor
  \[ \tmop{Mod}_A (\tmop{Sp}^{g G}) \longrightarrow
     \tmop{Mod}_{\underline{\mathbb{Z}}} (\tmop{Sp}^{g G})
     \xrightarrow{(-)^{\theta_{\underline{\mathbb{Z}}} G} \infixor
     (-)^{\eta_{\underline{\mathbb{Z}}} G}} D (\mathbb{Z}) . \]
\end{corollary}

\begin{proof}
  Note that the $\infty$-category $\tmop{Mod}_A (\tmop{Sp}^{g G})$ is
  generated by objects of the form $M \otimes_{\underline{\mathbb{Z}}}^L A$
  for finite permutation $G$-modules $M$, which is $t$-bounded since $M$ is
  $\underline{\mathbb{Z}}$-flat, and the result follows from
  \Cref{cor:fin-cyc-asm-eq-bnded}.
\end{proof}

We are ready to establish the comparison between the polynomial periodic
(resp. negative) cyclic homology and the periodic (resp. negative) cyclic
homology on smooth algebras:

\begin{proposition}
  \label{prop:asm-eq-smooth}Let $k$ be a $t$-bounded animated ring, and $R$ an
  $\mathbb{E}_1$-$k$-algebra. Then the commutative diagram
  \[ \begin{array}{ccc}
       \tmop{HC}^{-, \tmop{poly}} (R / k) & \longrightarrow & \tmop{HC}^- (R /
       k)\\
       \longdownarrow &  & \longdownarrow\\
       \tmop{HP}^{\tmop{poly}} (R / k) & \longrightarrow & \tmop{HP} (R / k)
     \end{array} \]
  as an instance of \eqref{eq:asm-fixed-pts-Tate} is Cartesian. If $R$ is
  $p$-completely smooth as an $\mathbb{E}_1$-$k$-algebra, and has bounded
  $\tmop{Tor}$-amplitude in $D (k)_p^{\wedge}$, then the horizontal assembly
  maps are equivalences after $p$-completion.
\end{proposition}

We first give a proof in the special case of $R$ being a $p$-completely smooth
animated $k$-algebra (i.e. with commutativity), since the general case needs
knowledge on polygonic spectra in {\cite{Krause2023}}, and the commutative
case is sufficient for this section.

\begin{proof*}{Proof of \Cref{prop:asm-eq-smooth} with commutativity}
  By \Cref{rem:asm-htpy-fixed-pts-Tate-cart}, this commutative diagram is
  Cartesian, thus for $p$-completely smooth animated $k$-algebras $R$, it
  suffices to check that the map $\tmop{HP}^{\tmop{poly}} (R / k) \rightarrow
  \tmop{HP} (R / k)$ is an equivalence after modulo $p$, which is equivalent
  to base change along $k \rightarrow k \otimes_{\mathbb{Z}}^{\mathbb{L}}
  \mathbb{F}_p$, thus we may assume that $k$ is a $t$-bounded animated
  $\mathbb{F}_p$-algebra. We can further check it after modulo $p$ again. By
  \Cref{rem:Tate-Cn-mod-n}, it reduces to check that the assembly map
  \[ \tmop{HH} (R / \underline{k})^{\theta_{\underline{k}} C_p}
     \longrightarrow \tmop{HH} (R / k)^{t C_p} \]
  is an equivalence. By
  \Cref{cor:fin-cyc-asm-eq-bnded,cor:fin-cyc-bnd-base-chg}, it suffices to
  check that the cyclonic spectrum $\tmop{HH} (R / \underline{k}) \in
  \tmop{Mod}_{\underline{k}} (\tmop{Sp}^{g^{<} \mathbb{T}})$ is $t$-bounded
  (thus so after forgetting to $\tmop{Sp}^{g C_p}$). We check in two steps.
  \begin{description}
    \item[Polynomial case] When $R = P \otimes_{\mathbb{F}_p}^{\mathbb{L}} k$
    where $P$ is a (finite) polynomial $\mathbb{F}_p$-algebra. Then we have
    $\tmop{HH} (R / \underline{k}) = \tmop{HH} (P / \underline{\mathbb{F}_p})
    \otimes_{\mathbb{F}_p}^{\mathbb{L}} k$, and by $t$-boundedness of $k$, it
    suffices to show that $\tmop{HH} (P / \underline{\mathbb{F}_p})$ is
    $t$-bounded. This follows from {\cite[2.2.4 \& 2.2.5]{Hesselholt1996}}.
    
    \item[General case] By passing to a Zariski cover, we may assume that
    there exists an étale map $S \rightarrow R$ where $S$ is a (finite)
    polynomial $k$-algebra. Then by {\cite[Add~3.2]{Hesselholt2020a}} (along
    with {\cite[Thm~B]{Borger2011}}, which is used in their proof), the map
    $\tmop{THH} (S) \rightarrow \tmop{THH} (R)$ is flat (even étale) in
    $\tmop{CAlg} (\tmop{Sp}^{g^{<} \mathbb{T}})$, thus so is the map
    $\tmop{HH} (S / \underline{k}) \rightarrow \tmop{HH} (R / \underline{k})$,
    and the result follows.
  \end{description}
\end{proof*}

Now we give the proof for \Cref{prop:asm-eq-smooth} in full generality. As in
the the first proof of \Cref{prop:asm-eq-smooth}, it reduces to check that the
assembly map
\[ \tmop{HH} (R / \underline{k})^{\theta_{\underline{k}} C_p} \longrightarrow
   \tmop{HH} (R / k)^{t C_p} \]
is an equivalence when $k$ is a $t$-bounded animated $\mathbb{F}_p$-algebra,
and $R$ is $p$-completely smooth as an $\mathbb{E}_1$-$k$-algebra with bounded
$\tmop{Tor}$-amplitude in $D (k)$. We prove a generalized version with
coefficients.

Let $R$ be an $\mathbb{E}_1$-ring spectrum, and $M$ an $R$-$R$-bimodule in
$\tmop{Sp}$. Then by {\cite[§6]{Krause2023}}, the topological Hochschild
homology $\tmop{THH} (R ; M)$ carries a canonical $p$-polygonic structure,
given by the sequence $\left( \tmop{THH} (R ; M^{\otimes_R^{\mathbb{L}} p^r})
\in \tmop{Sp}^{B C_{p^r}} \right)_{r \in \mathbb{N}}$, along with polygonic
Frobenius maps
\[ \tmop{THH} (R ; M^{\otimes_R^{\mathbb{L}} p^r}) \longrightarrow \tmop{THH}
   (R ; M^{\otimes_R^{\mathbb{L}} p^{r + 1}})^{t C_p} \]
in $\tmop{Sp}^{B C_{p^r}}$. Moreover, when all $\tmop{THH} (R ;
M^{\otimes_R^{\mathbb{L}} p^r})$ in question are bounded below (this is the
case when $M$ is a perfect $R$-$R$-bimodule in $\tmop{Sp}$), we get a sequence
$\left( \tmop{THH} (R ; M^{\otimes_R^{\mathbb{L}} p^r}) \in \tmop{Sp}^{g
C_{p^r}} \right)_{r \in \mathbb{N}}$ with equivalences
\[ \tmop{THH} (R ; M^{\otimes_R^{\mathbb{L}} p^{r + 1}})^{\Phi C_p}
   \xrightarrow{\simeq} \tmop{THH} (R ; M^{\otimes_R^{\mathbb{L}} p^r}) \]
of genuine $C_{p^r}$-spectra. This construction extends to the case without
bounded-below-ness, by a forthcoming work by Harpaz--Nikolaus--Saunier. All of
these constructions are functorial in $(R ; M)$.

There is also a forgetful functor from cyclotomic spectra to polygonic
spectra, and $\tmop{THH} (R ; R)$ as a polygonic spectrum is the same as the
underlying polygonic spectrum of the cyclotomic spectrum $\tmop{THH} (R)$.

{\construction{Let $k$ be an animated ring, $R$ an $\mathbb{E}_1$-$k$-algebra,
and $M$ an $R$-$R$-bimodule in $D (k)$. Then for every $r \in \mathbb{N}$, we
get a genuine equivariant enhancement $\tmop{HH} ((R ; M) / \underline{k})$ of
Hochschild homology $\tmop{HH} ((R ; M^{\otimes_R^{\mathbb{L}} p^r}) / k)$
with coefficients given by
\[ \tmop{HH} ((R ; M^{\otimes_R^{\mathbb{L}} p^r}) / \underline{k}) \assign
   \tmop{THH} (R ; M^{\otimes_R^{\mathbb{L}} p^r}) \otimes_{\tmop{THH}
   (k)}^{\mathbb{L}} \underline{k} . \]}}

\begin{lemma}
  \label{lem:asm-polygonic-HH-eq}Let $k$ be a $t$-bounded animated ring, $R$
  an $\mathbb{E}_1$-$k$-algebra with bounded $\tmop{Tor}$-amplitude in $D
  (k)$, and $M$ a perfect $R$-$R$-bimodule in $D (k)$. Then the assembly map
  \[ \tmop{HH} ((R ; M^{\otimes_R^{\mathbb{L}} p}) /
     \underline{k})^{\theta_{\underline{k}} C_p} \longrightarrow \tmop{HH} ((R
     ; M^{\otimes_R^{\mathbb{L}} p}) / k)^{t C_p} \]
  is an equivalence.
\end{lemma}

\begin{proof}
  The proof of {\cite[Prop~III.1.1]{Nikolaus2018}} (or
  {\cite[Prop~2.2.3]{Lurie2011a}}) implies that the functor $\tmop{HH} ((R ;
  (-)^{\otimes_R^{\mathbb{L}} p}) / k)^{t C_p}$ is exact. By
  \Cref{lem:decpl-Cp-Tate}, the symmetric monoidal structure on $(-)^{\Phi
  C_p}$, and the equivalence $\tmop{THH} (R ; M^{\otimes_R^{\mathbb{L}}
  p})^{\Phi C_p} \xrightarrow{\simeq} \tmop{THH} (R ; M)$, we see that the
  functor $\tmop{HH} ((R ; (-)^{\otimes_R^{\mathbb{L}} p}) /
  \underline{k})^{\Phi C_p}$ is exact as well. Thus, to see that the assembly
  map in question is an equivalence, it suffices to show that it is an
  equivalence when $M = R \otimes_k^{\mathbb{L}} R$, the free $R$-$R$-bimodule
  in $D (k)$ of rank $1$. In this case, for every $r \in \mathbb{N}$, we have
  an equivalence
  \[ \tmop{THH} (R ; (R \otimes_k^{\mathbb{L}} R)^{\otimes_R^{\mathbb{L}}
     p^r}) \simeq \tmop{THH} (k ; R^{\otimes_k^{\mathbb{L}} p^r}) \]
  in $\tmop{Mod}_{\tmop{THH} (k)} \left( \tmop{Sp}^{g C_{p^r}} \right)$, and
  after base change along $\tmop{THH} (k) \rightarrow \underline{k}$, it
  becomes the relative $C_{p^r}$-norm $R^{\otimes_{\underline{k}} C_{p^r}}$.
  Under this identification, the assembly map in question becomes the assembly
  map
  \[ (R^{\otimes_{\underline{k}}^{\mathbb{L}} C_p})^{\theta_{\underline{k}}
     C_p} \longrightarrow (R^{\otimes_k^{\mathbb{L}} C_p})^{t C_p}, \]
  which is an equivalence by \Cref{lem:p-th-pow-decpl-Tate}.
\end{proof}

\begin{lemma}
  \label{lem:p-th-pow-decpl-Tate}Let $k$ be a $t$-bounded animated ring, and
  $M$ a $k$-module spectrum of bounded $\tmop{Tor}$-amplitude in $D (k)$. Then
  the assembly map
  \[ (M^{\otimes_{\underline{k}}^{\mathbb{L}} C_p})^{\theta_{\underline{k}}
     C_p} \longrightarrow (M^{\otimes_k^{\mathbb{L}} C_p})^{t C_p} \]
  is an equivalence.
\end{lemma}

\begin{proof}
  Since both sides are exact in $M$, we may assume that $M$ is a flat
  $k$-module spectrum, which is a filtered colimit $\tmop{colim}_i M_i$ of
  finite free $k$-module spectra by Lazard's theorem
  {\cite[Thm~7.2.2.15]{Lurie2017}}. In this case, each derived cohomological
  $C_{p^r}$-Mackey functor $M_i^{\otimes_{\underline{k}}^{\mathbb{L}} C_p}$ is
  $t$-bounded by the $t$-boundedness of $k$, thus so is their filtered colimit
  $M^{\otimes_{\underline{k}}^{\mathbb{L}} C_p}$. The result follows from
  \Cref{cor:fin-cyc-asm-eq-bnded,cor:fin-cyc-bnd-base-chg}.
\end{proof}

\begin{proof*}{Proof of \Cref{prop:asm-eq-smooth} in general}
  It follows from \Cref{lem:asm-polygonic-HH-eq} by setting $M = R$.
\end{proof*}

\begin{question}
  Is there any categorical generalization of \Cref{lem:asm-polygonic-HH-eq},
  namely, replacing $R$ by dualizable presentable stable $k$-linear
  $\infty$-category $\mathcal{C}$ which is smooth with bounded
  $\tmop{Tor}$-amplitude\footnote{An attempt for this definition: a dualizable
  presentable stable $k$-linear $\infty$-category $\mathcal{C}$ has bounded
  $\tmop{Tor}$-amplitude if the coevaluation functor $\mathcal{C}^{\vee}
  \otimes_{D (k)} \mathcal{C} \rightarrow D (k)$ sends compact objects to
  objects of bounded $\tmop{Tor}$-amplitude in $D (k)$.}, and replacing $M$ by
  an colimit-preserving $k$-linear endofunctor $\mathcal{C} \rightarrow
  \mathcal{C}$?
\end{question}

Now we construct an HKR filtration on the polynomial periodic (resp. negative)
cyclic homology whose associated graded pieces are shifts of
(non-Hodge-completed) derived de Rham cohomology (resp. its Hodge-filtered
pieces), Hodge-de-completing {\cite[Thm~1.17]{Bhatt2018}} and
{\cite{Antieau2019}}.

{\construction{\label{cons:HKR-fil}{\dueto{HKR filtration}}Let $k$ be a
$t$-bounded animated ring. Then by {\cite{Raksit2020}} (which generalizes
{\cite{Antieau2019}}), for every smooth $k$-algebra $R$, there is an
exhaustive filtration $\tmop{Fil}_{\tmop{HKR}}^{\ast}$, functorial in $R$, on
the canonical $p$-completed map
\[ \tmop{HC}^- (R / k)_p^{\wedge} \longrightarrow \tmop{HP} (R / k)_p^{\wedge}
\]
whose $i$-th associated graded piece $\tmop{gr}_{\tmop{HKR}}^i$ is given by
the canonical map
\[ \tmop{Fil}_H^i \tmop{dR}_{R / k} [2 i]_p^{\wedge} \rightarrow \tmop{dR}_{R
   / k} [2 i]_p^{\wedge} . \]
By \Cref{prop:asm-eq-smooth} and sifted-colimit-preservation of the functors
$\tmop{CAlg}_k^{\tmop{an}} \rightarrow D (k), R \mapsto
\tmop{HP}^{\tmop{poly}} (R / k)$ (resp. $R \mapsto \tmop{HC}^{-, \tmop{poly}}
(R / k)$), for every animated $k$-algebra $R$, we get an exhaustive filtration
$\tmop{Fil}_{\tmop{HKR}}^{\ast}$ on the Cartesian square
\[ \begin{array}{ccc}
     \tmop{HC}^{-, \tmop{poly}} (R / k)_p^{\wedge} & \longrightarrow &
     \tmop{HC}^- (R / k)_p^{\wedge}\\
     \longdownarrow &  & \longdownarrow\\
     \tmop{HP}^{\tmop{poly}} (R / k)_p^{\wedge} & \longrightarrow & \tmop{HP}
     (R / k)_p^{\wedge}
   \end{array} \]
whose $i$-th associated graded piece $\tmop{gr}_{\tmop{HKR}}^i$ is given by
\[ \begin{array}{ccc}
     \tmop{Fil}_H^i \tmop{dR}_{R / k} [2 i]_p^{\wedge} & \longrightarrow &
     \tmop{Fil}_H^i  \widehat{\tmop{dR}}_{R / k} [2 i]_p^{\wedge}\\
     \longdownarrow &  & \longdownarrow\\
     \tmop{dR}_{R / k} [2 i]_p^{\wedge} & \longrightarrow &
     \widehat{\tmop{dR}}_{R / k} [2 i]_p^{\wedge}
   \end{array} \]
where $\widehat{\tmop{dR}}_{R / k}$ is the Hodge-completed derived de Rham
cohomology of $R / k$.}}

\begin{question}
  Is the HKR filtration in Construction~\ref{cons:HKR-fil} complete?
\end{question}

It is very natural to extend our definition to de-Nygaard-complete topological
periodic cyclic homology.

\begin{definition}
  \label{def:TP-poly}Let $k$ be a $t$-bounded $\mathbb{E}_{\infty}$-ring, and
  $\mathcal{C}$ a dualizable presentable stable $k$-linear $\infty$-category.
  Then the {\tmdfn{topological $k$-polynomial periodic cyclic homology}}
  $\tmop{TP}^{\tmop{poly} / k} (\mathcal{C})$ (resp. the {\tmdfn{topological
  $k$-polynomial negative cyclic homology $\tmop{TC}^{-, \tmop{poly} / k}
  (\mathcal{C})$}}) is defined to be the $\tmop{THH} (k)$-de-completed
  $\mathbb{T}$-Tate construction $\tmop{THH} (\mathcal{C})^{\theta_{\tmop{THH}
  (k)} \mathbb{T}}$ (resp. the $\tmop{THH} (k)$-de-completed homotopy
  $\mathbb{T}$-fixed points $\tmop{THH} (\mathcal{C})^{\eta_{\tmop{THH} (k)}
  \mathbb{T}}$).
\end{definition}

Topological $k$-polynomial periodic (resp. negative) cyclic homology, even
after $p$-completion, seems intractable in general, partially due to the
global nature of its prismatization. The situation is drastically simpler when
the ring $k = S$ is $p$-complete and perfectoid, thanks to the Bökstedt
periodicity of $\tmop{THH} (S)_p^{\wedge}$. We recollect some notations and
computations in {\cite[Prop~6.2~\&~6.3]{Bhatt2018}}. However, we view
$\tmop{THH} (S)^{t C_p}$ non-equivariantly as a $\tmop{TC}^- (S)$-module,
which follows more closely to the convention in
{\cite[Lem~2.1]{Riggenbach2022}}.

\begin{remark}[{\cite[Prop~6.2~\&~6.3]{Bhatt2018}}]
  \label{rem:THH-perfd}Let $S$ be a perfectoid ring, $A \assign A_{\inf} (S) =
  W (S^{\flat})$ with Frobenius endomorphism $\varphi \of A \rightarrow A$,
  and $\xi$ a chosen generator of the kernel $\ker (\theta)$ of Fontaine's map
  $\theta \of A \twoheadrightarrow S$. Then the commutative square
  \[ \begin{array}{ccc}
       \tmop{TC}^- (S)_p^{\wedge} & \xrightarrow{\varphi_p^{h\mathbb{T}}} &
       \tmop{TP} (S)_p^{\wedge}\\
       \longdownarrow &  & \longdownarrow\\
       \tmop{THH} (S)_p^{\wedge} & \xrightarrow{\varphi_p} & \tmop{THH} (S)^{t
       C_p}
     \end{array} \]
  is a pushout diagram of $\mathbb{E}_{\infty}$-rings, and its homotopy groups
  are given by
  \[ \begin{array}{ccc}
       A [u, v] / (uv - \xi) & \xrightarrow[v \mapsto \varphi (\xi) \sigma^{-
       1}]{u \mapsto \sigma} & A [\sigma^{\pm}]\\
       \longdownarrow &  & \longdownarrow\\
       R [u] = (A / \xi) [u] & \xrightarrow{u \mapsto \sigma} & (A / \varphi
       (\xi)) [\sigma^{\pm}]
     \end{array}, \]
  where $| u | = | \sigma | = 2$ and $| v | = - 2$, the vertical maps are
  $A$-linear, and the horizontal maps are $\varphi$-linear. The homotopy
  groups of the canonical map $\tmop{TC}^- (S)_p^{\wedge} \rightarrow
  \tmop{TP} (S)_p^{\wedge}$ is given by the $A$-linear map
  \[ A [u, v] / (uv - \xi) \xrightarrow[v \mapsto \sigma^{- 1}]{u \mapsto \xi
     \sigma} A [\sigma^{\pm}] . \]
\end{remark}

In particular, the $\tmop{THH} (S)$-module $S = \tmop{THH} (S) / u$ in
$\tmop{Sp}^{B\mathbb{T}}$ is perfect, which implies that

\begin{lemma}
  Let $S$ be a perfectoid ring, and $M$ a $\tmop{THH} (S)$-module in
  $\tmop{Sp}^{g^{<} \mathbb{T}}$. Then the canonical map
  \[ M^{\theta_{\tmop{THH} (S)} \mathbb{T}} \otimes_{\tmop{TP}
     (S)}^{\mathbb{L}} S^{t\mathbb{T}} \longrightarrow (M \otimes_{\tmop{THH}
     (S)} S)^{\theta_S \mathbb{T}} \]
  is an equivalence after $p$-completion.
\end{lemma}

\begin{corollary}
  \label{cor:asm-T-eq-sm}Let $S$ be a perfectoid ring, and $\mathcal{C}$ a
  dualizable presentable stable $S$-linear $\infty$-category. Then the
  commutative diagram
  \[ \begin{array}{ccc}
       \tmop{TC}^{-, \tmop{poly} / S} (\mathcal{C}) & \longrightarrow &
       \tmop{TC}^- (\mathcal{C})\\
       \longdownarrow &  & \longdownarrow\\
       \tmop{TP}^{\tmop{poly} / S} (\mathcal{C}) & \longrightarrow & \tmop{TP}
       (\mathcal{C})
     \end{array} \]
  as an instance of \eqref{eq:asm-fixed-pts-Tate} is Cartesian. If the
  assembly map $\tmop{HP}^{\tmop{poly}} (\mathcal{C}/ S) \rightarrow \tmop{HP}
  (\mathcal{C}/ S)$ is an equivalence after $p$-completion\footnote{By
  \Cref{prop:asm-eq-smooth}, this is the case when $\mathcal{C}= D (R)$ for
  some $p$-completely smooth $\mathbb{E}_1$-$S$-algebra $R$.}, then the
  horizontal assembly maps are equivalences after $(p, \ker
  (\theta))$-completion.
\end{corollary}

\begin{proof}
  By \Cref{rem:asm-htpy-fixed-pts-Tate-cart}, this commutative diagram is
  Cartesian. When the assembly map $\tmop{HP}^{\tmop{poly}} (\mathcal{C}/ S)
  \rightarrow \tmop{HP} (\mathcal{C}/ S)$ is $p$-completely an equivalence,
  then by \Cref{rem:THH-perfd}, it is the $\pmod{\ker (\theta)}$ reduction of
  the map $\tmop{TP}^{\tmop{poly} / S} (\mathcal{C}) \rightarrow \tmop{TP}
  (\mathcal{C})$, thus the later is $(p, \ker (\theta))$-completely an
  equivalence, and the result follows.
\end{proof}

Similarly to Construction~\ref{cons:HKR-fil}, we have

{\construction{\label{cons:mot-fil}{\dueto{Motivic filtration}}Let $S$ be a
perfectoid ring, and let $A \assign A_{\inf} (S)$. Then by {\cite{Bhatt2018}},
for every smooth $S$-algebra $R$, there is an exhaustive filtration
$\tmop{Fil}_M^{\ast}$, functorial in $R$, on the canonical $(p, \ker
(\theta))$-completed\footnote{The $p$-completed $\tmop{TC}^- (R)_p^{\wedge}$
and $\tmop{TP} (R)_p^{\wedge}$ are automatically $\ker (\theta)$-complete
(since they are Nygaard-complete), but it is conceptually better to phrase it
after $(p, \ker (\theta))$-completion, since the $p$-completed polynomial
versions might not be $\ker (\theta)$-complete.} map
\[ \tmop{TC}^- (R)_{(p, \ker (\theta))}^{\wedge} \longrightarrow \tmop{TP}
   (R)_{(p, \ker (\theta))}^{\wedge} \]
whose $i$-th associated graded piece $\tmop{gr}_M^i$ is given by the canonical
map
\[ \tmop{Fil}_N^i \varphi_A^{\ast}  \Prism_{R / A} \longrightarrow
   \varphi_A^{\ast}  \Prism_{R / A} \]
for the Frobenius twisted prismatic cohomology. By
\Cref{prop:asm-eq-smooth,cor:asm-T-eq-sm}, and sifted-colimit-preservation of
the functors $\tmop{CAlg}_S^{\tmop{an}} \rightarrow D (\tmop{TC}^- (S)), R
\mapsto \tmop{TP}^{\tmop{poly} / S} (R)$ (resp. $R \mapsto \tmop{TC}^{-,
\tmop{poly} / S} (R)$), for every animated $k$-algebra $R$, we get an
exhaustive filtration $\tmop{Fil}_M^{\ast}$ on the Cartesian square
\[ \begin{array}{ccc}
     \tmop{TC}^{-, \tmop{poly} / S} (R)_{(p, \ker (\theta))}^{\wedge} &
     \longrightarrow & \tmop{TC}^- (R)_{(p, \ker (\theta))}^{\wedge}\\
     \longdownarrow &  & \longdownarrow\\
     \tmop{TP}^{\tmop{poly} / S} (R)_{(p, \ker (\theta))}^{\wedge} &
     \longrightarrow & \tmop{TP} (R)_{(p, \ker (\theta))}^{\wedge}
   \end{array} \]
whose $i$-th associated graded piece is given by
\[ \begin{array}{ccc}
     \tmop{Fil}_N^i \varphi_A^{\ast}  \Prism_{R / A} & \longrightarrow &
     \tmop{Fil}_N^i \varphi_A^{\ast}  \widehat{\Prism}_{R / A}\\
     \longdownarrow &  & \longdownarrow\\
     \varphi_A^{\ast}  \Prism_{R / A} & \longrightarrow & \varphi_A^{\ast} 
     \widehat{\Prism}_{R / A}
   \end{array} . \]}}

\begin{remark}
  The construction of HKR filtration and motivic filtration on polynomial
  cyclic theories is quite formal: one only needs a proposition similar to
  \Cref{prop:asm-eq-smooth,cor:asm-T-eq-sm}. In particular, the above
  construction also adapts to the Breuil--Kisin case (as in
  {\cite[§11]{Bhatt2018}}) and the $q$-de Rham case.
\end{remark}

\begin{remark}
  The topological polynomial periodic cyclic homology should be comparable to
  $\tmop{TC}^{(- 1)}$ as introduced in {\cite[§1]{Manam2024}}. In
  \Cref{sec:THH-prism}, we will discuss and compare to the construction in the
  absolute case in an ongoing project of Devalapurkar--Hahn--Raksit--Yuan.
\end{remark}

A natural question is whether the previous picture extends to relative
prismatic cohomology over an arbitrary base prism? When the base prism is
transversal, we have the following expectation.

\begin{remark}
  In a companion paper {\cite{Mao2024}}, we defined {\tmdfn{prismatic
  Hochschild homology}} $\tmop{HH}^{\Prism} (\mathcal{C}/ A)$ for a
  transversal prism $(A, I)$ and a dualizable presentable stable $A /
  I$-linear $\infty$-category $\mathcal{C}$, and formulated an HKR-type
  conjecture for $p$-completely smooth $A / I$-algebras. If that conjecture
  holds, then the $p$-completed $\underline{(A, I)}^{C_{p^{r -
  1}}}$-de-completed $\mathbb{T}/ C_{p^{r - 1}}$-construction $\left(
  \tmop{HH}^{\Prism} (R / A)^{C_{p^{r - 1}}} \right)^{\theta_{\underline{(A,
  I)}^{C_{p^{r - 1}}}} (\mathbb{T}/ C_{p^{r - 1}})}$ for animated $(A /
  I)$-algebra would carry an exhaustive filtration with associated graded
  pieces equivalent to shifts of $\left( \varphi_A^{\ast}  \Prism_{R / A}
  \right) \otimes_A^{\mathbb{L}} (A / I_r)$. We will address this in the
  future.
\end{remark}

Up to our knowledge, it was not widely expected that the topological
Hochschild homology as a cyclotomic spectrum contains enough information for a
de-Nygaard-completion such as \Cref{def:TP-poly}. However, in view of Efimov's
rigidity of localizing motives, this is expected:

\begin{remark}[M.~{\tmname{Ramzi}}]
  \label{rem:rigid-loc-mot}Let $k$ be a commutative ring, and $G$ a finite
  group\footnote{As in \Cref{rem:de-cpl-Tate-symm-mon-struct}, this does not
  work for $G =\mathbb{T}$, as the symmetric monoidal $\infty$-category
  $\tmop{Sp}^{g^{<} \mathbb{T}}$ is not rigid, but only locally rigid.}, $A$
  an $\mathbb{E}_{\infty}$-algebra in $\tmop{Sp}^{g G}$, and $E \of
  \tmop{Cat}_k^{\tmop{perf}} \rightarrow \tmop{Mod}_A (\tmop{Sp}^{B G})$ a
  finitary symmetric monoidal localizing invariant. Then the symmetric
  monoidal functor $E$ factors uniquely through the presentably stable
  symmetric monoidal $\infty$-category $\tmop{Mot}_{\tmop{loc}, k}$, obtaining
  a functor $\tmop{Mot}_{\tmop{loc}, k} \rightarrow \tmop{Mod}_A (D (k)^{B
  G})$ in $\tmop{CAlg} (\Pr_{\tmop{St}}^L)$. Efimov's rigidity theorem, as
  mentioned in {\cite[Rem~4.3]{Efimov2024}}, tells us that the presentably
  stable symmetric monoidal $\infty$-category $\tmop{Mot}_{\tmop{loc}, k}$ is
  rigid, thus we get a unique strongly continuous functor
  $\tmop{Mot}_{\tmop{loc}, k} \rightarrow \tmop{Mod}_A (D (k)^{B
  G})^{\tmop{rig}}$.
  
  If we are in addition given a finitary symmetric monoidal factorization
  \[ \tmop{Cat}_k^{\tmop{perf}} \xrightarrow{\tilde{E}} \tmop{Mod}_A
     (\tmop{Sp}^{g G}) \longrightarrow \tmop{Mod}_A (D (k)^{B G}) \]
  of $E$, where the functor $\tilde{E}$ is localizing as well, then by
  \Cref{rem:borel-rigid}, we see that the strongly continuous functor
  $\tmop{Mot}_{\tmop{loc}, k} \rightarrow \tmop{Mod}_A (D (k)^{B
  G})^{\tmop{rig}}$ coincides with the composite functor
  \[ \tmop{Mot}_{\tmop{loc}, k} \xrightarrow{\tilde{E}} \tmop{Mod}_A
     (\tmop{Sp}^{g G}) \longrightarrow \tmop{Mod}_A (D (k)^{B G})^{\tmop{rig}}
     . \]
  Informally, this tells us that $\tilde{E}$ ``knows'' everything about
  refined $E$.
  
  Now we apply this to the finitary symmetric monoidal localizing invariant
  \[ \tilde{E} \assign \tmop{THH} \of \tmop{Cat}_k^{\tmop{perf}}
     \longrightarrow \tmop{CycSp}^{\tmop{gen}} \longrightarrow
     \tmop{Mod}_{\tmop{THH} (k)} (\tmop{Sp}^{g G}) \]
  for any finite cyclic group $G$. It follows that the refinement of
  $\tmop{Sp}^{B G}$-valued $\tmop{THH}$ is completely determined by the
  functor $\tilde{E}$, thus by the (genuine) cyclotomic $\tmop{THH}$. Roughly
  speaking, this implies that the cyclotomic $\tmop{THH}$ already contains
  ``all'' profinite or $p$-adic formal information, including any
  de-Nygaard-completion.
\end{remark}

\section{The absolute case}\label{sec:THH-prism}

In an ongoing work of Devalapurkar--Hahn--Raksit--Yuan, they propose another
de-completion $\tmop{THH}_{\Prism}$ of topological negative cyclic homology
for animated rings. In this section, we show that our techniques apply to this
situation as well, which leads to {\tmdfn{de-completed Borel completion}}
(\Cref{def:de-cpl-Borel-cpl}), giving rise to a localizing invariant on
$\mathbb{Z}$-linear categories corresponding to {\tmdfn{absolute prismatic
cohomology}}. We also compare de-completed Borel completion with de-completed
homotopy fixed points and de-completed Tate construction
(\Cref{rem:decpl-htpy-fixed-pts-Borel-cpl,rem:decpl-Tate-Borel-cpl,prop:de-cpl-Borel-Tate-t-bnded,cor:de-cpl-Borel-Tate-THH-Z}),
and the construction in this section with previous constructions
(\Cref{cor:de-cpl-THH-S-THH-Z,cor:de-cpl-Z-THH-Z}). We first recall their
definition. Actually, our argument shows that $\tmop{THH}_{\Prism}$ extends to
a localizing invariant for $\tmop{BP} \langle n \rangle$-linear categories
(\Cref{rem:THH-Prism-BP-n}), which is independent of choice of $n$.

\begin{notation}
  Let $G$ be a finite group (resp. $\mathbb{T}$). We will denote by
  \[ (-)^h \of \tmop{Sp}^{g G} \longrightarrow \tmop{Sp}^{g G} \]
  (resp.
  \[ (-)^h \of \tmop{Sp}^{g^{<} \mathbb{T}} \longrightarrow \tmop{Sp}^{g^{<}
     \mathbb{T}} \]
  ) the {\tmdfn{Borel completion}}, which has a lax symmetric monoidal
  structure.
\end{notation}

{\construction{{\dueto{Devalapurkar--Hahn--Raksit--Yuan}}The functor
\[ \tmop{THH}_{\Prism} \of \tmop{CAlg}_{\mathbb{Z}}^{\tmop{an}}
   \longrightarrow \tmop{Mod}_{\tmop{THH} (\mathbb{Z})^h} (\tmop{Sp}^{g^{<}
   \mathbb{T}}) \]
is defined to be the left derived functor of the functor
\begin{eqnarray*}
  \tmop{Poly}_{\mathbb{Z}}^{\tmop{fg}} & \longrightarrow &
  \tmop{Mod}_{\tmop{THH} (\mathbb{Z})^h} (\tmop{Sp}^{g^{<} \mathbb{T}})\\
  R & \longmapsto & \tmop{THH} (R)^h .
\end{eqnarray*}
By construction, we have a canonical map
\[ \tmop{THH}_{\Prism} (R) \longrightarrow \tmop{THH} (R)^h \]
in $\tmop{Mod}_{\tmop{THH} (\mathbb{Z})^h} (\tmop{Sp}^{g^{<} \mathbb{T}})$,
which is an equivalence when $R$ is a finitely generated polynomial ring.}}

\begin{remark}[Devalapurkar--Hahn--Raksit--Yuan]
  It follows immediately from the discussion of {\cite[§6.2]{Bhatt2022}}
  (which explains the result in {\cite{Bhatt2018}} for animated rings) that
  there exists a motivic filtration on $\tmop{THH} (R)^{C_{p^{\infty}}}$
  (resp. $(\tmop{THH} (R)^{\Phi C_p})^{C_{p^{\infty}} / C_p}$) whose
  associated graded pieces are given by twisted absolute Nygaard cohomology
  $\tmop{Fil}_N^{\ast}  \Prism_R \{ \ast \} \left[ 2 \mathord{\ast} \right]$
  (resp. twisted absolute prismatic cohomology $\Prism_R \{ \ast \} \left[ 2
  \mathord{\ast} \right]$).
\end{remark}

To extend $\tmop{THH}_{\Prism}$ to a localizing invariant, as in
\Cref{def:HP-poly}, the key is to de-complete Borel completion.

\begin{definition}
  \label{def:de-cpl-Borel-cpl}Let $G$ be a finite group, and $A$ an
  $\mathbb{E}_1$-algebra in the symmetric monoidal $\infty$-category
  $\tmop{Sp}^{g G}$ of $G$-spectra. Then the {\tmdfn{$A$-de-completed Borel
  completion}}, denoted by $(-)^{\eta_A}$, is the filtered-colimit-preserving
  approximation
  \[ \tmop{RMod}_A (\tmop{Sp}^{g G}) \longrightarrow \tmop{RMod}_{A^h}
     (\tmop{Sp}^{g G}) \]
  of the Borel completion. The same if we replace $G$ by $\mathbb{T}$, and the
  symmetric monoidal $\infty$-category $\tmop{Sp}^{g G}$ of $G$-spectra by
  that $\tmop{Sp}^{g^{<} \mathbb{T}}$ of cyclonic spectra.
\end{definition}

\begin{remark}
  \label{rem:de-cpl-Borel-cpl-depends-on-Ah}Let $G$ be a finite group (resp.
  $\mathbb{T}$), and $A$ an $\mathbb{E}_1$-algebra in the symmetric monoidal
  $\infty$-category $\tmop{Sp}^{g G}$ (resp. $\tmop{Sp}^{g^{<} \mathbb{T}}$)
  of $G$-spectra (resp. cyclonic spectra). As in \Cref{rem:borel-cpl-asm-eq},
  the $A$-de-completed Borel completion factors canonically through the base
  change along the Borel completion map $A \rightarrow A^h$.
\end{remark}

Actually, the de-completed Borel completion is simply a base change.

\begin{lemma}
  \label{lem:de-cpl-Borel-cpl-base-chg}Let $G$ be a finite group, and $A$ an
  $\mathbb{E}_1$-algebra in the symmetric monoidal $\infty$-category
  $\tmop{Sp}^{g G}$ of $G$-spectra. Then the $A$-de-completed Borel completion
  $(-)^{\eta_A}$ coincides with the base change functor
  \[ (-) \otimes_A^{\mathbb{L}} A^h \of \tmop{RMod}_A (\tmop{Sp}^{g G})
     \longrightarrow \tmop{RMod}_{A^h} (\tmop{Sp}^{g G}) . \]
  The same if we replace $G$ by $\mathbb{T}$, and $\tmop{Sp}^{g G}$ by
  $\tmop{Sp}^{g^{<} \mathbb{T}}$.
\end{lemma}

\begin{proof}
  We start with the case that $G$ is a finite group. By
  \Cref{rem:de-cpl-Borel-cpl-depends-on-Ah}, without loss of generality, we
  may assume that $A$ is Borel complete, and by definition, it suffices to
  check that, for every $G$-orbit $[G / H]$, the tensor product $[G / H]
  \otimes A$ is Borel complete. This follows from the dualizability of $[G /
  H]$ in $\tmop{Sp}^{g G}$, cf. the proof of {\cite[Ex~2.15]{Brantner2021}}.
  
  The case for $G =\mathbb{T}$ is similar: we have to show that $[\mathbb{T}/
  C_n] \otimes A$ is Borel complete, and it suffices to restrict to
  $\tmop{Sp}^{g C_m}$ for every $m \in \mathbb{N}_{> 0}$, and $[\mathbb{T}/
  C_n]$ is dualizable there as well.
\end{proof}

It then follows from \Cref{rem:borel-cpl-asm-eq} that de-completed homotopy
fixed points and de-completed Tate construction factor through de-completed
Borel completion, thus the latter contains more information than the former.
When the group is finite, this factorization takes a very simple form:

\begin{remark}
  \label{rem:decpl-htpy-fixed-pts-Borel-cpl}Let $G$ be a finite group, and $A$
  an $\mathbb{E}_1$-algebra in the symmetric monoidal $\infty$-category
  $\tmop{Sp}^{g G}$ of $G$-spectra. By restricting to compact objects in
  $\tmop{RMod}_A (\tmop{Sp}^{g G})$, we see that the $A$-de-completed homotopy
  $G$-fixed points $(-)^{\eta_A G}$ factors as
  \[ \tmop{RMod}_A (\tmop{Sp}^{g G}) \xrightarrow{(-)^{\eta_A}}
     \tmop{RMod}_{A^h} (\tmop{Sp}^{g G}) \xrightarrow{(-)^G} \tmop{RMod}_{A^{h
     G}} (\tmop{Sp}) . \]
\end{remark}

\begin{remark}
  \label{rem:decpl-Tate-Borel-cpl}Let $G$ be a finite group, and $A$ an
  $\mathbb{E}_1$-algebra in the symmetric monoidal $\infty$-category
  $\tmop{Sp}^{g G}$ of $G$-spectra. By restricting to compact objects in
  $\tmop{RMod}_A (\tmop{Sp}^{g G})$, we see that the $A$-de-completed $G$-Tate
  construction $(-)^{\theta_A G}$ factors as
  \[ \tmop{RMod}_A (\tmop{Sp}^{g G}) \xrightarrow{(-)^{\eta_A}}
     \tmop{RMod}_{A^h} (\tmop{Sp}^{g G}) \xrightarrow{(-)^{\Phi_{\{ e \}} G}}
     \tmop{RMod}_{A^{h G}} (\tmop{Sp}), \]
  where $(-)^{\Phi_{\{ e \}} G}$ is the geometric fixed points with respect to
  the family $\{ e \}$ of subgroups {\cite[Def~5.7]{Ayala2021}}. When $G$ is a
  cyclic $p$-group, we have $(-)^{\Phi_{\{ e \}} G} = (-)^{\Phi C_p}$.
\end{remark}

Actually, this factorization also holds when $G =\mathbb{T}$ for certain bases
which we establish in
\Cref{prop:de-cpl-Borel-Tate-t-bnded,cor:de-cpl-Borel-Tate-THH-Z}. Now we
formulate the main comparison theorem

\begin{theorem}
  \label{thm:THH-prism-comp}Let $R$ be a smooth
  $\mathbb{E}_1$-$\mathbb{Z}$-algebra of bounded
  $\tmop{Tor}$-amplitude\footnote{Since the projective dimension of
  $\mathbb{Z}$ is finite, this is the same as $t$-boundedness.}. Then the
  assembly map
  \[ \tmop{THH} (R)^{\eta_{\tmop{THH} (\mathbb{Z})}} \longrightarrow
     \tmop{THH} (R)^h \]
  in $\tmop{Mod}_{\tmop{THH} (\mathbb{Z})^h} (\tmop{Sp}^{g^{<} \mathbb{T}})$
  is an equivalence after $p$-completion.
\end{theorem}

Before giving its proof, we mention an immediate consequence: since finitely
generated polynomial rings are smooth, it follows that

\begin{corollary}
  \label{cor:THH-prism-noncomm}After $p$-completion, the functor
  $\tmop{THH}_{\Prism}$ factors as
  \[ \tmop{CAlg}_{\mathbb{Z}}^{\tmop{an}} \xrightarrow{\tmop{Perf} (-)}
     \tmop{Cat}_{\mathbb{Z}}^{\tmop{perf}} \longrightarrow
     \tmop{Mod}_{\tmop{THH} (\mathbb{Z})^h} (\tmop{Sp}^{g^{<}
     \mathbb{T}})_p^{\wedge}, \]
  where the second functor is the composite functor
  \[ \tmop{Cat}_{\mathbb{Z}}^{\tmop{perf}} \xrightarrow{\tmop{THH}}
     \tmop{Mod}_{\tmop{THH} (\mathbb{Z})} (\tmop{Sp}^{g^{<} \mathbb{T}})
     \xrightarrow{((-)^{\eta_{\tmop{THH} (\mathbb{Z})}})_p^{\wedge}}
     \tmop{Mod}_{\tmop{THH} (\mathbb{Z})^h} (\tmop{Sp}^{g^{<}
     \mathbb{T}})_p^{\wedge}, \]
  being the proposed localizing invariant extending $\tmop{THH}_{\Prism}$.
\end{corollary}

The strategy to prove \Cref{thm:THH-prism-comp} is similar to that for
\Cref{cor:asm-T-eq-sm}. However, the cyclonic spectrum $\tmop{THH}
(\mathbb{Z})$ is more complicated than $\tmop{THH} (S)$ for perfectoid rings
$S$. Nevertheless, the following two computational inputs, that we learnt from
Robert {\tmname{Burklund}} and Jeremy {\tmname{Hahn}}, are sufficient.

\begin{lemma}[{\cite{Boekstedt1994}}]
  \label{lem:THH-Z-mod-p}Let $p$ be a prime. Then we have
  \[ \pi_{\ast} (\tmop{THH} (\mathbb{Z}) / p) \cong \mathbb{F}_p [\mu_1]
     \otimes \Lambda (\lambda_1) \]
  with $\deg \mu_1 = 2 p$ and $\deg \lambda_1 = 2 p - 1$. Moreover, the
  cyclotomic Frobenius $\tmop{THH} (\mathbb{Z}) \rightarrow \tmop{THH}
  (\mathbb{Z})^{t C_p}$ realizes $\pi_{\ast} (\tmop{THH} (\mathbb{Z})^{t C_p}
  / p)$ as $\pi_{\ast} (\tmop{THH} (\mathbb{Z}) / p) [\mu_1^{- 1}]$.
\end{lemma}

\begin{remark}
  \label{rem:Segal-conj-BP-n}More generally, for every $n \in
  \mathbb{Z}_{\geqslant - 1}$, {\cite[Prop~2.7]{AngeliniKnoll2024}} shows
  that, for the $n$-th truncated Brown--Peterson spectrum $\tmop{BP} \langle n
  \rangle$, we have
  \[ \pi_{\ast} (\tmop{THH} (\tmop{BP} \langle n \rangle ; \mathbb{F}_p))
     \cong \mathbb{F}_p [\mu_{n + 1}] \otimes \Lambda (\lambda_1, \ldots,
     \lambda_{n + 1}) \]
  where $\deg \mu_{n + 1} = 2 p^{n + 1}$ and $\deg \lambda_i = 2 p^i - 1$ for
  $i = 1, \ldots, n + 1$. It follows from {\cite[Thm~C]{Hahn2022}} (explained
  in the paragraphs following {\cite[Thm~6.1.4]{Hahn2022a}} as well) that the
  cyclotomic Frobenius $\tmop{THH} (\tmop{BP} \langle n \rangle) \rightarrow
  \tmop{THH} (\tmop{BP} \langle n \rangle)^{t C_p}$ realizes the homotopy
  groups of the target modulo $(p, v_1, \ldots, v_n)$ as the localization of
  the source modulo $(p, v_1, \ldots, v_n)$ at $\mu_{n + 1}$.
\end{remark}

\begin{lemma}[{\cite{Boekstedt1994}}]
  \label{lem:mu-1-equivar-lift}Let $p$ be a prime. Then the class $\mu_1 \in
  \pi_{2 p} (\tmop{THH} (\mathbb{Z}) / p)$ lifts to a $\mathbb{T}$-equivariant
  class $\tilde{\mu}_1$ satisfying the relation $a_{\lambda}  \tilde{\mu}_1 =
  v_1$, where $a_{\lambda}$ is the Euler class in
  {\cite[Ex~6.1.3]{Bhatt2022}}, and $v_1$ is the canonical class in $\pi_{2 p
  - 2} (\tmop{TC}^- (\mathbb{Z}) / p)$\footnote{When $p$ is an odd prime, it
  comes from $v_1 \in \pi_{2 p - 2} (\mathbb{S}/ p)$. When $p = 2$, over the
  sphere, only $v_1^4 \in \pi_8 (\mathbb{S}/ 2)$ is well-defined.}.
\end{lemma}

We give a proof of \Cref{prop:vanish-geom-Cpr-fixed-pts}, since we need
similar techniques to establish our main results.

\begin{lemma}
  \label{lem:bnd-zero-proper-Tate}Let $M$ be a $t$-bounded
  $\mathbb{T}$-equivariant spectrum. Then for every composite number $n \in
  \mathbb{N}_{> 0}$, the proper $C_n$-Tate construction\footnote{See
  {\cite[Def~4.5]{Ayala2021}}.} $M^{\tau C_n}$ is contractible.
\end{lemma}

\begin{proof}
  We may assume that $M$ is concentrated in degree $0$. Then by the lax
  symmetric monoidal structure on $(-)^{\tau C_n}$, we may assume that $M
  =\mathbb{Z}$. Then this is established in {\cite[Prop~5.24]{Mathew2019}}
  (cf.~{\cite[Rem~2.18]{Ayala2017a}}).
\end{proof}

\begin{proposition}[Devalapurkar--Hahn--Raksit--Yuan]
  \label{prop:vanish-geom-Cpr-fixed-pts}The spectrum $(\tmop{THH}
  (\mathbb{Z})^h)^{\Phi C_n}$ is contractible after profinite completion, for
  any composite number $n \in \mathbb{N}_{> 0}$.
\end{proposition}

\begin{proof}
  Equivalently, we have to show that the proper $C_n$-Tate construction
  $\tmop{THH} (\mathbb{Z})^{\tau C_n}$ is contractible after profinite
  completion, for every composite number $n \in \mathbb{N}_{> 0}$. Let $p$ be
  a prime number. Note that $\tmop{THH} (\mathbb{Z})$ is equipped with proper
  Tate diagonal maps
  \[ \tmop{THH} (\mathbb{Z}) \longrightarrow \tmop{THH} (\mathbb{Z})^{\tau
     C_n} \]
  as explained in {\cite[§0.4]{Ayala2017}}, thus $\tmop{THH}
  (\mathbb{Z})^{\tau C_n}$ acquires a $\mathbb{Z}$-module structure, and
  therefore we may check $\tmop{THH} (\mathbb{Z})^{\tau C_n} = 0$ after modulo
  $(p, v_1)$. It follows from \Cref{lem:bnd-zero-proper-Tate} that
  \[ (\tmop{THH} (\mathbb{Z}) / (p, \tilde{\mu}_1))^{\tau C_n} = 0, \]
  and by \Cref{lem:mu-1-equivar-lift}, we are done.
\end{proof}

\begin{corollary}[Devalapurkar--Hahn--Raksit--Yuan]
  Let $M$ be a $\tmop{THH} (\mathbb{Z})^h$-module in $\tmop{Sp}^{g^{<}
  \mathbb{T}}$. Then the genuine $C_p$-fixed points $M^{C_p}$ and the
  geometric $C_p$-fixed points $M^{\Phi C_p}$, viewed as objects in
  $\tmop{Mod}_{\tmop{THH} (\mathbb{Z})^{h C_p}} (\tmop{Sp}^{g^{<} 
  (\mathbb{T}/ C_p)})$, are Borel complete.
\end{corollary}

\Cref{prop:vanish-geom-Cpr-fixed-pts} allows us to reduce to analyze
de-completed $C_p$-Tate construction, by joint conservativity of $\{ (-)^{\Phi
C_n} \barsuchthat n \in \mathbb{N}_{> 0} \}$ and
\Cref{rem:decpl-Tate-Borel-cpl}:

\begin{corollary}
  \label{cor:asm-de-cpl-Borel-Cp-Tate}Let $M$ be a $\tmop{THH}
  (\mathbb{Z})^h$-module in $\tmop{Sp}^{g^{<} \mathbb{T}}$. Then the assembly
  map
  \[ M^{\eta_{\tmop{THH} (\mathbb{Z})^h}} \longrightarrow M^h \]
  is an equivalence after $p$-completion if and only if the assembly map
  \[ M^{\eta_{\tmop{THH} (\mathbb{Z})^h} C_p} \longrightarrow M^{t C_p} \]
  is an equivalence after $p$-completion.
\end{corollary}

\begin{proof*}{Proof of \Cref{thm:THH-prism-comp}}
  It follows from
  \Cref{cor:asm-de-cpl-Borel-Cp-Tate,lem:asm-polygonic-THH-eq}.
\end{proof*}

It remains to establish the following version of
\Cref{lem:asm-polygonic-HH-eq} for $\tmop{THH}$:

\begin{lemma}
  \label{lem:asm-polygonic-THH-eq}Let $R$ be an
  $\mathbb{E}_1$-$\mathbb{Z}$-algebra with bounded $\tmop{Tor}$-amplitude in
  $D (\mathbb{Z})$, and $M$ a perfect $R$-$R$-bimodule in $D (\mathbb{Z})$.
  Then the assembly map
  \[ \tmop{THH} (R ; M^{\otimes_R^{\mathbb{L}} C_p})^{\theta_{\tmop{THH}
     (\mathbb{Z})} C_p} \longrightarrow \tmop{THH} (R ;
     M^{\otimes_R^{\mathbb{L}} C_p})^{t C_p} \]
  is an equivalence after $p$-completion.
\end{lemma}

\begin{proof}
  We mimic the proof of \Cref{lem:asm-polygonic-HH-eq}. By the exactness of
  both sides in $M$, we may assume that $M$ is the free $R$-$R$-bimodule $R
  \otimes R$ of rank $1$. Then the map in question becomes the assembly map
  \[ \tmop{THH} (\mathbb{Z}; R^{\otimes_{\mathbb{Z}}^{\mathbb{L}}
     p})^{\theta_{\tmop{THH} (\mathbb{Z})} C_p} \longrightarrow \tmop{THH}
     (\mathbb{Z}; R^{\otimes_{\mathbb{Z}}^{\mathbb{L}} C_p})^{t C_p} \]
  which follows from \Cref{lem:polygon-THH-Z-decpl-Tate}, being a polygonic
  $\tmop{THH}$ version of \Cref{lem:p-th-pow-decpl-Tate}.
\end{proof}

\begin{lemma}
  \label{lem:polygon-THH-Z-decpl-Tate}Let $M$ be a $\mathbb{Z}$-module
  spectrum of bounded $\tmop{Tor}$-amplitude. Then the assembly map
  \[ \tmop{THH} (\mathbb{Z}; M^{\otimes_{\mathbb{Z}}^{\mathbb{L}}
     p})^{\theta_{\tmop{THH} (\mathbb{Z})} C_p} \longrightarrow \tmop{THH}
     (\mathbb{Z}; M^{\otimes_{\mathbb{Z}}^{\mathbb{L}} p})^{t C_p} \]
  is an equivalence after $p$-completion.
\end{lemma}

\begin{proof}
  Again, since both sides are exact in $M$, we may assume that $M$ is a flat
  $\mathbb{Z}$-module, and the result is true for finite free
  $\mathbb{Z}$-modules $M$. By Lazard's theorem, the flat $\mathbb{Z}$-module
  $M$ is a filtered colimit of finite free modules, thus it suffices to show
  that the functor
  \begin{eqnarray*}
    \tmop{Mod}_{\mathbb{Z}}^{\flat} & \longrightarrow & \tmop{Sp}_p^{\wedge}\\
    M & \longmapsto & (\tmop{THH} (\mathbb{Z};
    M^{\otimes_{\mathbb{Z}}^{\mathbb{L}} p})^{t C_p})_p^{\wedge}
  \end{eqnarray*}
  preserves filtered colimits, where $\tmop{Mod}_{\mathbb{Z}}^{\flat}$ is the
  category of flat $\mathbb{Z}$-modules. Note that $\tmop{THH} (\mathbb{Z};
  M^{\otimes_{\mathbb{Z}}^{\mathbb{L}} p})^{t C_p}$ acquires a $\tmop{THH}
  (\mathbb{Z})^{t C_p}$-module structure, and in particular, a
  $\mathbb{Z}$-module structure, thus $v_1 = 0$ on it. Now by
  \Cref{lem:mu-1-equivar-lift}, we have
  \[ \tmop{THH} (\mathbb{Z}; M^{\otimes_{\mathbb{Z}}^{\mathbb{L}} p})^{t C_p}
     / (p, v_1) \simeq (\tmop{THH} (\mathbb{Z}; M^{\otimes_p^{\mathbb{L}} p})
     / (p, \tilde{\mu}_1))^{t C_p} \]
  and as (non-equivariant) spectra, we have
  \[ \tmop{THH} (\mathbb{Z}; M^{\otimes_{\mathbb{Z}}^{\mathbb{L}} p}) / (p,
     \tilde{\mu}_1) \simeq (\tmop{THH} (\mathbb{Z}) / (p, \widetilde{\mu_1}))
     \otimes_{\mathbb{Z}}^{\mathbb{L}} M^{\otimes_{\mathbb{Z}}^{\mathbb{L}} p}
  \]
  which is bounded, and preserves filtered colimits in $M$. The result the
  follows from the fact that $(-)^{t C_p}$ preserves uniformly bounded
  filtered colimits.
\end{proof}

Finally, we relate de-completed homotopy fixed points and de-completed Tate
construction to de-completed Borel completion when the base is $\tmop{THH}
(\mathbb{Z})^h$. We would thank Georg {\tmname{Tamme}} for such a question.
The point is that, as in the proofs of
\Cref{prop:vanish-geom-Cpr-fixed-pts,lem:asm-polygonic-THH-eq}, after $(p,
v_1)$-completion, the base $\tmop{THH} (\mathbb{Z})^h$ is close to be
$t$-bounded. We start with some simple lemmas.

\begin{lemma}
  \label{lem:zero-ring-mods}Let $V \in \tmop{Sp}^{\omega}$ be a finite
  spectrum, and $A$ an $\mathbb{E}_1$-ring such that $A \otimes V = 0$. Then
  for every right $A$-module $M$, we have $M \otimes V = 0$.
\end{lemma}

\begin{proof}
  Note that $M$ is a retract of $M \otimes A$ as a spectrum, thus $M \otimes
  V$ is a retract of $M \otimes A \otimes V = 0$ as a spectrum.
\end{proof}

\begin{remark}
  In our applications, we mainly take $V$ to be $\mathbb{S}/ (p, v_1^4)$ (or
  equivalently, $\mathbb{S}/ (p, v_1)$ when $p = 2$).
\end{remark}

\begin{lemma}
  \label{lem:htpy-fixed-pts-pres-colim}Let $V \in \tmop{Sp}^{\omega}$ be a
  finite spectrum, and $A_0 \rightarrow A$ a map of $\mathbb{E}_1$-algebras in
  $\tmop{Sp}^{B\mathbb{T}}$ such that $A_0 \otimes V$ is $t$-bounded below.
  Suppose that the spectrum $A_0^{t C_p} \otimes V$, which is in fact
  equivalent to $(A_0 \otimes V)^{t C_p}$, is contractible. Then
  \begin{enumerate}
    \item for every right $A$-module $M$, the spectrum
    $(M^{t\mathbb{T}})_p^{\wedge} \otimes V$ is contractible; and
    
    \item the functor
    \[ ((-)^{h\mathbb{T}})_p^{\wedge} \otimes V \of \tmop{RMod}_A
       (\tmop{Sp}^{B\mathbb{T}}) \longrightarrow \tmop{Sp}_p^{\wedge} \]
    preserves filtered colimits.
  \end{enumerate}
\end{lemma}

\begin{proof}
  For the contractibility of $(M^{t\mathbb{T}})_p^{\wedge} \otimes V$, or
  equivalently, that of $M^{t\mathbb{T}} \otimes (V / p)$, by
  \Cref{lem:zero-ring-mods} and the lax symmetric monoidal structure on
  $(-)^{t\mathbb{T}}$, it suffices to show the contractibility of
  $(A_0^{t\mathbb{T}})_p^{\wedge} \otimes V$. But then
  \begin{eqnarray*}
    (A_0^{t\mathbb{T}})_p^{\wedge} \otimes V & \simeq & ((A_0 \otimes
    V)^{t\mathbb{T}})_p^{\wedge}\\
    & \simeq & ((A_0 \otimes V)^{t C_p})^{h (\mathbb{T}/ C_p)}\\
    & \simeq & 0,
  \end{eqnarray*}
  where the second equivalence follows from the $t$-bounded-belowness of $A_0
  \otimes V$ and {\cite[Lem~II.4.2]{Nikolaus2018}}.
  
  The filtered-colimit-preservation of the functor
  $((-)^{h\mathbb{T}})_p^{\wedge} \otimes V$ follows directly from the
  canonical fiber sequence
  \[ \Sigma (-)_{h\mathbb{T}} \xrightarrow{\tmop{Nm}_{\mathbb{T}}}
     (-)^{h\mathbb{T}} \longrightarrow (-)^{t\mathbb{T}} \]
  and the Tate vanishing established above.
\end{proof}

\begin{proposition}
  \label{prop:de-cpl-Borel-Tate-t-bnded}Let $A_0 \rightarrow A$ be a map of
  $\mathbb{E}_1$-algebras in $\tmop{Sp}^{B\mathbb{T}}$, and $V \in
  \tmop{Sp}^{\omega}$ a finite spectrum such that $A_0 \otimes V$ is
  $t$-bounded. Then
  \begin{enumerate}
    \item every $A^{C_p}$-module in $\tmop{Sp}^{g^{<}  (\mathbb{T}/ C_p)}$ is
    $p$-completely\footnote{Equivalently, it is saying that every
    $A^{C_p}$-module in $\tmop{Sp}^{g^{<} (\mathbb{T}/ C_p)}$ is
    Borel-complete up to Bousfield $V / p$-localization.} Borel-complete up to
    Bousfield $V$-localization; and
    
    \item the $p$-completed de-completed homotopy fixed points functor
    \[ ((-)^{\eta_A \mathbb{T}})_p^{\wedge} \of \tmop{RMod}_A
       (\tmop{Sp}^{g^{<} \mathbb{T}}) \longrightarrow \tmop{Sp}_p^{\wedge} \]
    coincides with the $p$-completed genuine $C_{p^{\infty}}$-fixed points
    functor up to Bousfield $V$-localization, and the $p$-completed
    de-completed Tate construction functor
    \[ ((-)^{\theta_A \mathbb{T}})_p^{\wedge} \of \tmop{RMod}_A
       (\tmop{Sp}^{g^{<} \mathbb{T}}) \longrightarrow \tmop{Sp}_p^{\wedge} \]
    coincides with the composite functor
    \begin{equation}
      \tmop{RMod}_A (\tmop{Sp}^{g^{<} \mathbb{T}}) \xrightarrow{(-)^{\Phi
      C_p}} \tmop{RMod}_{A^{\Phi C_p}} (\tmop{Sp}^{g^{<}  (\mathbb{T}/ C_p)})
      \xrightarrow{(-)^{C_{p^{\infty}} / C_p}} \tmop{Sp}_p^{\wedge}
    \end{equation}
    up to Bousfield $V$-localization.
  \end{enumerate}
\end{proposition}

\begin{proof}
  We start with the first statement. Firstly, we check that, for every $r \geq
  1$, we have $(A_0^{C_p})^{\Phi C_{p^r}} \otimes (V /^{\mathbb{L}} p) \simeq
  0$. By ``genuine'' Tate-orbit lemma {\cite[Lem~5.3]{Antieau2018}}, the
  canonical map
  \[ (A_0^{C_p})^{\Phi C_{p^r}} \otimes (V / p) = ((A_0 \otimes V /
     p)^{C_p})^{\Phi C_{p^r}} \longrightarrow ((A_0 \otimes V / p)^{\Phi
     C_p})^{\Phi C_{p^r}} = (A_0 \otimes V / p)^{\Phi C_{p^{r + 1}}} \]
  is an equivalence. By \Cref{lem:bnd-zero-proper-Tate}, we have $(A_0 \otimes
  V / p)^{\Phi C_{p^{r + 1}}} \simeq 0$, and by \Cref{lem:zero-ring-mods}, we
  see that $(A_0 \otimes V / p)^{\tau C_{p^{r + 1}}} \simeq 0$.
  
  It follows that $A_0^{C_p}$ is $p$-completely Borel-complete up to Bousfield
  $V$-localization, as by $p$-completeness, it suffices to check the
  isomorphisms $(-)^{\Phi H} \rightarrow (-)^{\tau H}$ on cyclic $p$-groups
  $H$. Then it follows from \Cref{lem:zero-ring-mods} and the lax symmetric
  monoidal structure on $(-)^{C_m}$ and on the map $(-)^{\Phi C_m} \rightarrow
  (-)^{\tau C_m}$ that
  \begin{enumerate}
    \item $(A^{C_p})^{\Phi C_m} \otimes V / p = 0$, and thus
    
    \item for every right $A^{C_p}$-module $M$ in $\tmop{Sp}^{g^{<} 
    (\mathbb{T}/ C_p)}$, we have $M^{\tau C_m} \otimes V / p = M^{\Phi C_m}
    \otimes V / p = 0.$
  \end{enumerate}
  In particular, every right $A^{C_p}$-module in $\tmop{Sp}^{g^{<} 
  (\mathbb{T}/ C_p)}$ is $p$-completely Borel-complete up to Bousfield
  $V$-localization.
  
  Now we prove the second statement. By construction, on compact objects, we
  have the desired results, essentially by {\cite[Prop~5.9]{Ayala2021}} and
  \Cref{prop:Tate-orb-bnd-mods,prop:lim-formula-bnd-mods-tT}. It suffices to
  show that the functors $(-)^{C_{p^{\infty}}}$ and $((-)^{\Phi
  C_p})^{C_{p^{\infty}} / C_p}$ preserve filtered colimits.
  
  Note that $(-)^{C_{p^{\infty}}} = ((-)^{C_p})^{C_{p^{\infty}} / C_p}$, that
  is to say, we may write both functors $(-)^{C_{p^{\infty}}}$ and $((-)^{\Phi
  C_p})^{C_{p^{\infty}} / C_p}$ as a composite
  \[ \tmop{RMod}_A (\tmop{Sp}^{g^{<} \mathbb{T}}) \longrightarrow
     \tmop{RMod}_{A^{C_p}} (\tmop{Sp}^{g^{<}  (\mathbb{T}/ C_p)})
     \xrightarrow{(-)^{C_{p^{\infty}} / C_p}} \tmop{Sp}_p^{\wedge}, \]
  where we are tacitly using the lax symmetric monoidal structure on the
  natural transformation $(-)^{C_p} \Rightarrow (-)^{\Phi C_p}$. Note that the
  first functor $(-)^{\Phi C_p}$ or $(-)^{C_p}$ preserves filtered colimits,
  thus it suffices to see that the second functor $(-)^{C_{p^{\infty}} / C_p}$
  preserves filtered colimits up to Bousfield $V$-localization.
  
  We are now reduced to show that the functor
  \[ \tmop{RMod}_{A^{h C_p}} (\tmop{Sp}^{B (\mathbb{T}/ C_p)})
     \xrightarrow{(-)^{h (\mathbb{T}/ C_p)} \otimes V} \tmop{Sp}_p^{\wedge} \]
  preserves filtered colimits. By examining the fiber sequence
  \[ \Sigma (-)_{h (\mathbb{T}/ C_p)} \xrightarrow{\tmop{Nm}_{\mathbb{T}/
     C_p}} (-)^{h (\mathbb{T}/ C_p)} \longrightarrow (-)^{t (\mathbb{T}/ C_p)}
  \]
  It remains to show that, for every right $A^{h C_p}$-module $M$ in
  $\tmop{Sp}^{B (\mathbb{T}/ C_p)}$, we have $(M^{t (\mathbb{T}/
  C_p)})_p^{\wedge} \otimes V = 0$. Actually, for every right $A_0^{h
  C_p}$-module $M$, the lax symmetric monoidal structure on $(-)^{t
  (\mathbb{T}/ C_p)}$ gives rise to a $((A_0^{h C_p})^{t (\mathbb{T}/
  C_p)})_p^{\wedge}$-module structure on $(M^{t (\mathbb{T}/
  C_p)})_p^{\wedge}$, and since $A_0 \otimes V$ is $t$-bounded, by
  \Cref{lem:S1-Tate-fixed-pt}, we have
  \begin{eqnarray*}
    ((A_0^{h C_p})^{t (\mathbb{T}/ C_p)})_p^{\wedge} \otimes V & \simeq &
    (((A_0 \otimes V)^{h C_p})^{t (\mathbb{T}/ C_p)})_p^{\wedge}\\
    & \simeq & 0
  \end{eqnarray*}
  and thus $(M^{t (\mathbb{T}/ C_p)})_p^{\wedge} \otimes V = 0$ by
  \Cref{lem:zero-ring-mods}.
\end{proof}

We need the following consequence of the Tate fixed point lemma.

\begin{lemma}
  \label{lem:S1-Tate-fixed-pt}Let $M$ be a $t$-bounded
  $\mathbb{T}$-equivariant spectrum. Then $((M^{h C_p})^{t (\mathbb{T}/
  C_p)})_p^{\wedge} = 0$.
\end{lemma}

\begin{proof}
  By dévissage, we reduce to the case that $M$ is concentrated in degree $0$.
  In this case, we have
  \begin{eqnarray*}
    ((M^{h C_p})^{t (\mathbb{T}/ C_p)}) / p & \simeq & (M^{h C_p})^{t (C_{p^2}
    / C_p)}\\
    & \simeq & 0,
  \end{eqnarray*}
  by the Tate fixed point lemma and {\cite[Lem~IV.4.12]{Nikolaus2018}}.
\end{proof}

\begin{example}
  Let $k$ be an animated ring, and $\mathcal{C}$ a presentable stable
  $k$-linear $\infty$-category. Then by \Cref{prop:de-cpl-Borel-Tate-t-bnded},
  we have
  \begin{eqnarray*}
    \tmop{HC}^{-, \tmop{poly}} (\mathcal{C}/ k)_p^{\wedge} & = & \left(
    (\tmop{HH} (\mathcal{C}/
    \underline{k})^{\eta_{\underline{k}}})^{C_{p^{\infty}}}
    \right)_p^{\wedge},\\
    \tmop{HP}^{\tmop{poly}} (\mathcal{C}/ k)_p^{\wedge} & = & \left(
    ((\tmop{HH} (\mathcal{C}/ \underline{k})^{\eta_{\underline{k}}})^{\Phi
    C_p})^{h (C_{p^{\infty}} / C_p)} \right)_p^{\wedge}\\
    & = & ((\tmop{THH} (\mathcal{C}) \otimes_{\tmop{THH} (k)}^{\mathbb{L}}
    k^{t C_p})^{h (\mathbb{T}/ C_p)})_p^{\wedge} .
  \end{eqnarray*}
  In retrospect, the object $(\tmop{THH} (\mathcal{C}) \otimes_{\tmop{THH}
  (k)}^{\mathbb{L}} k^{t C_p})^{h (\mathbb{T}/ C_p)}$ appeared
  {\cite[Rem~3.5.6]{Devalapurkar2023}} in the even more general setup of $k$
  being an $\mathbb{E}_{\infty}$-ring.
\end{example}

Now we deduce the comparison over $\tmop{THH} (R)^h$ for every
$\mathbb{E}_{\infty}$-$\mathbb{Z}$-algebra $R$.

\begin{corollary}
  \label{cor:de-cpl-Borel-Tate-THH-Z}Let $A$ be an $\mathbb{E}_1$-$\tmop{THH}
  (\mathbb{Z})$-algebra in $\tmop{Sp}^{B\mathbb{T}}$. Then the $(p,
  v_1)$-completed de-completed homotopy fixed points functor
  \[ ((-)^{\eta_A \mathbb{T}})_{(p, v_1)}^{\wedge} \of \tmop{RMod}_A
     (\tmop{Sp}^{g^{<} \mathbb{T}}) \longrightarrow \tmop{Sp}_{(p,
     v_1)}^{\wedge} \]
  coincides with the $(p, v_1)$-completed genuine $C_{p^{\infty}}$-fixed
  points functor, and the $(p, v_1)$-completed de-completed Tate construction
  functor
  \[ ((-)^{\theta_A \mathbb{T}})_{(p, v_1)}^{\wedge} \of \tmop{Mod}_A
     (\tmop{Sp}^{g^{<} \mathbb{T}}) \longrightarrow \tmop{Sp}_{(p,
     v_1)}^{\wedge} \]
  coincides with the composite functor
  \[ \tmop{Mod}_A (\tmop{Sp}^{g^{<} \mathbb{T}}) \xrightarrow{(-)^{\Phi C_p}}
     \tmop{Mod}_{A^{\Phi C_p}} (\tmop{Sp}^{g^{<}  (\mathbb{T}/ C_p)})
     \xrightarrow{(-)^{C_{p^{\infty}} / C_p}} \tmop{Sp}_{(p, v_1)}^{\wedge} .
  \]
\end{corollary}

\begin{proof}
  As a very special case of {\cite[Thm~G]{Hahn2022}}, the spectrum $\tmop{TR}
  (\mathbb{Z}) / (p, v_1)$ is $t$-bounded. The result then follows from
  \Cref{prop:de-cpl-Borel-Tate-t-bnded} by taking $V =\mathbb{S}/ (p, v_1^4)$
  (or $\mathbb{S}/ (p, v_1)$ when $p > 2$), and $A_0 = \tmop{TR}
  (\mathbb{Z})^h$.
\end{proof}

\begin{remark}
  In \Cref{cor:de-cpl-Borel-Tate-THH-Z}, we can replace $\mathbb{Z}$ by
  $\tmop{BP} \langle n \rangle$, and $(p^{\alpha_0}, v_1^{\alpha_1}, \ldots,
  v_{n + 1}^{\alpha_{n + 1}})$ in place of $(p, v_1^4)$ for some $(\alpha_0,
  \ldots, \alpha_{n + 1}) \in \mathbb{N}_{> 0}^{n + 1}$, for which we need the
  full {\cite[Thm~G]{Hahn2022}}.
\end{remark}

\begin{remark}
  \label{rem:mu1-THH-perfd}Let $S$ be a perfectoid ring. Then the image of
  $\mu_1 \in \pi_{2 p} (\tmop{THH} (\mathbb{Z}) / p)$ under the map
  $\tmop{THH} (\mathbb{Z}) \rightarrow \tmop{THH} (S)$ is $\sigma^p$ up to a
  multiplier of a unit, where $\sigma \in \pi_2 (\tmop{THH} (S)_p^{\wedge})$
  is the Bökstedt element. To see this, when $S =\mathbb{F}_p$, this is
  established in {\cite{Boekstedt1994}}, which implies the case for $S$ being
  a perfect $\mathbb{F}_p$-algebra, since the map $\tmop{THH} (\mathbb{Z})
  \rightarrow \tmop{THH} (S)$ factors through $\tmop{THH} (\mathbb{F}_p)$. The
  general case follows follows by considering the composite
  \[ \tmop{THH} (\mathbb{Z}) \longrightarrow \tmop{THH} (S) \longrightarrow
     \tmop{THH} \left( S / \sqrt{pS} \right) \]
  and note that $S / \sqrt{pS}$ is a perfect $\mathbb{F}_p$-algebra.
\end{remark}

\begin{example}
  \label{ex:TP-poly-de-cpl-Borel-cpl}Let $S$ be a perfectoid ring, and
  $\mathcal{C}$ a dualizable presentable stable $S$-linear $\infty$-category.
  Then we have
  \begin{eqnarray*}
    \tmop{TC}^{-, \tmop{poly} / S} (\mathcal{C}/ S)_{(p, \xi)}^{\wedge} & = &
    \left( (\tmop{THH} (\mathcal{C})^{\eta_{\tmop{THH} (S)}})^{C_{p^{\infty}}}
    \right)_{(p, \xi)}^{\wedge},\\
    \tmop{TP}^{\tmop{poly} / S} (\mathcal{C}/ S)_{(p, \xi)}^{\wedge} & = &
    \left( ((\tmop{THH} (\mathcal{C})^{\eta_{\tmop{THH} (S)}})^{\Phi
    C_p})^{C_{p^{\infty}} / C_p} \right)_{(p, \xi)}^{\wedge} .
  \end{eqnarray*}
  Indeed, they follow from
  \Cref{cor:de-cpl-Borel-Tate-THH-Z,rem:mu1-THH-perfd} by taking $A \assign
  \tmop{THH} (S)$ (where $a_{\lambda}  \tilde{\sigma} = \xi$ for the generator
  $\tilde{\sigma} \in \pi_2 (\tmop{TC}^- (S)_p^{\wedge})$ up to a complex
  orientation, and $a_{\lambda}  \tilde{\mu}_1 = v_1$).
\end{example}

Finally, we show that de-completed Borel completion over $\tmop{THH} (S)$ for
perfectoid rings $S$, and over $\mathbb{Z}$, are ``the same'' as over
$\tmop{THH} (\mathbb{Z})$.

\begin{proposition}
  \label{prop:asm-THH-S-over-THH-Z}The assembly map $\tmop{THH}
  (S)^{\eta_{\tmop{THH} (\mathbb{Z})}} \rightarrow \tmop{THH} (S)^h$ in
  $\tmop{Mod}_{\tmop{THH} (\mathbb{Z})^h} (\tmop{Sp}^{g^{<} \mathbb{T}})$ is
  an equivalence after $p$-completion.
\end{proposition}

\begin{proof}
  By
  \Cref{cor:asm-de-cpl-Borel-Cp-Tate,lem:decpl-Cp-Tate,lem:de-cpl-Borel-cpl-base-chg},
  it is equivalent to show that the map
  \[ \tmop{THH} (S) \otimes_{\tmop{THH} (\mathbb{Z})}^{\mathbb{L}} \tmop{THH}
     (\mathbb{Z})^{t C_p} \longrightarrow \tmop{THH} (S)^{t C_p} \]
  induced by the cyclotomic Frobenius $\tmop{THH} (S) \rightarrow \tmop{THH}
  (S)^{t C_p}$ is an equivalence after modulo $p$, but this follows
  immediately from \Cref{rem:mu1-THH-perfd}.
\end{proof}

It then follows from \Cref{lem:de-cpl-Borel-cpl-base-chg} that

\begin{corollary}
  \label{cor:de-cpl-THH-S-THH-Z}Let $M$ be a $\tmop{THH} (S)$-module in
  $\tmop{Sp}^{g^{<} \mathbb{T}}$. Then the assembly map
  \[ M^{\eta_{\tmop{THH} (\mathbb{Z})}} \longrightarrow M^{\eta_{\tmop{THH}
     (S)}} \]
  in $\tmop{Mod}_{\tmop{THH} (\mathbb{Z})^h} (\tmop{Sp}^{g^{<} \mathbb{T}})$
  is an equivalence after $p$-completion.
\end{corollary}

\begin{remark}
  \label{rem:THH-Prism-BP-n}More generally, we may replace $\mathbb{Z}$ by any
  truncated Brown--Peterson spectrum. Indeed, for $- 1 \leq n < m$, the
  assembly map
  \[ \tmop{THH} (\tmop{BP} \langle n \rangle)^{\eta_{\tmop{THH} (\tmop{BP}
     \langle m \rangle)}} \longrightarrow \tmop{THH} (\tmop{BP} \langle n
     \rangle)^h \]
  in $\tmop{RMod}_{\tmop{THH} (\tmop{BP} \langle m \rangle)^h}
  (\tmop{Sp}^{g^{<} \mathbb{T}})$ is an equivalence after $(p, v_1, \ldots,
  v_n)$-completion, as a consequence of \Cref{rem:Segal-conj-BP-n}. In
  particular, this shows that $\tmop{THH}_{\Prism}$ extends to a localizing
  invariant
  \[ \tmop{Cat}_{\tmop{BP} \langle n \rangle}^{\tmop{perf}}
     \xrightarrow{\tmop{THH}} \tmop{RMod}_{\tmop{THH} (\tmop{BP} \langle n
     \rangle)} (\tmop{Sp}^{g^{<} \mathbb{T}})
     \xrightarrow{(-)^{\eta_{\tmop{THH} (\tmop{BP} \langle n \rangle)}}}
     \tmop{RMod}_{\tmop{THH} (\tmop{BP} \langle n \rangle)^h}
     (\tmop{Sp}^{g^{<} \mathbb{T}})_{(p, v_1, \ldots, v_n)}^{\wedge} \]
  which is independent of choice of $n \in \mathbb{Z}_{\geqslant - 1}$ in the
  sense of compatibility with forgetful functors.
\end{remark}

On the other hand, we can also deduce the result for $\mathbb{Z}$ from the
case $S =\mathbb{F}_p$:

\begin{proposition}
  The assembly map $\underline{\mathbb{Z}}^{\eta_{\tmop{THH} (\mathbb{Z})}}
  \longrightarrow \mathbb{Z}$ in $\tmop{Mod}_{\tmop{THH} (\mathbb{Z})^h}
  (\tmop{Sp}^{g^{<} \mathbb{T}})$ is an equivalence.
\end{proposition}

\begin{proof}
  Since $\underline{\mathbb{Q}} \in \tmop{CAlg} (\tmop{Sp}^{g^{<}
  \mathbb{T}})$ is Borel, it suffices to show this equivalence after modulo
  every prime $p$. By
  \Cref{cor:asm-de-cpl-Borel-Cp-Tate,lem:decpl-Cp-Tate,lem:de-cpl-Borel-cpl-base-chg},
  we are reduced to show that the map
  \[ \tau_{\geqslant 0} (\mathbb{F}_p^{t C_p}) \otimes_{\tmop{THH}
     (\mathbb{Z})}^{\mathbb{L}} \tmop{THH} (\mathbb{Z})^{t C_p}
     \longrightarrow \mathbb{F}_p^{t C_p} \]
  is an equivalence, where the $\mathbb{E}_{\infty}$-map $\tmop{THH}
  (\mathbb{Z}) \rightarrow \tau_{\geqslant 0} (\mathbb{F}_p^{t C_p})$ is the
  composite
  \[ \tmop{THH} (\mathbb{Z}) \xrightarrow{\simeq} \tau_{\geqslant 0}
     (\tmop{THH} (\mathbb{Z})^{t C_p}) \longrightarrow \tau_{\geqslant 0}
     (\mathbb{Z}^{t C_p}) \longrightarrow \tau_{\geqslant 0} (\mathbb{F}_p^{t
     C_p}) \]
  where the first equivalence is the Segal conjecture for $\tmop{THH}
  (\mathbb{Z})$. Now this composite can be rewritten as the composite
  \[ \tmop{THH} (\mathbb{Z}) \xrightarrow{\simeq} \tau_{\geqslant 0}
     (\tmop{THH} (\mathbb{Z})^{t C_p}) \longrightarrow \tau_{\geqslant 0}
     (\tmop{THH} (\mathbb{F}_p)^{t C_p}) \longrightarrow \tau_{\geqslant 0}
     (\mathbb{F}_p^{t C_p}), \]
  and we have a composite pushout diagram
  \[ \begin{array}{ccccc}
       \tau_{\geqslant 0} (\tmop{THH} (\mathbb{Z})^{t C_p}) & \longrightarrow
       & \tau_{\geqslant 0} (\tmop{THH} (\mathbb{F}_p)^{t C_p}) &
       \longrightarrow & \tau_{\geqslant 0} (\mathbb{F}_p^{t C_p})\\
       \longdownarrow &  & \longdownarrow &  & \longdownarrow\\
       \tmop{THH} (\mathbb{Z})^{t C_p} & \longrightarrow & \tmop{THH}
       (\mathbb{F}_p)^{t C_p} & \longrightarrow & \mathbb{F}_p^{t C_p}
     \end{array} \]
  of $\mathbb{E}_{\infty}$-rings, where the left pushout is
  \Cref{prop:asm-THH-S-over-THH-Z}, and the right pushout is simply detected
  by both vertical maps are inverting the Bökstedt element $\mu_0 \in \pi_2
  (\tmop{THH} (\mathbb{F}_p)^{t C_p})$.
\end{proof}

\begin{remark}
  \label{rem:THH-Z-Cp-Tate-THH-Fp}Let $p$ be an odd prime.
  {\cite[Prop~5.2.1]{Devalapurkar2025}} implies that the map $\tmop{THH}
  (\mathbb{Z})^{t C_p} \rightarrow \mathbb{Z}^{t C_p}$ coincides with the map
  $\tmop{THH} (\mathbb{Z})^{t C_p} \rightarrow \tmop{THH} (\mathbb{F}_p)^{t
  C_p}$ induced by the ring map $\mathbb{Z} \rightarrow \mathbb{F}_p$.
\end{remark}

\begin{question}
  Is the statement of \Cref{rem:THH-Z-Cp-Tate-THH-Fp} true when $p = 2$?
\end{question}

\begin{corollary}
  \label{cor:de-cpl-Z-THH-Z}Let $M$ be a $\underline{\mathbb{Z}}$-module in
  $\tmop{Sp}^{g^{<} \mathbb{T}}$. Then the assembly map
  \[ M^{\eta_{\tmop{THH} (\mathbb{Z})}} \longrightarrow
     M^{\eta_{\underline{\mathbb{Z}}}} \]
  in $\tmop{Mod}_{\tmop{THH} (\mathbb{Z})^h} (\tmop{Sp}^{g^{<} \mathbb{T}})$
  is an equivalence.
\end{corollary}

\section{Noncommutative crystalline--de Rham comparison revisited}

In this section, we explain that original ideas in {\cite{Petrov2019}}, along
with an observation by A.~{\tmname{Raksit}}, adapt to de-completed case, which
leads to a comparison (\Cref{cor:Z-lin-cat-crys-comp}), which
``authentically'' corresponds to the crystalline--de Rham comparison. As
explained in \Cref{rem:crys-dR-comp-odd-primes}, the argument in
{\cite{Devalapurkar2025}} is enough for this comparison when $p$ is an odd
prime, so our result is stronger only when $p = 2$.

The key is to produce a map $\tmop{HH} (\mathbb{F}_p
/\mathbb{Z})^{\theta_{\mathbb{Z}} C_p} \rightarrow \mathbb{Z}^{t C_p}$
(Construction~\ref{cons:map-HH-PD-struct}). As A.~{\tmname{Raksit}} observed,
up to $(-)^{t C_p}$, we can replace $\tmop{HH} (\mathbb{Z} [t] /\mathbb{Z})$
by $\tmop{dR}_{\mathbb{Z} [t] /\mathbb{Z}}$, and consequently, we can identify
$\tmop{HH} (\mathbb{F}_p /\mathbb{Z})^{\theta_{\mathbb{Z}} C_p}$ with
$\tmop{dR}_{\mathbb{F}_p /\mathbb{Z}}^{t C_p}$. Combining with the map
$\tmop{dR}_{\mathbb{F}_p /\mathbb{Z}} \rightarrow \mathbb{Z}_{(p)}$ induced by
the PD-structure on $(p) \subseteq \mathbb{Z}_{(p)}$, we get the map we want.
As a consequence, we get a map $\tmop{THH} (\mathbb{F}_p)^{t C_p} \rightarrow
\tmop{HH} (\mathbb{F}_p /\mathbb{Z})^{\theta_{\mathbb{Z}} C_p} \rightarrow
\mathbb{Z}^{t C_p}$. We then show that this is an equivalence
(\Cref{prop:THH-Fp-tCp-Z-tCp-eq}), the blueprint of whose proof can be found
after Construction~\ref{cons:THH-Fp-Z-tCp}.

We start with A.~{\tmname{Raksit}}'s observation. Let $k$ be a commutative
ring, and $R$ a smooth commutative $k$-algebra. In general, we are only
equipped with an HKR-filtration on the periodic cyclic homology whose
associated graded pieces are equivalent to shifts of Hodge-completed derived
de Rham cohomology of $R / k$. However, when $R$ is of relative
dimension$\nosymbol \leq 1$, the filtration in question, along with the
commutative ring structure\footnote{The key is the commutative ring
structure.}, is simply extended from the filtration on $k^{t\mathbb{T}}$,
which was observed by A.~{\tmname{Raksit}} as in
Construction~\ref{cons:map-dR-HH}. A slight improvement of his argument gives
a stronger result: we compare associated $[0, 1]$-graded pieces of
HKR-filtered Hochschild homology and Hodge-filtered de Rham cohomology in
\Cref{cor:HKR-dR-shear}. For this purpose, we first review some facts on
homotopy coherent cochain complexes with weights in $[0, 1]$, which is
systematically studied in K.~{\tmname{Magidson}}'s work {\cite{Magidson2024}}.

\begin{notation}
  Let $k$ be a commutative ring. We denote by $\tmop{DG}_- \tmop{DAlg}_k$ the
  $\infty$-category of $h_-$-differential graded derived commutative
  $k$-algebras {\cite[Nota~5.3.1]{Raksit2020}}, and $\tmop{DG}_-^{[0, 1]}
  \tmop{DAlg}_k \subseteq \tmop{DG}_- \tmop{DAlg}_k$ the full subcategory of
  objects whose weights are concentrated in degrees $[0, 1]$.
\end{notation}

\begin{remark}[{\cite[§3.3]{Magidson2024}}]
  \label{rem:DG-alg-sq-zero-ext}Let $k$ be a commutative ring. Then the
  $\infty$-category $\tmop{DG}_-^{[0, 1]} \tmop{DAlg}_k$ is equivalent to the
  $\infty$-category of triples $(A, M, \eta)$ where
  \begin{itemize}
    \item $A$ is a derived commutative $k$-algebra;
    
    \item $M$ is an $A$-module; and
    
    \item $\eta \of A \rightarrow M$ is a derivation (or equivalently, a map
    $L_{A / k} \rightarrow M$ of $A$-modules).
  \end{itemize}
  Given a triple $(A, M, \eta)$ as above,the corresponding object in
  $\tmop{DG}_-^{[0, 1]} \tmop{DAlg}_k$ is equivalent to the square-zero
  extension of $A (0)$ by $M (1) [- 1]$ determined by the derivation $\eta$,
  whose underlying $k$-module spectrum is equivalent to $\tmop{fib} (A
  \rightarrow M)$.
\end{remark}

The key construction is the following:

{\construction{\label{cons:circle-shearing}Let $k$ be a commutative ring. We
construct a functor
\[ \tmop{DG}_-^{[0, 1]} \tmop{DAlg}_k \longrightarrow \tmop{DAlg} (D
   (k)^{B\mathbb{T}}) \]
as follows. Let $(A, M, \eta)$ be a triple corresponding to an object of
$\tmop{DG}_-^{[0, 1]} \tmop{DAlg}_k$ as in \Cref{rem:DG-alg-sq-zero-ext}.
Equip $A$ and $M$ with trivial $\mathbb{T}$-action, we get a derivation $\eta
\of A \rightarrow M$ in $\tmop{DAlg} (D (k)^{B\mathbb{T}})$. Now we compose
this derivation with the map $M = M \otimes \mathbb{S} \xrightarrow{M \otimes
a_{\lambda}} M \otimes \mathbb{S}^{\lambda} = \Sigma^{\lambda} M$, obtaining
another derivation $\tau \of A \rightarrow \Sigma^{\lambda} M$ in $\tmop{DAlg}
(D (k)^{B\mathbb{T}})$. We take the square-zero extension $\left[ A
\xrightarrow{\tau} \Sigma^{\lambda} M \right]$ of $A$ by $\Sigma^{\lambda} M
[- 1]$ determined by the derivation $\tau$. This gives rise to a functor
\begin{eqnarray*}
  \tmop{DG}_-^{[0, 1]} \tmop{DAlg}_k & \longrightarrow & \tmop{DAlg} (D
  (k)^{B\mathbb{T}})\\
  (A, M, \eta) & \longmapsto & \left[ A \xrightarrow{\tau} \Sigma^{\lambda} M
  \right] .
\end{eqnarray*}
Moreover, the map $M \otimes a_{\lambda} \of M \rightarrow \Sigma^{\lambda} M$
induces a map $\left[ A \xrightarrow{\eta} M \right] \rightarrow \left[ A
\xrightarrow{\tau} \Sigma^{\lambda} M \right]$ of derived commutative algebras
in $D (k)^{B\mathbb{T}}$, where $\left[ A \xrightarrow{\eta} M \right]$ is
equipped with trivial $\mathbb{T}$-action.}}

\begin{remark}
  In fact, in Construction~\ref{cons:circle-shearing}, we can endow the target
  with a filtered circle action. Since this notion is not used in the text, we
  will not discuss this filtered enhancement.
\end{remark}

We need the following categorification of Tate construction.

\begin{definition}[{\cite[§2]{Petrov2019}}]
  Let $k$ be a commutative ring spectrum. The {\tmdfn{$\mathbb{T}$-Tate
  category}} $D (k)^{t\mathbb{T}}$ is the quotient category
  \[ D (k)^{B\mathbb{T}} / \langle M \otimes \mathbb{T} \barsuchthat M \in D
     (k) \rangle \]
  where $\langle M \otimes \mathbb{T} \barsuchthat M \in D (k) \rangle
  \subseteq D (k)^{B\mathbb{T}}$ is the stable subcategory generated by
  induced $\mathbb{T}$-representations.
\end{definition}

The magic of Construction~\ref{cons:circle-shearing} is that it does not
change the Tate construction, which is formally formulated as follows.

\begin{proposition}
  \label{prop:Tate-sq-zero}Let $k$ be a commutative ring, and $(A, M, \eta)$ a
  triple corresponding to an object of $\tmop{DG}_-^{[0, 1]} \tmop{DAlg}_k$.
  Then the map
  \[ \left[ A \xrightarrow{\eta} M \right] \longrightarrow \left[ A
     \xrightarrow{\tau} \Sigma^{\lambda} M \right] \]
  in $\tmop{DAlg} (D (k)^{B\mathbb{T}})$ constructed in
  Construction~\ref{cons:circle-shearing}, where $\left[ A \xrightarrow{\eta}
  M \right]$ is equipped with trivial $\mathbb{T}$-action, is an equivalence
  after passing along the quotient functor $D (k)^{B\mathbb{T}} \rightarrow D
  (k)^{t\mathbb{T}}$.
\end{proposition}

\begin{proof}
  The point is that $\tmop{fib} (a_{\lambda} \of \mathbb{S} \rightarrow
  \mathbb{S}^{\lambda}) \simeq \mathbb{S} \otimes \mathbb{T}$ is an induced
  $\mathbb{T}$-representation.
\end{proof}

Let $k \rightarrow R$ be a smooth map of commutative algebras. The next goal,
\Cref{cor:HKR-dR-shear}, is to show that, if we apply this construction to the
truncated de Rham complex
\[ 0 \longrightarrow R \xrightarrow{\mathd} \Omega_{R / k}^1 \longrightarrow 0
\]
we get the associated $[0, 1]$-graded piece of the Hochschild homology
$\tmop{HH} (R / k)$. We explain this more formally for maps $k \rightarrow R$
of derived commutative rings as follows.

\begin{remark}
  \label{rem:forget-T-equivar-deriv-symm-mon}Let $k$ be a derived commutative
  ring. The forgetful functor $D (k)^{B\mathbb{T}} = \tmop{coAlg}_{k
  [\mathbb{T}]^{\vee}} (D (k)) \rightarrow D (k)$ maps derived commutative
  algebras in $D (k)^{B\mathbb{T}}$ to derived commutative $k$-algebras, and
  it carries the adic filtration, as defined in
  {\cite[Cons~3.2.24]{Magidson2024}}, on a map $f$ in $\tmop{DAlg} (D
  (k)^{B\mathbb{T}})$ to the adic filtration on the underlying map $f$ in
  $\tmop{DAlg}_k$.
\end{remark}

\begin{lemma}
  \label{lem:HKR-ad}Let $k \rightarrow A$ be a map of derived commutative
  rings. Then the HKR-filtration on $\tmop{HH} (A / k)$, as a filtered derived
  commutative algebra in $D (k)^{B\mathbb{T}}$, coincides with the adic
  filtration on the map $\tmop{HH} (A / k) \rightarrow A$ in $\tmop{DAlg} (D
  (k)^{B\mathbb{T}})$.
\end{lemma}

\begin{proof}
  By the universal property of adic filtration (defined as a left adjoint in
  {\cite[Cons~3.2.24]{Magidson2024}})\footnote{In the special case of a
  surjection of animated commutative rings, this universal property is
  explained in {\cite[Prop~13.3]{Clausen2022}} as well.}, we get a map
  \[ \tmop{Fil}_{\tmop{ad}} (\tmop{HH} (A / k) \rightarrow A) \longrightarrow
     \tmop{Fil}_{\tmop{HKR}} \tmop{HH} (A / k) \]
  in $\tmop{Fil} \tmop{DAlg} (D (k)^{B\mathbb{T}})$. To show that this is an
  equivalence, by conservativity, it suffices to show the equivalence after
  passing along the forgetful functor $\tmop{Fil} \tmop{DAlg} (D
  (k)^{B\mathbb{T}}) \rightarrow \tmop{Fil} \tmop{DAlg}_k$, and then it
  follows from \Cref{rem:forget-T-equivar-deriv-symm-mon} that the left hand
  side is still the adic filtration of the map $\tmop{HH} (A / k) \rightarrow
  A$.
  
  As both the source and the target preserve geometric realizations in $A \in
  \tmop{DAlg}_k$, and the $\infty$-category $\tmop{DAlg}_k$ is generated by
  $\{ \tmop{LSym}_k (M) \barsuchthat M \in D (k) \}$ under geometric
  realizations, we reduce to the special case that $A = \tmop{LSym}_k (M)$ for
  some $M \in D (k)$. As in the proof of {\cite[Lem~6.2.7]{Raksit2020}}, we
  have
  \[ \tmop{Fil}_{\tmop{HKR}} \tmop{HH} (A / k) = \tmop{LSym}_A^{\geqslant
     \mathord{\ast}} (M [1] \otimes_k^{\mathbb{L}} A) \]
  which does coincide with the adic filtration of the map $\tmop{HH} (A / k)
  \rightarrow A$.
\end{proof}

\begin{remark}
  \label{rem:adic-01-sq-zero}Let $\mathcal{C}$ be $D (k)$ or $D
  (k)^{B\mathbb{T}}$, and $A \rightarrow B$ be a map of derived commutative
  algebras in $\mathcal{C}$. Then the derived commutative algebra
  $\tmop{gr}_{\tmop{ad}}^{[0, 1]} (A \rightarrow B)$ can be realized as the
  square-zero extension
  \[ \left[ B \xrightarrow{\mathd} L_B \xrightarrow{\tmop{can}} L_{B / A}
     \right] \]
  of $B$ by $L_{B / A} [- 1]$ in $\mathcal{C}$.
\end{remark}

\begin{lemma}
  \label{lem:cot-cplx-HH}Let $k$ be a commutative ring, and $A$ a derived
  commutative $k$-algebra. Then the canonical map
  \[ L_{\tmop{HH} (A / k) / k} \otimes_{\tmop{HH} (A / k)}^{\mathbb{L}} A
     \longrightarrow L_{A / k} \]
  in $\tmop{Mod}_A (D (k)^{B\mathbb{T}})$ can be identified with the map
  \[ L_{A / k} \otimes_k^{\mathbb{L}} k [\mathbb{T}] \longrightarrow L_{A / k}
  \]
  induced by the $\mathbb{T}$-equivariant projection $\mathbb{T} \rightarrow
  \mathord{\ast}$.
\end{lemma}

\begin{proof}
  Let $X \in \tmop{An}$ be an anima, and $B_X \assign
  A^{\otimes_k^{\mathbb{L}} X}$ the $X$-th tensor power of $A$. Then we have
  an equivalence
  \[ L_{B_X / k} \otimes_{B_X}^{\mathbb{L}} A \simeq L_{A / k}
     \otimes_k^{\mathbb{L}} k [X] \]
  in $\tmop{Mod}_A (D (k))$, which is functorial in $X \in \tmop{An}$. Now we
  restrict this to the full subcategory $(B\mathbb{T})^{\vartriangleright}
  \subseteq \tmop{An}$ spanned by $\left\{ \mathbb{T}, \mathord{\ast}
  \right\}$, we obtained the identification that we want.
\end{proof}

\begin{proposition}
  \label{prop:HKR-sqzero}Let $k$ be a commutative ring, and $A$ an animated
  commutative $k$-algebra. Then the derived commutative algebra
  $\tmop{gr}_{\tmop{HKR}}^{[0, 1]} \tmop{HH} (A / k)$ can be identified with
  the square-zero extension
  \[ \left[ A \xrightarrow{\mathd} L_{A / k} \xrightarrow{L_{A / k} \otimes
     a_{\lambda}} \Sigma^{\lambda} L_{A / k} \right] \]
  of $A$ by $\Sigma^{\lambda - 1} L_{A / k} \simeq L_{A / k} [1]$ (via a
  complex orientation of $k$) in $\tmop{DAlg} (D (k)^{B\mathbb{T}})$.
\end{proposition}

\begin{proof}
  By \Cref{rem:adic-01-sq-zero,lem:HKR-ad}, it suffices to identify the map
  \[ f \of L_{A / k} \longrightarrow L_{A / \tmop{HH} (A / k)} \]
  in $\tmop{Mod}_A (D (k)^{B\mathbb{T}})$ with the map
  \[ L_{A / k} \xrightarrow{L_{A / k} \otimes a_{\lambda}} \Sigma^{\lambda}
     L_{A / k} . \]
  We examine the transitivity sequence
  \[ L_{\tmop{HH} (A / k) / k} \otimes_{\tmop{HH} (A / k)}^{\mathbb{L}} A
     \longrightarrow L_{A / k} \xrightarrow{f} L_{A / \tmop{HH} (A / k)} \]
  in $\tmop{Mod}_A (D (k)^{B\mathbb{T}})$ associated to the composite map
  \[ k \longrightarrow \tmop{HH} (A / k) \longrightarrow A \]
  in $\tmop{DAlg} (D (k)^{B\mathbb{T}})$. By \Cref{lem:cot-cplx-HH}, the map
  $f$ fits into a fiber sequence
  \[ L_{A / k} \otimes_k^{\mathbb{L}} k [\mathbb{T}] \longrightarrow L_{A / k}
     \xrightarrow{f} L_{A / \tmop{HH} (A / k)} \]
  in $\tmop{Mod}_A (D (k)^{B\mathbb{T}})$, and the result follows from the
  fiber sequence
  \[ \mathbb{S} \otimes \mathbb{T} \longrightarrow \mathbb{S}
     \xrightarrow{a_{\lambda}} \mathbb{S}^{\lambda} . \]
\end{proof}

\begin{corollary}
  \label{cor:HKR-dR-shear}Let $k$ be a commutative ring, and $A$ an animated
  commutative $k$-algebra. Then the derived commutative algebra
  $\tmop{gr}_{\tmop{HKR}}^{[0, 1]} \tmop{HH} (A / k) \in \tmop{DAlg} (D
  (k)^{B\mathbb{T}})$ can be identified with the image of $L \Omega_{A /
  k}^{[0, 1]} \in \tmop{DG}_-^{[0, 1]} \tmop{DAlg}_k$ under the functor in
  Construction~\ref{cons:circle-shearing}.
\end{corollary}

Now we restrict to smooth maps of dimension$\nosymbol \leq 1$. In this case,
both HKR-filtration and Hodge-filtration are concentrated in weights $[0, 1]$,
thus we have

{\construction{\label{cons:map-dR-HH}{\dueto{Raksit}}Let $k$ be a commutative
ring, and $A$ a smooth $k$-algebra of dimension$\nosymbol \leq 1$. Then by
Construction~\ref{cons:circle-shearing} and \Cref{cor:HKR-dR-shear}, we get a
map
\[ \tmop{dR}_{R / k} \longrightarrow \tmop{HH} (R / k) \]
in $\tmop{DAlg} (D (k)^{B\mathbb{T}})$, which becomes an equivalence in $D
(k)^{t\mathbb{T}}$ by \Cref{prop:Tate-sq-zero}. This construction is
functorial in smooth maps $k \rightarrow R$ of relative dimension$\nosymbol
\leq 1$.}}

Now we want to apply this to analyze the de-completed $C_p$-Tate construction
$\tmop{HH} (\mathbb{F}_p /\mathbb{Z})^{\theta_{\mathbb{Z}} C_p}$.

\begin{remark}
  \label{rem:Borel-cpl-decpl-Tate}The composite lax symmetric monoidal functor
  \[ D (\mathbb{Z})^{B\mathbb{T}} \xrightarrow{(-)^h} \tmop{Mod}_{\mathbb{Z}}
     (\tmop{Sp}^{g^{<} \mathbb{T}}) \xrightarrow{(-)^{\theta_{\mathbb{Z}}
     C_p}} \tmop{Mod}_{\mathbb{Z}^{t C_p}} (\tmop{Sp}^{B (\mathbb{T}/ C_p)})
  \]
  factors canonically through $D (\mathbb{Z})^{t\mathbb{T}}$.
\end{remark}

{\construction{\label{cons:map-HH-PD-struct}Applying
Construction~\ref{cons:map-dR-HH} to the coCartesian diagram
\[ \begin{array}{ccc}
     k [t] & \xrightarrow{t \mapsto p} & k\\
     \longdownarrow &  & \longdownarrow\\
     k & \longrightarrow & k \otimes_{\mathbb{Z}}^{\mathbb{L}} \mathbb{F}_p
   \end{array}, \]
of animated $k$-algebras, we get a map
\[ \tmop{dR}_{(k \otimes_{\mathbb{Z}}^{\mathbb{L}} \mathbb{F}_p) / k}
   \longrightarrow \tmop{HH} (k \otimes_{\mathbb{Z}}^{\mathbb{L}} \mathbb{F}_p
   / k)^{\eta_k} \]
of derived commutative algebras in $\tmop{Mod}_{\mathbb{Z}} (\tmop{Sp}^{g^{<}
\mathbb{T}})$, which becomes an equivalence after taking $(-)^{\Phi C_p}$ by
\Cref{rem:Borel-cpl-decpl-Tate}. In particular, when $k =\mathbb{Z}_{(p)}$,
the map
\[ \tmop{dR}_{\mathbb{F}_p /\mathbb{Z}_{(p)}} = D_{\mathbb{Z}_{(p)}} (p)
   \longrightarrow \mathbb{Z}_{(p)} \]
of augmented (derived) commutative $\mathbb{Z}_{(p)}$-algebras (where
$D_{\mathbb{Z}_{(p)}} (p)$ is the PD-envelope of $(p) \subseteq
\mathbb{Z}_{(p)}$), induced by the PD-structure on $(\mathbb{Z}_{(p)}, (p))$,
gives rise to a map
\[ \tmop{HH} (\mathbb{F}_p /\mathbb{Z})^{\theta_{\mathbb{Z}} C_p} \rightarrow
   \mathbb{Z}^{t C_p} \]
of $\mathbb{E}_{\infty}$-$\mathbb{Z}^{t C_p}$-algebras in $\tmop{Sp}^{B
(\mathbb{T}/ C_p)}$.}}

{\construction{\label{cons:THH-Fp-Z-tCp}Composing with the map $\tmop{THH}
(\mathbb{F}_p)^h \rightarrow \tmop{HH} (\mathbb{F}_p
/\mathbb{Z})^{\eta_{\mathbb{Z}}}$ of $\mathbb{E}_{\infty}$-$\tmop{THH}
(\mathbb{Z})^h$-algebras in $\tmop{Sp}^{g^{<} \mathbb{T}}$, we get a composite
map
\[ \tmop{THH} (\mathbb{F}_p)^{t C_p} \longrightarrow \tmop{HH} (\mathbb{F}_p
   /\mathbb{Z})^{\theta_{\mathbb{Z}} C_p} \longrightarrow \mathbb{Z}^{t C_p}
\]
of $\mathbb{E}_{\infty}$-$\tmop{THH} (\mathbb{Z})^{t C_p}$-algebras in
$\tmop{Sp}^{B (\mathbb{T}/ C_p)}$.}}

\begin{remark}
  The map in Construction~\ref{cons:THH-Fp-Z-tCp} lives in the Borel
  equivariant world, while our construction relies on genuine equivariant
  homotopy theory.
\end{remark}

Our next goal is to show that the map in Construction~\ref{cons:THH-Fp-Z-tCp}
is an equivalence. We check it after modulo $p$. The point is that, although
we do not know whether this map on the nose acquires a $\mathbb{Z}$-linear
structures, its modulo $p$ has an $\mathbb{F}_p$-linear structure. More
precisely, our argument follows the following strategy\footnote{We thank
A.~{\tmname{Petrov}} for explaining this.} in {\cite{Petrov2018}}:
\begin{enumerate}
  \item Give the composite map $\tmop{THH} (\mathbb{F}_p)^{t C_p} \rightarrow
  \mathbb{Z}^{t C_p} \rightarrow \mathbb{F}_p^{t C_p}$ a $\mathbb{Z}$-linear
  structure (Construction~\ref{cons:Z-lin-struct-THH-Fp-Fp-tCp});
  
  \item Show that, under this $\mathbb{Z}$-linear structure, the induced map
  $\tmop{THH} (\mathbb{F}_p)^{t C_p} \otimes_{\mathbb{Z}}^{\mathbb{L}}
  \mathbb{F}_p \rightarrow \mathbb{F}_p^{t C_p}$ is an equivalence
  (\Cref{lem:THH-Fp-tCp-base-chg-equiv});
  
  \item Show that the induced map above coincides with the modulo $p$
  reduction of the map $\tmop{THH} (\mathbb{F}_p)^{t C_p} \rightarrow
  \mathbb{Z}^{t C_p}$ in Construction~\ref{cons:THH-Fp-Z-tCp}
  (\Cref{lem:mod-p-non-Z-lin-map}).
\end{enumerate}
\begin{remark}
  Concretely, the map $D_{\mathbb{Z}_{(p)}} (p) \rightarrow \mathbb{Z}_{(p)}$
  is given by the map
  \begin{eqnarray*}
    \mathbb{Z}_{(p)} \left[ t, \frac{t^2}{2!}, \frac{t^3}{3!}, \ldots \right]
    / (t - p) & \longrightarrow & \mathbb{Z}_{(p)}\\
    \frac{t^n}{n!} & \longmapsto & \frac{p^n}{n!},
  \end{eqnarray*}
  and therefore, after (derived) modulo $p$, it becomes the augmentation map
  \begin{eqnarray*}
    \mathbb{F}_p \left[ t, \frac{t^p}{p}, \frac{t^{p^2}}{p^{p + 1}}, \ldots
    \right] / (t) & \longrightarrow & \mathbb{F}_p\\
    \frac{t^{p^r}}{p^{p^{r - 1} + \cdots + 1}} & \longmapsto & 0.
  \end{eqnarray*}
  It follows that, the map $\tmop{dR}_{\mathbb{F}_p /\mathbb{Z}_{(p)}}
  \otimes_{\mathbb{Z}_{(p)}}^{\mathbb{L}} \mathbb{F}_p \rightarrow
  \mathbb{F}_p$ induced by the PD-structure on $(p) \subseteq
  \mathbb{Z}_{(p)}$ coincides with the map $\tmop{dR}_{k
  \otimes_{\mathbb{Z}}^{\mathbb{L}} \mathbb{F}_p / k} \rightarrow \tmop{dR}_{k
  / k} = k$ induced by the multiplication map $k
  \otimes_{\mathbb{Z}}^{\mathbb{L}} \mathbb{F}_p \rightarrow k$, for derived
  commutative $\mathbb{F}_p$-algebras $k$.
\end{remark}

We get a de-completed version of {\cite[Cor~2.7]{Petrov2019}} (without
restricting to odd primes):

\begin{corollary}
  \label{cor:mod-p-cris-map}The modulo $p$ reduction of the map $\tmop{HH}
  (\mathbb{F}_p /\mathbb{Z})^{\Phi C_p} \rightarrow \mathbb{Z}^{t C_p}$ of
  $\mathbb{E}_{\infty}$-$\mathbb{Z}^{t C_p}$-algebras in $\tmop{Sp}^{B
  (\mathbb{T}/ C_p)}$ in Construction~\ref{cons:map-HH-PD-struct} becomes the
  map $\tmop{HH} (k \otimes_{\mathbb{Z}}^{\mathbb{L}} \mathbb{F}_p / k)^{\Phi
  C_p} \rightarrow \tmop{HH} (k / k)^{\Phi C_p}$ of
  $\mathbb{E}_{\infty}$-$\mathbb{F}_p^{t C_p}$-algebras in $\tmop{Sp}^{B
  (\mathbb{T}/ C_p)}$ induced by the multiplication map $k
  \otimes_{\mathbb{Z}}^{\mathbb{L}} \mathbb{F}_p \rightarrow k$, for $k
  =\mathbb{F}_p$\footnote{Although we only apply to the case that $k
  =\mathbb{F}_p$, this notation allows us to distinguish different
  $\mathbb{F}_p$-algebra structures on $\mathbb{F}_p
  \otimes_{\mathbb{Z}}^{\mathbb{L}} \mathbb{F}_p$.}.
\end{corollary}

{\construction{\label{cons:Z-lin-struct-THH-Fp-Fp-tCp}Let $k =\mathbb{F}_p$.
The composite map in Construction~\ref{cons:THH-Fp-Z-tCp} fits into a
commutative diagram
\[ \begin{array}{ccccc}
     \tmop{THH} (k)^{t C_p} & \longrightarrow & \tmop{HH} (k
     /\mathbb{Z})^{\theta_{\mathbb{Z}} C_p} & \longrightarrow & \mathbb{Z}^{t
     C_p}\\
     &  & \longdownarrow &  & \longdownarrow\\
     &  & \tmop{HH} (k \otimes_{\mathbb{Z}}^{\mathbb{L}} \mathbb{F}_p
     /\mathbb{F}_p)^{\theta_{\mathbb{F}_p} C_p} & \longrightarrow &
     \mathbb{F}_p^{t C_p}
   \end{array} \]
of $\mathbb{E}_{\infty}$-$\tmop{THH} (\mathbb{Z})^{t C_p}$-algebras in
$\tmop{Sp}^{B (\mathbb{T}/ C_p)}$. It follows from \Cref{cor:mod-p-cris-map}
that the composite map $\tmop{THH} (k)^{t C_p} \rightarrow \mathbb{F}_p^{t
C_p}$ is equivalent to applying $(-)^{t C_p}$ to the map
\[ \tmop{THH} (\mathbb{F}_p) \longrightarrow \tmop{HH} (\mathbb{F}_p
   /\mathbb{F}_p) \]
induced by the map $(\mathbb{S} \rightarrow \mathbb{F}_p) \rightarrow
(\mathbb{F}_p \rightarrow \mathbb{F}_p)$ of maps of
$\mathbb{E}_{\infty}$-rings, which acquires a structure of maps of
$\mathbb{E}_{\infty}$-$\mathbb{Z}^{t C_p}$-algebras in $\tmop{Sp}^{B
(\mathbb{T}/ C_p)}$. This map is also the same as the multiplication map
$\tmop{THH} (\mathbb{F}_p) \rightarrow \mathbb{F}_p$.}}

\begin{lemma}
  \label{lem:THH-Fp-tCp-base-chg-equiv}The map $\tmop{THH} (\mathbb{F}_p)^{t
  C_p} \otimes_{\mathbb{Z}}^{\mathbb{L}} \mathbb{F}_p \rightarrow
  \mathbb{F}_p^{t C_p}$ induced by the $\mathbb{Z}^{t C_p}$-linear map
  $\tmop{THH} (\mathbb{F}_p)^{t C_p} \rightarrow \mathbb{F}_p^{t C_p}$ in
  Construction~\ref{cons:Z-lin-struct-THH-Fp-Fp-tCp} is an equivalence.
\end{lemma}

\begin{proof}
  The composite map $\mathbb{Z} \rightarrow \tmop{THH} (\mathbb{F}_p)
  \rightarrow \mathbb{F}_p$ is the canonical map since it is a map of
  commutative algebras in the heart of $D (\mathbb{Z})^{B\mathbb{T}}$. By
  {\cite[Cor~IV.4.13]{Nikolaus2018}}, the first map becomes an equivalence
  after taking $(-)^{t C_p}$, and the result follows.
\end{proof}

Finally, we are in the following slightly tricky situation: given two
$\mathbb{Z}$-modules $M$ and $N$, along with a (non-$\mathbb{Z}$-linear) map
$f \of M \rightarrow N$. Suppose that the composite map $g \of M \rightarrow N
\rightarrow N \otimes_{\mathbb{Z}}^{\mathbb{L}} \mathbb{F}_p$ carries a
$\mathbb{Z}$-linear structure, which induces a map $M
\otimes_{\mathbb{Z}}^{\mathbb{L}} \mathbb{F}_p \rightarrow N
\otimes_{\mathbb{Z}}^{\mathbb{L}} \mathbb{F}_p$. We claim that this map
coincides with the modulo $p$ reduction of $f$. We break this into two steps.
The first is to produce the map $M /^{\mathbb{L}} \mathbb{F}_p \rightarrow N
\otimes_{\mathbb{Z}}^{\mathbb{L}} \mathbb{F}_p$ without the
$\mathbb{Z}$-linear structure on $g$.

{\construction{\label{cons:map-to-Fp-mod-induce}Let $\mathcal{C}$ be a
presentable stable $\infty$-category, and $M \in \mathcal{C}$ an object, and
$L \in \tmop{Mod}_{\mathbb{F}_p} (\mathcal{C})$ an $\mathbb{F}_p$-module in
$\mathcal{C}$, and $f \of M \rightarrow L$ a map in $\mathcal{C}$. This gives
rise to a composite map
\[ M / p \longrightarrow L / p \simeq L \otimes_{\mathbb{Z}}^{\mathbb{L}}
   \mathbb{F}_p \xrightarrow{\tmop{mult}} L. \]}}

\begin{lemma}
  \label{lem:spectra-map-to-Fp-mod}Let $\mathcal{C}$ be a presentable stable
  $\infty$-category, and $M \in \tmop{Mod}_{\mathbb{Z}} (\mathcal{C})$ a
  $\mathbb{Z}$-module in $\mathcal{C}$, and $L \in \tmop{Mod}_{\mathbb{F}_p}
  (\mathcal{C})$ an $\mathbb{F}_p$-module in $\mathcal{C}$, and $f \of M
  \rightarrow L$ a $\mathbb{Z}$-linear map in $\mathcal{C}$. Then the map $M
  \otimes_{\mathbb{Z}}^{\mathbb{L}} \mathbb{F}_p \rightarrow L$ in
  $\tmop{Mod}_{\mathbb{F}_p} (\mathcal{C})$, as a map in $\mathcal{C}$,
  coincides with the map $M / p \rightarrow L$ in
  Construction~\ref{cons:map-to-Fp-mod-induce}.
\end{lemma}

\begin{proof}
  The map $f \otimes_{\mathbb{Z}}^{\mathbb{L}} \mathbb{F}_p \of M
  \otimes_{\mathbb{Z}}^{\mathbb{L}} \mathbb{F}_p \rightarrow L
  \otimes_{\mathbb{Z}}^{\mathbb{L}} \mathbb{F}_p$ in
  $\tmop{Mod}_{\mathbb{F}_p} (\mathcal{C})$, as a map in $\mathcal{C}$,
  coincides with the map $f / p \of M / p \rightarrow L / p$. The result
  follows.
\end{proof}

\begin{lemma}
  \label{lem:mod-p-non-Z-lin-map}Let $\mathcal{C}$ be a presentable stable
  $\infty$-category, with $M$ and $N$ two $\mathbb{Z}$-modules in
  $\mathcal{C}$, and $f \of M \rightarrow N$ a map in $\mathcal{C}$. We give
  the composite map $g \of M \xrightarrow{f} N \rightarrow N
  \otimes_{\mathbb{Z}}^{\mathbb{L}} \mathbb{F}_p$ a $\mathbb{Z}$-linear
  structure, which induces a map $h \of M \otimes_{\mathbb{Z}}^{\mathbb{L}}
  \mathbb{F}_p \rightarrow N \otimes_{\mathbb{Z}}^{\mathbb{L}} \mathbb{F}_p$
  in $\tmop{Mod}_{\mathbb{F}_p} (\mathcal{C})$. Then as a map in
  $\mathcal{C}$, the map $h$ coincides with the map $f / p$.
\end{lemma}

\begin{proof}
  By \Cref{lem:spectra-map-to-Fp-mod}, the map $h$, as a map in $\mathcal{C}$,
  coincides with the composite
  \[ M / p \xrightarrow{g / p} (N \otimes_{\mathbb{Z}}^{\mathbb{L}}
     \mathbb{F}_p) / p \simeq N \otimes_{\mathbb{Z}}^{\mathbb{L}} \mathbb{F}_p
     \otimes_{\mathbb{Z}}^{\mathbb{L}} \mathbb{F}_p \xrightarrow{\tmop{id}_N
     \otimes \tmop{mult}} N \otimes_{\mathbb{Z}}^{\mathbb{L}} \mathbb{F}_p .
  \]
  Thus it suffices to show that the composite
  \[ N / p \rightarrow (N \otimes_{\mathbb{Z}}^{\mathbb{L}} \mathbb{F}_p) / p
     \simeq N \otimes_{\mathbb{Z}}^{\mathbb{L}} \mathbb{F}_p
     \otimes_{\mathbb{Z}}^{\mathbb{L}} \mathbb{F}_p \xrightarrow{\tmop{id}_N
     \otimes \tmop{mult}} N \otimes_{\mathbb{Z}}^{\mathbb{L}} \mathbb{F}_p \]
  is equivalent to identity. Note that the first map can be identified with
  \[ N \otimes_{\mathbb{Z}}^{\mathbb{L}} \mathbb{F}_p \xrightarrow{\tmop{id}_N
     \otimes 1 \otimes \tmop{id}_{\mathbb{F}_p}} N
     \otimes_{\mathbb{Z}}^{\mathbb{L}} \mathbb{F}_p
     \otimes_{\mathbb{Z}}^{\mathbb{L}} \mathbb{F}_p \]
  and the result follows from the fact that the composite map
  \[ \mathbb{F}_p \xrightarrow{1 \otimes \tmop{id}_{\mathbb{F}_p}}
     \mathbb{F}_p \otimes_{\mathbb{Z}}^{\mathbb{L}} \mathbb{F}_p
     \xrightarrow{\tmop{mult}} \mathbb{F}_p \]
  of animated rings is identity.
\end{proof}

\begin{proposition}
  \label{prop:THH-Fp-tCp-Z-tCp-eq}The map $\tmop{THH} (\mathbb{F}_p)^{t C_p}
  \rightarrow \mathbb{Z}^{t C_p}$ of $\mathbb{E}_{\infty}$-$\tmop{THH}
  (\mathbb{Z})^{t C_p}$-algebras in Construction~\ref{cons:THH-Fp-Z-tCp} is an
  equivalence.
\end{proposition}

\begin{proof}
  Since both the source and the target is $p$-complete by
  {\cite[Lem~I.2.9]{Nikolaus2018}}, it suffices to check that this map is an
  equivalence after modulo $p$. Now the result follows from
  \Cref{lem:THH-Fp-tCp-base-chg-equiv,lem:mod-p-non-Z-lin-map}.
\end{proof}

\begin{corollary}
  \label{cor:THH-Z-mod-crys-comp}Let $M$ be a $\tmop{THH}
  (\mathbb{Z})^h$-module in $\tmop{Sp}^{g^{<} \mathbb{T}}$. Then there exists
  a lax symmetric monoidal (in $M$) equivalence
  \[ \left( (M \otimes_{\tmop{THH} (\mathbb{Z})^h}^{\mathbb{L}} \tmop{THH}
     (\mathbb{F}_p)^h)^{\theta_{\tmop{THH} (\mathbb{F}_p)} \mathbb{T}}
     \right)_p^{\wedge} \rightarrow ((M \otimes_{\tmop{THH}
     (\mathbb{Z})^h}^{\mathbb{L}} \mathbb{Z})^{\theta_{\mathbb{Z}}
     \mathbb{T}})_p^{\wedge} \]
  of $\tmop{TP} (\mathbb{Z})$-modules.
\end{corollary}

\begin{proof}
  By \Cref{cor:de-cpl-Borel-Tate-THH-Z}, it suffices to produce a symmetric
  monoidal equivalence
  \[ (M \otimes_{\tmop{THH} (\mathbb{Z})^h}^{\mathbb{L}} \tmop{THH}
     (\mathbb{F}_p)^h)^{\Phi C_p} \longrightarrow (M \otimes_{\tmop{THH}
     (\mathbb{Z})^h}^{\mathbb{L}} \mathbb{Z})^{\Phi C_p} \]
  in $\tmop{Mod}_{\tmop{THH} (\mathbb{Z})^{t C_p}} (\tmop{Sp}^{B (\mathbb{T}/
  C_p)})$, which is an immediate consequence of
  \Cref{prop:THH-Fp-tCp-Z-tCp-eq}.
\end{proof}

\begin{remark}
  \label{rem:crys-dR-comp-odd-primes}When $p$ is an odd prime,
  {\cite[Prop~5.2.1]{Devalapurkar2025}} shows that the map $\mathbb{Z}
  \rightarrow \tmop{THH} (\mathbb{F}_p)$ induces an equivalence $\mathbb{Z}^{t
  C_p} \rightarrow \tmop{THH} (\mathbb{F}_p)^{t C_p}$ of
  $\mathbb{E}_{\infty}$-$\tmop{THH} (\mathbb{Z})^{t C_p}$-algebras, which
  implies \Cref{cor:THH-Z-mod-crys-comp} (although a priori, the map there
  might be different from the one here), and the map even has a
  $\mathbb{Z}^{t\mathbb{T}}$-module structure.
\end{remark}

\begin{corollary}
  \label{cor:Z-lin-cat-crys-comp}Let $\mathcal{C}$ be a dualizable presentable
  stable $\mathbb{Z}$-linear $\infty$-category. Then there exists a lax
  symmetric monoidal (in $\mathcal{C}$) equivalence
  \[ \tmop{TP}^{\tmop{poly} /\mathbb{F}_p} (\mathcal{C}
     \otimes_{\mathbb{Z}}^{\mathbb{L}} \mathbb{F}_p)_p^{\wedge}
     \longrightarrow \tmop{HP}^{\tmop{poly}}
     (\mathcal{C}/\mathbb{Z})_p^{\wedge} \]
  of $\tmop{TP} (\mathbb{Z})$-modules.
\end{corollary}

\begin{proof}
  Apply \Cref{cor:THH-Z-mod-crys-comp} to $M = \tmop{THH} (\mathcal{C})
  \otimes_{\tmop{THH} (\mathbb{Z})}^{\mathbb{L}} \tmop{THH} (\mathbb{Z})^h$.
\end{proof}

\section{Comparison to $\tmop{THH}$}\label{sec:comp-THH}

Recall that the {\tmdfn{Cartier isomorphism}} identifies algebraic de Rham
cohomology groups of smooth $\mathbb{F}_p$-algebras with their algebraic
differential forms. We give two noncommutative analogues. In this section, we
discuss one of them, which compares polynomial co-periodic cyclic homology of
dualizable presentable stable $\mathbb{F}_p$-linear $\infty$-categories with
their topological Hochschild homology. This comparison was proved in
{\cite[Cor~11.15]{Kaledin2020}} for associative $\mathbb{F}_p$-algebras, using
Goodwillie derivative of (Hochschild--)Witt trace theory. We give two
arguments. Although the second argument is much shorter, the first argument
gives us more information, which is used in \Cref{sec:conj-fil}.

The key to the first argument is the observation that the $\mathbb{T}/
C_p$-equivariant $\mathbb{E}_1$-$\mathbb{Z}^{t C_p}$-module $\mathbb{F}_p^{t
C_p}$ is co-induced.

{\construction{\label{cons:rescale-fiber-seq}Let $n \in \mathbb{N}_{> 0}$.
There is a map
\[ \mathbb{Z}^{\mathbb{T}/ C_n} \longrightarrow \mathbb{Z}/ n \]
in $D (\mathbb{Z})^{B\mathbb{T}}$, where $\mathbb{Z}/ n$ is equipped with the
trivial action, i.e. $\varpi_1^{\ast} (\mathbb{Z}/ n)$. Indeed, this map is
taken to be represented\footnote{The $\mathbb{E}_1$-monoidal equivalence of
the monoidal $\infty$-category $D (\mathbb{Z})^{B\mathbb{T}}$ and the monoidal
$\infty$-category of mixed complexes up to quasi-isomorphism is explained in
details in {\cite[§5.4]{Leip2022}}.} by the following surjective map
\[ \begin{array}{ccc}
     0 &  & - 1\\
     \mathbb{Z} & \underset{n}{\overset{0}{\longrightleftarrows}} &
     \mathbb{Z}\\
     \longdownarrow &  & \longdownarrow\\
     \mathbb{Z}/ n & \underset{0}{\overset{0}{\longrightleftarrows}} & 0
   \end{array} \]
of mixed complexes concentrated in homological degrees $[- 1, 0]$. Since it is
surjective, it is a fibration in the projective model structure on mixed
complexes, with kernel being the mixed complex $\left( n\mathbb{Z}
\underset{n}{\overset{0}{\longrightleftarrows}} \mathbb{Z} \right)$, which
represents $\mathbb{Z}^{\mathbb{T}} \in D (\mathbb{Z})^{B\mathbb{T}}$.
Consequently, we get a fiber sequence
\[ \mathbb{Z}^{\mathbb{T}} \longrightarrow \mathbb{Z}^{\mathbb{T}/ C_n}
   \longrightarrow \mathbb{Z}/ n \]
in $D (\mathbb{Z})^{B\mathbb{T}}$.}}

\begin{remark}
  One can also construct the map and the fiber sequence in
  Construction~\ref{cons:rescale-fiber-seq} by the {\tmdfn{gold relation}}
  $a_{\lambda^n} u_{\lambda^m} = (n / m) a_{\lambda^m} u_{\lambda^n}$ for $(m,
  n) \in \mathbb{N}^2_{> 0}$ with $m \divides n$ in {\cite[Lem~3.6]{Hill2017}}
  (where $u_{\lambda^m}$'s become equivalence after forgetting to the Borel
  equivariant objects), where $\lambda^n$ is the complex $S^1$-representation
  $S^1 \xrightarrow{(-)^n} \mathbb{C}^{\times}$ viewed as a real
  representation, and $a_{\lambda^n}$ is the {\tmdfn{Euler class}}
  $\mathbb{S}^0 \rightarrow \mathbb{S}^{\lambda^n}$, i.e. the map of Thom
  spectra of the inclusion $\{ 0 \} \subseteq \lambda^n$ of representations.
  Indeed, there is a fiber sequence
  \[ \Sigma_{\mathbb{T}}^{\infty}  [\mathbb{T}/ C_n]_+ \longrightarrow
     \mathbb{S}^0 \longrightarrow \mathbb{S}^{\lambda^n} \]
  of cyclonic spectra (and thus, of $\mathbb{T}$-equivariant spectra),
  cf.~{\cite[Obs~2.32]{Sulyma2020}}. Taking $\mathbb{Z}$-linear dual, and
  applying the ``octahedral axiom'' to the gold relation $a_{\lambda^n} \sim
  na_{\lambda^1}$, we get the fiber sequence in
  Construction~\ref{cons:rescale-fiber-seq}.
  
  This allows us to establish results beyond char $p$, as observed by Yuri
  {\tmname{Sulyma}} in {\cite[Lem~4.9]{Sulyma2020}}.
\end{remark}

\begin{corollary}
  \label{cor:Z-mod-n-Cn-Tate}The map $\mathbb{Z}^{\mathbb{T}/ C_n} \rightarrow
  \mathbb{Z}/ n$ in $D (\mathbb{Z})^{B\mathbb{T}}$ as in
  Construction~\ref{cons:rescale-fiber-seq} induces an equivalence
  \[ \mathbb{Z}^{t C_n} \otimes_{\mathbb{Z}}^{\mathbb{L}}
     \mathbb{Z}^{\mathbb{T}/ C_n} \longrightarrow (\mathbb{Z}/ n)^{t C_n} \]
  in $\tmop{Mod}_{\mathbb{Z}^{t C_n}} (D (\mathbb{Z})^{B (\mathbb{T}/ C_n)})$.
\end{corollary}

\begin{proof}
  Applying $(-)^{t C_n}$ to the fiber sequence in
  Construction~\ref{cons:rescale-fiber-seq}, we see that the induced map
  $(\mathbb{Z}^{\mathbb{T}/ C_n})^{t C_n} \rightarrow (\mathbb{Z}/ n)^{t C_n}$
  in $D (\mathbb{Z})^{B (\mathbb{T}/ C_n)}$ is an equivalence (which even has
  an $\mathbb{E}_1$-structure). It suffices to establish an equivalence
  \[ \mathbb{Z}^{t C_n} \otimes_{\mathbb{Z}}^{\mathbb{L}}
     \mathbb{Z}^{\mathbb{T}/ C_n} \longrightarrow (\mathbb{Z}^{\mathbb{T}/
     C_n})^{t C_n} \]
  in $\tmop{Mod}_{\mathbb{Z}^{t C_n}} (D (\mathbb{Z})^{B (\mathbb{T}/ C_n)})$.
  This is given by equivalences
  \begin{eqnarray*}
    \mathbb{Z}^{t C_n} \otimes_{\mathbb{Z}}^{\mathbb{L}}
    \mathbb{Z}^{\mathbb{T}/ C_n} & \simeq & (\mathbb{Z}^{t C_n}
    \otimes_{\mathbb{Z}}^{\mathbb{L}} \mathbb{Z} [\mathbb{T}/ C_n]) [- 1]\\
    & \simeq & (\mathbb{Z}^{t C_n} \otimes (\mathbb{T}/ C_n)) [- 1]\\
    & \simeq & (\mathbb{Z} [\mathbb{T}/ C_n])^{t C_n} [- 1]\\
    & \simeq & (\mathbb{Z}^{\mathbb{T}/ C_n})^{t C_n}
  \end{eqnarray*}
  in $\tmop{Mod}_{\mathbb{Z}^{t C_n}} (D (\mathbb{Z})^{B (\mathbb{T}/ C_n)})$.
\end{proof}

\begin{warning}
  The map in \Cref{cor:Z-mod-n-Cn-Tate} does not carry any
  $\mathbb{E}_1$-structure when $n = 2$. Indeed, on the left hand side, the
  square-zero class $e \in \pi_{- 1} (\mathbb{Z}^{\mathbb{T}/ C_n})$ gives
  rise to a square-zero class $e \in \pi_{- 1} (\mathbb{Z}^{\mathbb{T}/ C_n}
  \otimes_{\mathbb{Z}}^{\mathbb{L}} \mathbb{Z}^{t C_n})$, while the homotopy
  ring $\pi_{\ast} (\mathbb{F}_2^{t C_2})$ is isomorphic to $\mathbb{F}_2
  \pparen{s}$ for a generator $s \in \pi_1 (\mathbb{F}_2^{t C_2})$, which is
  integral, thus any square-zero elements is necessarily zero.
\end{warning}

\begin{question}
  Is there any version of $\mathbb{E}_1$-enhancement of
  \Cref{cor:Z-mod-n-Cn-Tate}?
\end{question}

Now we compare polynomial periodic cyclic homology of dualizable presentable
stable $\mathbb{F}_p$-linear $\infty$-categories. Slightly more generally, we
consider $\tmop{THH} (\mathbb{F}_p)$-modules in cyclonic spectra. By
\Cref{lem:decpl-Tate-orb-lem-T} and the symmetric monoidal structure on
$(-)^{\Phi C_p}$, we have

\begin{lemma}
  \label{lem:poly-HP-formula}Let $k$ be a commutative algebra. Then the
  composite functor
  \[ \tmop{Mod}_{\tmop{THH} (k)} (\tmop{Sp}^{g^{<} \mathbb{T}})
     \xrightarrow{(-) \otimes_{\tmop{THH} (k)}^{\mathbb{L}} \underline{k}}
     \tmop{Mod}_{\underline{k}} (\tmop{Sp}^{g^{<} \mathbb{T}})
     \xrightarrow{((-)^{\theta_{\underline{k}} \mathbb{T}})_p^{\wedge}} D
     (k^{t\mathbb{T}})_p^{\wedge} . \]
  is equivalent to the composite functor
  \[ \begin{array}{ccccc}
       \tmop{Mod}_{\tmop{THH} (k)} (\tmop{Sp}^{g^{<} \mathbb{T}}) &
       \xrightarrow{(-)^{\Phi C_p}} & \tmop{Mod}_{\tmop{THH} (k)^{\Phi C_p}}
       (\tmop{Sp}^{g_p (\mathbb{T}/ C_p)}) &  & \\
       &  & \longdownarrow &  & \\
       &  & \tmop{Mod}_{k^{t C_p}} (\tmop{Sp}^{g_p (\mathbb{T}/ C_p)}) &
       \xrightarrow{(-)^{\eta_{k^{t C_p}} (\mathbb{T}/ C_p)}} & D
       (k^{t\mathbb{T}})
     \end{array} \]
  where the vertical arrow is the base change along the map $\tmop{THH}
  (k)^{\Phi C_p} \rightarrow k^{t C_p}$ of $\mathbb{T}/
  C_p$-$\mathbb{E}_{\infty}$-rings.
\end{lemma}

\begin{remark}
  \label{rem:geom-Cp-fixed-Tate}Recall that both maps in the composite map
  $\underline{\mathbb{Z}} \rightarrow \underline{\mathbb{Z}_p} = \tmop{TR}
  (\mathbb{F}_p) \rightarrow \tmop{THH} (\mathbb{F}_p)$ of
  $\mathbb{T}$-$\mathbb{E}_{\infty}$-rings become equivalences after taking
  geometric fixed points $(-)^{\Phi C_p}$ (cf.~{\cite[Rem~10.9]{Ayala2021}}).
  It follows that we can identify the $\mathbb{T}$-$\mathbb{E}_{\infty}$-ring
  $\tmop{THH} (\mathbb{F}_p)$ with $\underline{\mathbb{Z}}^{\Phi C_p}$.
  Moreover, by {\cite{Hesselholt1997}}, there exists a Bökstedt element
  $\sigma \in \pi_2 \tmop{TF} (\mathbb{F}_p)$ which maps to a Bökstedt
  element in $\pi_2 \tmop{TR}^r (\mathbb{F}_p)$ for every $r \in \mathbb{N}_{>
  0}$, and by {\cite[Lem~II.6.1]{Nikolaus2018}} and computational results
  about $\tmop{THH} (\mathbb{F}_p)$, we can identify the
  $\mathbb{T}$-$\mathbb{E}_{\infty}$-ring $\tmop{THH} (\mathbb{F}_p)
  [\sigma^{- 1}]$ with the Borel $\mathbb{T}$-$\mathbb{E}_{\infty}$-ring
  $\mathbb{Z}^{t C_p}$. Consequently, for every commutative ring $k$, the map
  $\underline{k}^{\Phi C_p} \rightarrow k^{t C_p}$ of
  $\mathbb{T}$-$\mathbb{E}_{\infty}$-rings is simply inverting $\sigma$.
\end{remark}

\begin{lemma}
  \label{lem:base-chg-Fp-tCp-decpl-fixed}The composite functor
  \[ \tmop{Mod}_{\tmop{THH} (\mathbb{F}_p)^{t C_p}} (\tmop{Sp}^{g_p
     (\mathbb{T}/ C_p)}) \longrightarrow \tmop{Mod}_{\mathbb{F}_p^{t C_p}}
     (\tmop{Sp}^{g_p (\mathbb{T}/ C_p)})
     \xrightarrow{(-)^{\eta_{\mathbb{F}_p^{t C_p}} (\mathbb{T}/ C_p)}} D
     (\mathbb{Z}^{t C_p}) \]
  coincides with the forgetful functor, where the first functor is the base
  change along the map $\tmop{THH} (\mathbb{F}_p)^{t C_p} \rightarrow
  \mathbb{F}_p^{t C_p}$ of Borel $\mathbb{E}_{\infty}$-$p$-cyclonic spectra.
\end{lemma}

Recall that the map $\mathbb{Z} \rightarrow \tmop{THH} (\mathbb{F}_p)$ as in
\Cref{rem:geom-Cp-fixed-Tate} also becomes an equivalence after taking $(-)^{t
C_p}$ by {\cite[Cor~IV.4.13]{Nikolaus2018}}, thus we can replace $\tmop{THH}
(\mathbb{F}_p)^{t C_p}$ by $\mathbb{Z}^{t C_p}$ in
\Cref{lem:base-chg-Fp-tCp-decpl-fixed}.

\begin{proof}
  Since both functors in question preserve filtered colimits, it suffices to
  check on compact objects of $\tmop{Mod}_{\mathbb{Z}^{t C_p}} (\tmop{Sp}^{g_p
  (\mathbb{T}/ C_p)})$. Note that the first base change functor preserves
  compact objects. The result follows from \Cref{cor:Z-mod-n-Cn-Tate}.
\end{proof}

\begin{remark}
  The target of the composite functor in
  \Cref{lem:base-chg-Fp-tCp-decpl-fixed} has an
  $\mathbb{F}_p^{t\mathbb{T}}$-module structure, and such a structure is yet
  to explore.
\end{remark}

Summarizing the above discussion, we get

\begin{proposition}
  The composite functor
  \[ \tmop{Mod}_{\tmop{THH} (\mathbb{F}_p)} (\tmop{Sp}^{g^{<} \mathbb{T}})
     \xrightarrow{(-) \otimes_{\tmop{THH} (\mathbb{F}_p)}^{\mathbb{L}}
     \underline{\mathbb{F}_p}} \tmop{Mod}_{\underline{\mathbb{F}_p}}
     (\tmop{Sp}^{g^{<} \mathbb{T}})
     \xrightarrow{(-)^{\theta_{\underline{\mathbb{F}_p}} \mathbb{T}}} D
     (\mathbb{Z}^{t\mathbb{T}}) \]
  coincides with with the composite functor
  \[ \tmop{Mod}_{\tmop{THH} (\mathbb{F}_p)} (\tmop{Sp}^{g^{<} \mathbb{T}})
     \xrightarrow{(-)^{\Phi C_p}} \tmop{Mod}_{\tmop{THH} (\mathbb{F}_p)^{\Phi
     C_p}} (\tmop{Sp}^{g_p (\mathbb{T}/ C_p)}) \longrightarrow D
     (\mathbb{Z}^{t\mathbb{T}}), \]
  where the second functor is inverting $\sigma$ to the underlying $\tmop{THH}
  (\mathbb{F}_p)^{\Phi C_p}$-module.
\end{proposition}

Applying to the $\tmop{THH} (\mathbb{F}_p)$-module spectrum $\tmop{THH}
(\mathcal{C})$, we get:

\begin{corollary}
  \label{cor:poly-HP-THH}Let $\mathcal{C}$ be a dualizable presentable stable
  $\mathbb{F}_p$-linear $\infty$-category. Then the polynomial periodic cyclic
  homology $\tmop{HP}^{\tmop{poly}} (\mathcal{C}/\mathbb{F}_p)$ is equivalent
  to $\tmop{THH} (\mathcal{C}) [\sigma^{- 1}]$ as
  $\mathbb{Z}^{t\mathbb{T}}$-module spectra.
\end{corollary}

\begin{remark}
  We do not compare the multiplicative structures on polynomial periodic
  cyclic homology (as a lax symmetric monoidal functor) with topological
  Hochschild homology in \Cref{cor:poly-HP-THH}.
\end{remark}

Now we give a second proof of \Cref{lem:base-chg-Fp-tCp-decpl-fixed}, inspired
by the proof of {\cite[Prop~2.15]{Mathew2020}}. The key is the following
lemma.

\begin{lemma}
  \label{lem:THH-Fp-tCp-mod-htpy-fixed-pts}The composite functor
  \[ \tmop{Mod}_{\tmop{THH} (\mathbb{F}_p)^{t C_p}} (\tmop{Sp}^{B (\mathbb{T}/
     C_p)}) \longrightarrow \tmop{Mod}_{\mathbb{F}_p^{t C_p}} (\tmop{Sp}^{B
     (\mathbb{T}/ C_p)}) \xrightarrow{(-)^{h (\mathbb{T}/ C_p)}} D (\tmop{TP}
     (\mathbb{F}_p)) \]
  coincides with the forgetful functor.
\end{lemma}

\begin{proof}
  Since the image of $u \in \pi_2 \tmop{TC}^- (\mathbb{F}_p)$ under the
  canonical map $\tmop{TC}^- (\mathbb{F}_p) \rightarrow \tmop{TP}
  (\mathbb{F}_p)$ is $up \in \pi_2 \tmop{TP} (\mathbb{F}_p)$ as reviewed in
  \Cref{rem:THH-perfd}, this composite functor coincides with the composite
  functor
  \[ \tmop{Mod}_{\tmop{THH} (\mathbb{F}_p)^{t C_p}} (\tmop{Sp}^{B (\mathbb{T}/
     C_p)}) \xrightarrow{(-)^{h (\mathbb{T}/ C_p)}} D (\tmop{TP}
     (\mathbb{F}_p)) \xrightarrow{(-) /^{\mathbb{L}} p} D (\tmop{TP}
     (\mathbb{F}_p)), \]
  which is subsequently identified with the forgetful functor by the proof of
  {\cite[Prop~6.4]{Bhatt2018}} (or more precisely, the first displayed formula
  there).
\end{proof}

\Cref{lem:base-chg-Fp-tCp-decpl-fixed} follows from the fact that the
composite functor in \Cref{lem:THH-Fp-tCp-mod-htpy-fixed-pts}, by virtue of
identification with the forgetful functor, preserves filtered colimits, and
that the first functor $\tmop{Mod}_{\tmop{THH} (\mathbb{F}_p)^{t C_p}}
(\tmop{Sp}^{g_p (\mathbb{T}/ C_p)}) \rightarrow \tmop{Mod}_{\mathbb{F}_p^{t
C_p}} (\tmop{Sp}^{g_p (\mathbb{T}/ C_p)})$ in
\Cref{lem:base-chg-Fp-tCp-decpl-fixed} preserves compact objects.

\section{Conjugate filtration}\label{sec:conj-fil}

As explained in the introduction of \Cref{sec:comp-THH}, another
noncommutative analogue of the Cartier isomorphism is the conjugate
filtration, which we will address in this section. As there, our version is
constructed out of $\tmop{THH} (k)$-modules for a base commutative ring $k$,
and by \Cref{lem:poly-HP-formula}, it is crucial to analyze the map
$\tmop{THH} (k)^{\Phi C_p} \rightarrow k^{t C_p}$, endowing the target $k^{t
C_p}$ a suitable filtration.

First, we note that, the homotopy $C_p$-fixed points $k^{h C_p}$ of a
commutative $\mathbb{F}_p$-algebra $k$ is a direct summand of the $C_p$-Tate
construction $k^{t C_p}$.

\begin{remark}
  Let $M$ be an $\mathbb{F}_p$-vector space. Then it follows from computations
  that the composite map
  \[ M^{h C_p} \longrightarrow M^{t C_p} \longrightarrow \tau_{\leqslant 0}
     M^{t C_p} \]
  of $\mathbb{F}_p^{h C_p}$-module spectra is an equivalence. In particular,
  let $k$ be a commutative $\mathbb{F}_p$-algebra. Then the composite map
  \[ k^{h C_p} \longrightarrow k^{t C_p} \longrightarrow \tau_{\leqslant 0}
     k^{t C_p} \]
  of $k^{h C_p}$-module spectra is an equivalence, where the first map has an
  $\mathbb{E}_{\infty}$-structure, and the second map has a $k^{t C_p}$-module
  structure.
\end{remark}

Recall that, in \Cref{sec:comp-THH}, a key fact is that the $\mathbb{T}/
C_p$-equivariant $\mathbb{Z}^{t C_p}$-module spectrum $\mathbb{F}_p^{t C_p}$
is co-induced. It is natural to ask whether the same holds for $k^{t C_p}$ for
any commutative $\mathbb{F}_p$-algebra $k$. We do not know the answer, but
expect it to be false, at least functorially in $k$ (see
\Cref{rem:Cp-Tate-CoInd-functor-mult}). However, note that the
$\mathbb{E}_{\infty}$-ring $k^{t C_p}$ is 2-periodic, and we show that a
``fundamental region'' of $k^{t C_p}$ with respect to its 2-periodicity is
$\mathbb{T}/ C_p$-equivariantly $k$-linearly co-induced.

\begin{notation}
  Let $n$ be a positive integer. Let $\tmop{CoInd}_e^{\mathbb{T}/ C_n}$ denote
  right adjoint $D (\mathbb{Z}) \rightarrow D (\mathbb{Z})^{B (\mathbb{T}/
  C_n)}$ to the (symmetric monoidal) forgetful functor $D (\mathbb{Z})^{B
  (\mathbb{T}/ C_n)} \rightarrow D (\mathbb{Z})$.
\end{notation}

{\construction{\label{cons:htpy-Cp-fixed-CoInd}Let $M$ be a spectrum. Endow
$M$ with trivial $\mathbb{T}$-action, we get a map $\varpi_p^{\ast} M
\rightarrow M^{h C_p}$ in $\tmop{Fun} (B (\mathbb{T}/ C_p), \tmop{Sp})$, where
$\varpi_p \of \mathbb{T}/ C_p \rightarrow \mathord{\ast}$ is the quotient map,
which gives rise to a lax symmetric monoidal lax transformation
\[ \varpi_p^{\ast} (-) \longrightarrow (-)^{h C_p} \]
between lax symmetric monoidal functors $\tmop{Sp} \rightarrow \tmop{Fun} (B
(\mathbb{T}/ C_p), \tmop{Sp})$. On the other hand, there is a non-$\mathbb{T}/
C_p$-equivariant map $M^{h C_p} \rightarrow M$ which, by adjunction, gives
rise to a map $M^{h C_p} \rightarrow \tmop{CoInd}_e^{\mathbb{T}/ C_p} M$ in
$\tmop{Fun} (B (\mathbb{T}/ C_p), \tmop{Sp})$ which is functorial in $M \in
\tmop{Sp}$, and this has a lax symmetric monoidal structure. In summary, we
have a composite lax symmetric natural transformation
\[ \varpi_p^{\ast} (-) \longrightarrow (-)^{h C_p} \longrightarrow
   \tmop{CoInd}_e^{\mathbb{T}/ C_p} (-) \]
of lax symmetric monoidal functors $\tmop{Sp} \rightarrow \tmop{Fun} (B
(\mathbb{T}/ C_p), \tmop{Sp})$.}}

\begin{lemma}
  \label{lem:trunc-htpy-fixed-pts-CoInd}Let $M$ be an $\mathbb{F}_p$-vector
  space. Then the composite map
  \[ \tau_{\geqslant - 1} M^{h C_p} \longrightarrow M^{h C_p} \longrightarrow
     \tmop{CoInd}_e^{\mathbb{T}/ C_p} M \]
  of $\mathbb{T}/ C_p$-equivariant $\mathbb{F}_p$-module\footnote{The
  $\mathbb{F}_p$-module structure comes from lax symmetric monoidal structures
  in Construction~\ref{cons:htpy-Cp-fixed-CoInd}.} spectra is an equivalence,
  where the second map is as in Construction~\ref{cons:htpy-Cp-fixed-CoInd}.
\end{lemma}

\begin{proof}
  Since the maps in question are already constructed, to check that it is an
  equivalence, we could pick a free $\mathbb{Z}$-lift $\tilde{M}$ of $M$,
  namely, a free abelian group $\tilde{M} =\mathbb{Z}^{\oplus I}$ with
  $\tilde{M} / p \cong M$. Then taking $I$-direct sum\footnote{Tate
  construction preserves infinite direct sum of uniformly $t$-bounded
  objects.} of the equivalence in \Cref{cor:Z-mod-n-Cn-Tate} and truncating at
  $[- 1, 0]$, we get an equivalence as an inverse to the composite map in
  question.
\end{proof}

It follows from lax symmetric monoidal structures in
Construction~\ref{cons:htpy-Cp-fixed-CoInd} and
Lemma~\ref{lem:trunc-htpy-fixed-pts-CoInd} that

\begin{corollary}
  \label{cor:0th-assoc-gr-odd-fil}Let $k$ be a commutative
  $\mathbb{F}_p$-algebra. Then there is a canonical equivalence
  \[ \tau_{\geqslant - 1} k^{h C_p} \longrightarrow
     \tmop{CoInd}_e^{\mathbb{T}/ C_p} k \]
  of $\mathbb{T}/ C_p$-equivariant $k$-modules.
\end{corollary}

\begin{remark}
  \label{rem:Cp-Tate-CoInd-functor-mult}Let $M$ be an $\mathbb{F}_p$-vector
  space. It is unclear whether we can {\tmem{functorially}} identify $M^{t
  C_p}$ with the co-induced $\mathbb{T}/ C_p$-equivariant $\mathbb{Z}$-module
  spectrum $\tmop{CoInd}_e^{\mathbb{T}/ C_p} (M^{t\mathbb{T}})$ (although we
  can do it functorially in free $\mathbb{Z}$-lifts $\tilde{M}$), let alone
  $\mathbb{F}_p$-linearly. The multiplicative structure is even more
  complicated.
\end{remark}

Now we describe the filtration on $k^{h C_p}$ and $k^{t C_p}$ as promised.

{\construction{\label{cons:odd-fil-Cp}Let $k$ be a commutative
$\mathbb{F}_p$-algebra. Then we consider the {\tmdfn{odd filtration}} on
$\mathbb{T}/ C_p$-equivariant $\tau_{\geqslant 0} (k^{t C_p})$-module spectra
$k^{h C_p}$, $k^{t C_p}$, and the canonical $\mathbb{T}/ C_p$-equivariant
$\tau_{\geqslant 0} (k^{t C_p})$-module map $k^{h C_p} \rightarrow k^{t C_p}$,
given by the odd parts of the Whitehead filtrations, i.e. for $i \in
\mathbb{Z}$, we have
\begin{eqnarray*}
  \tmop{Fil}_{\tmop{odd}}^i k^{h C_p} & \assign & \tau_{\geqslant 2 i - 1}
  k^{h C_p} ;\\
  \tmop{Fil}_{\tmop{odd}}^i k^{t C_p} & \assign & \tau_{\geqslant 2 i - 1}
  k^{t C_p} .
\end{eqnarray*}
We note that $\tmop{Fil}_{\tmop{odd}}^{> 0} k^{h C_p} = 0$, and the maps
$\tmop{gr}_{\tmop{odd}}^i k^{h C_p} \rightarrow \tmop{gr}_{\tmop{odd}}^i k^{t
C_p}$ are equivalences for $i \leq 0$.}}

\begin{remark}
  \label{rem:trunc-mod-str-desc}Note that $\tmop{gr}_{\tmop{odd}}^0 k^{t C_p}
  = \tau_{\geqslant - 1} k^{h C_p}$ is concentrated in degree $[- 1, 0]$, by
  {\cite[Prop~2.2.1.8]{Lurie2017}}, the $\tau_{\geqslant 0} (k^{t
  C_p})$-module structure on $\tmop{gr}_{\tmop{odd}}^0 k^{t C_p}$ descends
  canonically to a $\tau_{[0, 1]} (k^{t C_p})$-module structure. We note that
  the $\mathbb{T}/ C_p$-equivariant $\mathbb{E}_{\infty}$-$\mathbb{Z}$-algebra
  $\tau_{[0, 1]} (k^{t C_p})$ coincides with the
  $\mathbb{E}_{\infty}$-$\mathbb{Z}$-algebra $k
  \otimes_{\mathbb{Z}}^{\mathbb{L}} \mathbb{F}_p$ with trivial $\mathbb{T}/
  C_p$-action\footnote{This can be seen as follows. Let $M \in D
  (\mathbb{Z})_{\geqslant 0}$ be a connective $\mathbb{Z}$-module spectrum.
  The trivial action on $M$ gives rise to a $\mathbb{T}/ C_p$-equivariant map
  $M \rightarrow M^{h C_p} \rightarrow M^{t C_p}$, and thus a map $M
  \rightarrow \tau_{[0, 1]} (M^{t C_p})$ in $D (\mathbb{Z})^{B (\mathbb{T}/
  C_p)}$. This map gives rise to a lax symmetric monoidal natural
  transformation
  \[ (-) \longrightarrow \tau_{[0, 1]} ((-)^{t C_p}), \]
  and subsequently a lax symmetric monoidal natural transformation
  \[ (-) \otimes_{\mathbb{Z}}^{\mathbb{L}} \mathbb{F}_p \longrightarrow
     \tau_{[0, 1]} ((-)^{t C_p}) \]
  as $\tau_{[0, 1]} (\mathbb{Z}^{t C_p}) =\mathbb{F}_p$, which is an
  equivalence on abelian groups (where we could use the usual resolution for
  cyclic groups to compute $(-)^{t C_p}$). We apply this to the commutative
  $\mathbb{F}_p$-algebra $k$. Alternatively, as the anonymous referee pointed
  out, one could show that $k \otimes_{\mathbb{Z}}^{\mathbb{L}} \mathbb{F}_p$
  is the connective cover of $\tau_{[0, 1]} (k^{t C_p})^{h (\mathbb{T}/
  C_p)}$.}.
\end{remark}

We now produce the conjugate filtration for $\tmop{THH} (k)$-modules. Before
this, we need a base-independence result, which says that the de-completed
Tate construction over $k^{t C_p}$ coincides with that over $\tau_{\geqslant
0} (k^{t C_p})$.

\begin{lemma}
  Let $k$ be a commutative $\mathbb{F}_p$-algebra. Then the natural
  transformation
  \[ (-)^{\eta_{\tau_{\geqslant 0} (k^{t C_p})} (\mathbb{T}/ C_p)}
     \longrightarrow (-)^{\eta_{k^{t C_p}} (\mathbb{T}/ C_p)} \]
  of functors $\tmop{Mod}_{k^{t C_p}} (\tmop{Sp}^{g_p \mathbb{T}}) \rightarrow
  D (k)$, induced by the map $\tau_{\geqslant 0} (k^{t C_p}) \rightarrow k^{t
  C_p}$ in $\tmop{CAlg} (D (k)^{B (\mathbb{T}/ C_p)})$, is an equivalence.
\end{lemma}

\begin{proof}
  It suffices to show that, for every compact generator $M = k^{t C_p} \otimes
  [\mathbb{T}/ C_{p^r}] \in \tmop{Mod}_{k^{t C_p}} (\tmop{Sp}^{g_p
  \mathbb{T}})^{\aleph_0}$, the assembly map
  \[ M^{\eta_{\tau_{\geqslant 0} (k^{t C_p})} (\mathbb{T}/ C_p)}
     \longrightarrow M^{h (\mathbb{T}/ C_p)} \]
  is an equivalence. Recall that the Whitehead tower of $k^{t C_p} =
  (\tau_{\geqslant 0} k^{t C_p}) [\sigma^{- 1}]$ is a sequential colimit of
  shifts of $\tau_{\geqslant 0} (k^{t C_p})$, where $\sigma \in \pi_2 (k^{t
  C_p})$ is a generator. Then the result follows from the fact that $(-)^{h
  (\mathbb{T}/ C_p)}$ preserves weakly Whitehead towers
  (cf.~{\cite[Lem~I.2.6]{Nikolaus2018}}, or the dual of
  {\cite[Lem~3.3]{Bhatt2018}}).
\end{proof}

{\construction{{\dueto{Conjugate filtration}}Let $k$ be a commutative
$\mathbb{F}_p$-algebra, and $M$ a $\tmop{THH} (k)$-module in $\tmop{Sp}^{g^{<}
\mathbb{T}}$. Then by \Cref{lem:poly-HP-formula} and
Construction~\ref{cons:odd-fil-Cp}, the $\mathbb{Z}^{t\mathbb{T}}$-module
spectrum
\[ (M \otimes_{\tmop{THH} (k)}^{\mathbb{L}}
   \underline{k})^{\theta_{\underline{k}} \mathbb{T}} \simeq \left( M^{\Phi
   C_p} \otimes_{\tmop{THH} (k)^{\Phi C_p}}^{\mathbb{L}} \tau_{\geqslant 0}
   (k^{t C_p}) \otimes_{\tau_{\geqslant 0} (k^{t C_p})}^{\mathbb{L}} k^{t C_p}
   \right)^{\eta_{\tau_{\geqslant 0} (k^{t C_p})} (\mathbb{T}/ C_p)} \]
admits a filtration $\tmop{Fil}_{\tmop{conj}}^{\ast}$ induced by the odd
filtration on $k^{t C_p}$.}}

We now analyze the associated graded pieces of the conjugate filtration. It
suffices to analyze the zeroth associated piece:

\begin{remark}
  By 2-periodicity of $k^{t C_p}$, all $\tmop{gr}_{\tmop{odd}}^i k^{t C_p}$'s,
  as $k^{\Phi C_p}$-modules, are isomorphic up to a shift. Thus to study
  associated graded pieces of odd filtered $k^{t C_p}$, it suffices to study
  $\tmop{gr}_{\tmop{odd}}^0 k^{t C_p} \simeq \tmop{gr}_{\tmop{odd}}^i k^{h
  C_p}$, thus the same for the conjugate filtration.
\end{remark}

Since $k \otimes_{\mathbb{Z}}^{\mathbb{L}} \mathbb{F}_p = (\tau_{\geqslant 0}
k^{t C_p}) / \sigma$ is a perfect $\tau_{\geqslant 0} k^{t C_p}$-module in
$\tmop{Sp}^{B (\mathbb{T}/ C_p)}$, we have

\begin{lemma}
  \label{lem:decpl-Tate-k-tCp-trunc}Let $k$ be a commutative
  $\mathbb{F}_p$-algebra. Then the natural transformation
  \[ (-)^{\eta_{\tau_{\geqslant 0} (k^{t C_p})} (\mathbb{T}/ C_p)}
     \longrightarrow (-)^{\eta_{k \otimes_{\mathbb{Z}}^{\mathbb{L}}
     \mathbb{F}_p} (\mathbb{T}/ C_p)} \]
  of functors $\tmop{Mod}_{k \otimes_{\mathbb{Z}}^{\mathbb{L}} \mathbb{F}_p}
  (\tmop{Sp}^{g_p \mathbb{T}}) \rightarrow D (k)$, induced by the map
  $\tau_{\geqslant 0} (k^{t C_p}) \rightarrow k
  \otimes_{\mathbb{Z}}^{\mathbb{L}} \mathbb{F}_p$ in $\tmop{CAlg} (D (k)^{B
  (\mathbb{T}/ C_p)})$, is an equivalence (we tacitly used
  \Cref{rem:borel-cpl-asm-eq}).
\end{lemma}

Similarly, by perfectness of the $k$-module $k
\otimes_{\mathbb{Z}}^{\mathbb{L}} \mathbb{F}_p$, we have

\begin{lemma}
  \label{lem:decpl-Tate-k-mod-p-k}Let $k$ be a commutative
  $\mathbb{F}_p$-algebra. Then the natural transformation
  \[ (-)^{\eta_k (\mathbb{T}/ C_p)} \longrightarrow (-)^{\eta_{k
     \otimes_{\mathbb{Z}}^{\mathbb{L}} \mathbb{F}_p} (\mathbb{T}/ C_p)} \]
  of functors $\tmop{Mod}_{k \otimes_{\mathbb{Z}}^{\mathbb{L}} \mathbb{F}_p}
  (\tmop{Sp}^{g_p \mathbb{T}}) \rightarrow D (k)$, induced by the map $k
  \rightarrow k \otimes_{\mathbb{Z}}^{\mathbb{L}} \mathbb{F}_p$ in
  $\tmop{CAlg} (D (k)^{B (\mathbb{T}/ C_p)})$, is an equivalence (we tacitly
  used \Cref{rem:borel-cpl-asm-eq}).
\end{lemma}

It follows from
\Cref{rem:trunc-mod-str-desc,lem:decpl-Tate-k-tCp-trunc,lem:decpl-Tate-k-mod-p-k}
that

\begin{lemma}
  \label{lem:0th-assoc-gr-conj-fil}Let $k$ be a commutative
  $\mathbb{F}_p$-algebra, and $M$ a $\tmop{THH} (k)$-module in
  $\tmop{Sp}^{g^{<} \mathbb{T}}$. Then the $0$-th associated graded piece
  $\tmop{gr}_{\tmop{conj}}^0  (M \otimes_{\tmop{THH} (k)}^{\mathbb{L}}
  \underline{k})^{\theta_{\underline{k}} \mathbb{T}}$of the conjugate
  filtration is equivalent to
  \[ \left( \left( M^{\Phi C_p} \otimes_{\tmop{THH} (k)^{\Phi
     C_p}}^{\mathbb{L}} (k \otimes_{\mathbb{Z}}^{\mathbb{L}} \mathbb{F}_p)
     \right) \otimes_{k \otimes_{\mathbb{Z}}^{\mathbb{L}}
     \mathbb{F}_p}^{\mathbb{L}} \tmop{gr}_{\tmop{odd}}^0 k^{t C_p}
     \right)^{\eta_k (\mathbb{T}/ C_p)} \]
  in $D (k)$ (again, we tacitly used \Cref{rem:borel-cpl-asm-eq}).
\end{lemma}

Now we embark to understand the composite map $\tmop{THH} (k)
\xrightarrow{\simeq} \tmop{THH} (k)^{\Phi C_p} \rightarrow \underline{k
\otimes_{\mathbb{Z}}^{\mathbb{L}} \mathbb{F}_p}$ of
$\mathbb{T}$-$\mathbb{E}_{\infty}$-rings\footnote{For our purposes, thanks to
\Cref{rem:borel-cpl-asm-eq}, it suffices to understand the underlying
$\mathbb{T}$-equivariant $\mathbb{E}_{\infty}$-ring. However, here we deduce
slightly stronger results.}. In fact, this could be understood for any
commutative ring, not necessarily over $\mathbb{F}_p$.

\begin{remark}
  \label{rem:Tate-diag-geom-fixed-pts}Let $k$ be an animated ring. By
  construction (cf.~{\cite[§IV.2]{Nikolaus2018}}), we have a commutative
  diagram
  \[ \begin{array}{ccccc}
       k & \longrightarrow & \tmop{THH} (k) &  & \\
       \longdownarrow \nocomma \simeq &  & \longdownarrow \nocomma \simeq &  &
       \\
       (N_e^{C_p} k)^{\Phi C_p} & \longrightarrow & \tmop{THH} (k)^{\Phi C_p}
       & \longrightarrow & \underline{k}^{\Phi C_p}\\
       \longdownarrow &  & \longdownarrow &  & \longdownarrow\\
       (k^{\otimes_{\mathbb{S}} p})^{t C_p} & \longrightarrow & \tmop{THH}
       (k)^{t C_p} & \longrightarrow & k^{t C_p}
     \end{array} \]
  where the leftmost row are maps of $\mathbb{E}_{\infty}$-rings while the
  rest are $\mathbb{T}$-$\mathbb{E}_{\infty}$-rings. The composite map $k
  \rightarrow k^{t C_p}$ from the top left to the bottom right is the
  Tate-valued Frobenius.
\end{remark}

\begin{remark}
  Let $k$ be an animated ring. Note that we have a commutative diagram
  \[ \begin{array}{ccc}
       k & \longrightarrow & \underline{k}^{\Phi C_p}\\
       \longdownarrow &  & \longdownarrow\\
       \underline{k} & \xrightarrow{\underline{\tilde{\varphi}_k}} &
       \underline{k \otimes_{\mathbb{Z}}^{\mathbb{L}} \mathbb{F}_p}
     \end{array} \]
  in $\tmop{CAlg} (D (\mathbb{Z}))$, where the top arrow is found in
  \Cref{rem:Tate-diag-geom-fixed-pts}, the bottom arrow is induced by the
  ``extended'' Frobenius map $\tilde{\varphi}_k \of k \rightarrow k
  \otimes_{\mathbb{Z}}^{\mathbb{L}} \mathbb{F}_p \rightarrow k
  \otimes_{\mathbb{Z}}^{\mathbb{L}} \mathbb{F}_p$ of animated rings, and
  except the top left $k$, other terms are
  $\mathbb{T}$-$\mathbb{E}_{\infty}$-rings, and the arrows between
  $\mathbb{T}$-$\mathbb{E}_{\infty}$-rings have
  $\mathbb{T}$-$\mathbb{E}_{\infty}$-structures. By the universal property of
  $\tmop{THH}$ in {\cite{Angeltveit2018}}, we get a commutative diagram
  \[ \begin{array}{ccc}
       \tmop{THH} (k) & \longrightarrow & \underline{k}^{\Phi C_p}\\
       \longdownarrow &  & \longdownarrow\\
       \underline{k} & \xrightarrow{\underline{\tilde{\varphi}_k}} &
       \underline{k \otimes_{\mathbb{Z}}^{\mathbb{L}} \mathbb{F}_p}
     \end{array} \]
  of $\mathbb{T}$-$\mathbb{E}_{\infty}$-rings\footnote{Those who are not
  familiar with $\mathbb{T}$-$\mathbb{E}_{\infty}$-rings can simply replace
  them by $\mathbb{T}$-equivariant $\mathbb{E}_{\infty}$-rings, and the
  $\tmop{THH}$ also has an analogous universal property by
  McClure--Schwänzl--Vogt, cf.~{\cite[Prop~IV.2.2]{Nikolaus2018}}.}.
\end{remark}

So far, we have analyzed the part $M^{\Phi C_p} \otimes_{\tmop{THH} (k)^{\Phi
C_p}}^{\mathbb{L}} (k \otimes_{\mathbb{Z}}^{\mathbb{L}} \mathbb{F}_p)$ of the
expression in \Cref{lem:0th-assoc-gr-conj-fil}. Now we analyze the part $((-)
\otimes_{k \otimes_{\mathbb{Z}}^{\mathbb{L}} \mathbb{F}_p}^{\mathbb{L}}
\tmop{gr}_{\tmop{odd}}^0 k^{t C_p})^{\eta_k (\mathbb{T}/ C_p)}$ for
commutative $\mathbb{F}_p$-algebras $k$. It follows from
\Cref{cor:0th-assoc-gr-odd-fil} that this simplifies to $(-) \otimes_{k
\otimes_{\mathbb{Z}}^{\mathbb{L}} \mathbb{F}_p}^{\mathbb{L}} k$, where the $k
\otimes_{\mathbb{Z}}^{\mathbb{L}} \mathbb{F}_p$-module structure on $k$ is
simply given by the multiplication map $k \otimes_{\mathbb{Z}}^{\mathbb{L}}
\mathbb{F}_p \rightarrow k$. Recall that, for animated $\mathbb{F}_p$-algebras
$k$, the composite map
\[ k \xrightarrow{\tilde{\varphi}_k} k \otimes_{\mathbb{Z}}^{\mathbb{L}}
   \mathbb{F}_p \longrightarrow k \]
coincides with the usual Frobenius map $\varphi_k \of k \rightarrow k$. The
above discussion implies that

\begin{proposition}
  Let $k$ be a commutative $\mathbb{F}_p$-algebra, and $M$ a $\tmop{THH}
  (k)$-module in $\tmop{Sp}^{g^{<} \mathbb{T}}$. Then the $0$-th associated
  graded piece $\tmop{gr}_{\tmop{conj}}^0  (M \otimes_{\tmop{THH}
  (k)}^{\mathbb{L}} \underline{k})^{\theta_{\underline{k}} \mathbb{T}}$of the
  conjugate filtration is equivalent to the Frobenius twist
  \[ (M^{\Phi C_p} \otimes_{\tmop{THH} (k)}^{\mathbb{L}} k) \otimes_{k,
     \varphi_k}^{\mathbb{L}} k \]
  in $D (k)$, where the $\tmop{THH} (k)$-module structure on $M^{\Phi C_p}$ is
  induced by the equivalence $\tmop{THH} (k) \xrightarrow{\simeq} \tmop{THH}
  (k)^{\Phi C_p}$.
\end{proposition}

\begin{corollary}
  Let $k$ be a commutative $\mathbb{F}_p$-algebra, and $\mathcal{C}$ a
  dualizable presentable stable $\mathbb{F}_p$-linear $\infty$-category. Then
  the $0$-th associated graded piece $\tmop{gr}_{\tmop{conj}}^0
  \tmop{HP}^{\tmop{poly}} (\mathcal{C}/ k)$ is equivalent to $\tmop{HH}
  (\mathcal{C}/ k) \otimes_{k, \varphi_k}^{\mathbb{L}} k$ in $D (k)$.
\end{corollary}

Finally, we show that the conjugate filtration is complete. It suffices to
establish the following connectivity result.

\begin{lemma}
  \label{lem:connectivity-conj-fil}Let $k$ be a commutative
  $\mathbb{F}_p$-algebra, and $M$ a $\tmop{THH} (k)^{\Phi C_p}$-module in
  $\tmop{Sp}^{g^{<} (\mathbb{T}/ C_p)}$. Suppose that $M$ is connective. Then
  the spectrum
  \[ \left( M \otimes_{\tmop{THH} (k)^{\Phi C_p}}^{\mathbb{L}} \tau_{\geqslant
     - 1} k^{t C_p} \right)^{\eta_{\tau_{\geqslant 0} (k^{t C_p})}
     (\mathbb{T}/ C_p)} \]
  is connective.
\end{lemma}

Note that $k^{t C_p}$ is co-induced in $D (\mathbb{Z})$. Although it might not
be compatible with the $\tmop{THH} (k)^{\Phi C_p}$-module structure, this is
already enough for us to prove the connectivity.

\begin{proof}
  By the bar resolution of the relative tensor product, and the fact that the
  de-completed Tate construction preserves small colimits, it suffices to show
  that, for every $n \in \mathbb{N}$, the spectrum
  \[ (M \otimes_{\underline{\mathbb{Z}}}^{\mathbb{L}} (\tmop{THH} (k)^{\Phi
     C_p})^{\otimes_{\underline{\mathbb{Z}}}^{\mathbb{L}} n}
     \otimes_{\underline{\mathbb{Z}}}^{\mathbb{L}} \tau_{\geqslant - 1} k^{t
     C_p})^{\eta_{\tau_{\geqslant 0} (k^{t C_p})} (\mathbb{T}/ C_p)} \]
  is connective. We write $k$ as a direct sum $\mathbb{F}_p^{\oplus I}$, and
  applying \Cref{cor:Z-mod-n-Cn-Tate}, we see that this spectrum is equivalent
  to the direct sum of $I$'s copies of
  \[ M \otimes_{\mathbb{Z}}^{\mathbb{L}} (\tmop{THH} (k)^{\Phi
     C_p})^{\otimes_{\mathbb{Z}}^{\mathbb{L}} n}
     \otimes_{\mathbb{Z}}^{\mathbb{L}} \tau_{\geqslant 0} \mathbb{Z}^{t C_p}
  \]
  which is connective.
\end{proof}

By 2-periodicity of $k^{t C_p}$, and the fact that $(-)^{\Phi C_p}$ preserves
connectivity, we deduce from \Cref{lem:connectivity-conj-fil} that

\begin{corollary}
  Let $k$ be a commutative $\mathbb{F}_p$-algebra, and $M$ a connective
  $\tmop{THH} (k)$-module in $\tmop{Sp}^{g^{<} \mathbb{T}}$. Then the
  conjugate filtration on $(M \otimes_{\tmop{THH} (k)}^{\mathbb{L}}
  k)^{\eta_{\underline{k}} \mathbb{T}}$ is complete.
\end{corollary}

\begin{corollary}
  \label{cor:conj-fil-complete-bnd-below-cats}Let $k$ be a commutative
  $\mathbb{F}_p$-algebra, and $\mathcal{C}$ a dualizable presentable stable
  $\mathbb{F}_p$-linear $\infty$-category. Suppose that its topological
  Hochschild homology $\tmop{THH} (\mathcal{C})$ is bounded
  below\footnote{This is the case when $\mathcal{C}= D (R)$ for a $(-
  1)$-connective $\mathbb{E}_1$-$k$-algebra $R$, or $\mathcal{C}= D (X)$ for a
  quasicompact quasiseparated $k$-scheme $X$.}. Then the conjugate filtration
  on its polynomial periodic cyclic homology $\tmop{HP}^{\tmop{poly}}
  (\mathcal{C}/ k)$ is complete.
\end{corollary}

\appendix\section{Tate cohomology complex}\label{app:Tate-cohom-cplx}

We briefly show that our de-completed Tate cohomology on Lazard-semi-flat
chain complexes of cohomological Mackey functors can be computed by
{\tmdfn{Tate cohomology complex}} in {\cite[§6.2]{Kaledin2015a}}, or
{\cite[§2.1]{Petrov2018}}, thus it is a homotopy invariant version of the
latter.

Let $k$ be a commutative ring. Recall that a {\tmdfn{finitely generated
permutation $G$-module}} is a (left) $k [G]$-module of the form $k [X]$ for
some {\tmem{finite}} $G$-set $X$. We denote by $\tmop{Perm}_G (k) \subseteq
\tmop{LMod}_{k [G]}$ the full subcategory spanned by finitely generated
permutation $G$-modules. We refer to {\cite{Balmer2021}} for comparison of
different characterizations of cohomological Mackey functors, and
{\cite[Ex~2.5]{Brantner2021}} for the equivalence $D
(\tmop{Mack}_G^{\tmop{coh}} (k)) = \tmop{Mod}_{\underline{k}} (\tmop{Sp}^{g
G})$ between derived cohomological $G$-Mackey functors and
$\underline{k}$-modules in genuine $G$-spectra.

Now we start with a flatness of cohomological Mackey functors.

\begin{definition}
  Let $k$ be a commutative ring, and $G$ a finite group.
  \begin{itemize}
    \item We say that a cohomological Mackey functor is {\tmdfn{Lazard-flat}}
    if it is a filtered colimit of finitely generated permutation modules.
    
    \item We say that an (unbounded) chain complex of cohomological Mackey
    functors is {\tmdfn{Lazard-semi-flat}} if it is a filtered colimit of
    bounded chain complexes of finitely generated permutation modules. Compare
    with {\cite[Thm~1.1]{Christensen2015}}.
  \end{itemize}
\end{definition}

\begin{remark}
  \label{rem:ind-perm-coh-mack-funs}Let $k$ be a commutative ring, and $G$ a
  finite group. Then permutation $k [G]$-modules are compact in the category
  $\tmop{Mack}_G^{\tmop{coh}} (k)$ of cohomological $G$-Mackey functors.
  Consequently, we have a canonical fully faithful functor
  \[ \tmop{Ind} (\tmop{Perm}_G) \longhookrightarrow \tmop{Mack}_G^{\tmop{coh}}
     (k) . \]
  In particular, we can identify Lazard-flat (derived) cohomological Mackey
  functors as objects of $\tmop{Ind} (\tmop{Perm}_G (k))$.
  
  Similarly, recall that a chain complex in an additive category is compact if
  and only if it is bounded and degreewise compact
  {\cite[Thm~4.5]{Christensen2015}}. The preceding argument shows that we have
  a canonical fully faithful functor
  \[ \tmop{Ind} (\tmop{Ch}^b (\tmop{Perm}_G (k))) \longhookrightarrow
     \tmop{Ch} (\tmop{Mack}_G^{\tmop{coh}} (k)) \]
  and thus we can identify Lazard-semi-flat chain complexes of cohomological
  Mackey functors as objects of $\tmop{Ind} (\tmop{Ch}^b (\tmop{Perm}_G
  (k)))$.
\end{remark}

Lazard-flat cohomological Mackey functors (resp. Lazard-semi-flat chain
complexes of cohomological Mackey functors) are in fact $k [G]$-modules (resp.
chain complexes of $k [G]$-modules):

\begin{remark}
  \label{rem:ind-perm-G-mods}Let $k$ be a commutative ring, and $G$ a finite
  group. Then finitely generated permutation $k [G]$-modules are compact in
  the category $\tmop{LMod}_{k [G]}$ of left $G$-modules, thus the fully
  faithful functor
  \[ \tmop{Ind} (\tmop{Perm}_G (k)) \longhookrightarrow
     \tmop{Mack}_G^{\tmop{coh}} (k) \]
  in \Cref{rem:ind-perm-coh-mack-funs} factors through the inclusion
  $\tmop{LMod}_{k [G]} \hookrightarrow \tmop{Mack}_G^{\tmop{coh}} (k)$, and
  the fully faithful functor
  \[ \tmop{Ind} (\tmop{Ch}^b (\tmop{Perm}_G (k))) \longhookrightarrow
     \tmop{Ch} (\tmop{Mack}_G^{\tmop{coh}} (k)) \]
  factors through the inclusion $\tmop{Ch} (\tmop{LMod}_{k [G]})
  \hookrightarrow \tmop{Ch} (\tmop{Mack}_G^{\tmop{coh}} (k))$.
\end{remark}

Now we review the Tate cohomology complex. Let $k$ be a commutative ring, and
$G$ a finite group. A {\tmdfn{complete resolution}} {\cite[§2.1]{Petrov2018}}
of the left $k [G]$-module $k$ is an acyclic complex $P_{\ast} \in \tmop{Ch}
(\tmop{LMod}_{k [G]}^{\tmop{free}})$ of free left $k [G]$-modules along with
an isomorphism $\varepsilon \of \mathbb{Z} \rightarrow \ker (\mathd \of P_0
\rightarrow P_{- 1})$ of left $k [G]$-modules.

\begin{definition}
  Let $k$ be a commutative ring, $G$ a finite group, and $P_{\ast}$ a complete
  resolution of $k$. The functor of the {\tmdfn{Tate cohomology complex}} is
  defined to be the functor
  \begin{eqnarray*}
    \tmop{Ch} (\tmop{Mack}_G^{\tmop{coh}} (k)) & \longrightarrow & \tmop{Ch}
    (k),\\
    M_{\ast} & \longmapsto & (\tmop{Tot}^{\oplus} (M_{\ast} \otimes_k
    P_{\ast}))^G .
  \end{eqnarray*}
\end{definition}

On Lazard-semi-flat chain complexes of cohomological Mackey functors, the Tate
cohomology complex represents the de-completed Tate cohomology. More
precisely, we have

\begin{proposition}
  Let $k$ be a commutative ring, and $G$ a finite group. Then we have a
  commutative diagram
  \[ \begin{array}{ccc}
       \tmop{Ind} (\tmop{Ch}^b (\tmop{Perm}_G (k))) & \longrightarrow &
       \tmop{Ch} (k)\\
       \longdownarrow &  & \longdownarrow\\
       D (\tmop{Mack}_G^{\tmop{coh}} (k)) = \tmop{Mod}_{\underline{k}}
       (\tmop{Sp}^{g G}) & \xrightarrow{(-)^{\theta_{\underline{k}} G}} & D
       (k)
     \end{array} \]
  of additive $\infty$-categories, where the top horizontal arrow is the
  composite of the inclusion in \Cref{rem:ind-perm-coh-mack-funs} and the Tate
  cohomology complex, and the vertical arrows are localizations at
  quasi-isomorphisms.
\end{proposition}

\begin{proof*}{Sketch of proof}
  Since all functors in question preserve filtered colimits, we can replace
  the top left term $\tmop{Ind} (\tmop{Ch}^b (\tmop{Perm}_G (k)))$ by
  $\tmop{Ch}^b (\tmop{Perm}_G (k))$. By stability, we could further restrict
  to the full subcategory $\tmop{Ch}_{\geqslant 0}^b (\tmop{Perm}_G (k))$ of
  bounded chain complexes concentrated in non-negative degrees. This is the
  same as freely adjoining finitary-geometric-realizations to $\tmop{Perm}_G
  (k)$. All functors in question preserve finitary geometric realizations,
  thus we can replace the top left term simply by $\tmop{Perm}_G (k)$.
  
  In this case, the left vertical arrow is fully faithful, with essential
  image being compact in $\tmop{Mod}_{\underline{k}} (\tmop{Sp}^{g G})$.
  Therefore we may replace $(-)^{\theta_{\underline{k}} G}$ by $(-)^{t G}$. It
  remains to show that, for every finitely generated permutation $G$-module
  $M$, the complex $(M \otimes_k P_{\ast})^G$ represents $M^{t G}$, which
  follows from definition.
\end{proof*}

\section{Tate construction with respect to cyclic $p$-groups}

Recall that in {\cite[Lem~II.4.1~\&~II.4.2]{Nikolaus2018}}, they showed that,
for every $t$-bounded below $\mathbb{T}$-equivariant spectrum $X$, the
canonical map
\[ X^{t\mathbb{T}} \longrightarrow \lim_r X^{t C_{p^r}} \longrightarrow \lim_r
   (X^{t C_p})^{h (C_{p^r} / C_p)} \]
exhibits the target as the $p$-completion of the source. We need a variant of
this, which allows us to establish the equivalence for non-$t$-bounded objects
which are linear over a $t$-bounded object. We start with two simple
observations.

\begin{lemma}
  \label{lem:Tate-orb-Z-mods}Let $r \in \mathbb{N}_{> 0}$ be a natural number,
  and $X \in D (\mathbb{Z})^{B C_{p^r}}$ a $C_{p^r}$-equivariant
  $\mathbb{Z}$-module spectrum. Then the canonical map
  \[ X^{t C_{p^r}} \longrightarrow (X^{t C_p})^{h (C_{p^r} / C_p)} \]
  is an equivalence, and that $X^{t C_{p^r}}$ is $p$-complete.
\end{lemma}

\begin{proof}
  As in the proof of {\cite[Lem~II.4.1]{Nikolaus2018}}, we examine the map
  \[ \begin{array}{ccccc}
       X_{h C_{p^r}} & \longrightarrow & X^{h C_{p^r}} & \longrightarrow &
       X^{t C_{p^r}}\\
       \longdownarrow &  & \longdownarrow &  & \longdownarrow\\
       (X_{h C_p})^{h (C_{p^r} / C_p)} & \longrightarrow & (X^{h C_p})^{h
       (C_{p^r} / C_p)} & \longrightarrow & (X^{t C_p})^{h C_{p^{r - 1}}}
     \end{array} \]
  of horizontal fiber sequences, where the middle vertical map is an
  equivalence. The cofiber of the leftmost vertical map is $(X_{h C_p})^{t
  (C_{p^r} / C_p)}$, which is zero by (an iterated application of) the Tate
  orbit lemma for $\mathbb{T}$-equivariant $\mathbb{Z}$-module spectra,
  {\cite[the footnote in Ex~I.2.3]{Nikolaus2018}}.
  
  For the $p$-completeness of $X^{t C_{p^r}}$, by the first part, it suffices
  to show the $p$-completeness of $X^{t C_p}$. Note that it is a
  $\mathbb{Z}^{t C_p}$-module by lax symmetric monoidal structure on $(-)^{t
  C_p}$ as in {\cite[§I.3]{Nikolaus2018}}, and we know that there is an
  equivalence $\tmop{THH} (\mathbb{F}_p) \simeq \tau_{\geqslant 0}
  (\mathbb{Z}^{t C_p})$ of $\mathbb{E}_{\infty}$-rings by
  {\cite[Cor~IV.4.16]{Nikolaus2018}}. This gives rise to an
  $\mathbb{F}_p$-module structure on the spectrum $X^{t C_p}$ (which might not
  be compatible with the $\mathbb{Z}$-module structure), thus the spectrum
  $X^{t C_p}$ is $p$-complete.
\end{proof}

\begin{lemma}
  \label{lem:lim-formula-Z-mods-tT}Let $X \in D (\mathbb{Z})^{B\mathbb{T}}$ be
  a $\mathbb{T}$-equivariant $\mathbb{Z}$-module spectrum. Then the canonical
  map
  \[ X^{t\mathbb{T}} \longrightarrow \lim_r X^{t C_{p^r}} \]
  exhibits the target as the $p$-completion of the source.
\end{lemma}

\begin{proof}
  The proof of {\cite[Lem~IV.4.12]{Nikolaus2018}} shows that this map, as a
  map of spectra, can be rewritten as
  \[ X^{t\mathbb{T}} \longrightarrow \lim_r  (X^{t\mathbb{T}} / p^r) \]
  and the result follows immediately.
\end{proof}

Now we extend these result to modules over a $t$-bounded ring.

\begin{proposition}
  \label{prop:Tate-orb-bnd-mods}Let $V \in \tmop{Sp}^{\omega}$ be a finite
  spectrum, $r \in \mathbb{N}_{> 0}$ a natural number, and $A \in \tmop{Sp}^{B
  C_{p^r}}$ a $C_{p^r}$-equivariant $\mathbb{E}_1$-ring such that $A \otimes
  V$ is $t$-bounded above, and $M \in \tmop{RMod}_A \left( \tmop{Sp}^{B
  C_{p^r}} \right)$ a $C_{p^r}$-equivariant right $A$-module. Then the
  canonical map
  \[ M^{t C_{p^r}} \longrightarrow (M^{t C_p})^{h (C_{p^r} / C_p)} \]
  is an equivalence after Bousfield $V$-localization, and $M^{t C_{p^r}}$ is
  $p$-complete after Bousfield $V$-localization.
\end{proposition}

\begin{proof}
  First, we examine the fiber sequence
  \[ (\tau_{\geqslant 0} A) \otimes V \longrightarrow A \otimes V
     \longrightarrow (\tau_{< 0} A) \otimes V. \]
  The second term is $t$-bounded above by assumption, and so is the third term
  by finiteness of $V$. It follows that the first term is $t$-bounded above
  and thus $t$-bounded. We may replace $A$ by $\tau_{\geqslant 0} A$, and
  without loss of generality, we may assume that $A$ is connective.
  
  Let $\mathcal{F} \subseteq \tmop{Sp}^{B C_{p^r}}$ be the collection of $X
  \in \tmop{Sp}^{B C_{p^r}}$ such that the canonical map
  \[ X^{t C_{p^r}} \longrightarrow (X^{t C_p})^{h (C_{p^r} / C_p)} \]
  is an equivalence, and that $M^{t C_{p^r}}$ is $p$-complete. Since functors
  $(-)^{t C_{p^r}}$, $((-)^{t C_p})^{h (C_{p^r} / C_p)}$, and $(-)_p^{\wedge}$
  are exact, the collection $\mathcal{F}$ is closed under taking cofibers,
  translations, and retracts. By \Cref{lem:Tate-orb-Z-mods}, the collection
  $\mathcal{F}$ contains every $C_{p^r}$-equivariant spectrum which admits a
  $\mathbb{Z}$-module structure.
  
  We now show that $M \otimes V \in \mathcal{F}$. As in the proof of
  \Cref{lem:zero-ring-mods}, the $C_{p^r}$-equivariant spectrum $M \otimes V$
  is a retract of $M \otimes A \otimes V$, thus it suffices to show that $M
  \otimes A \otimes V \in \mathcal{F}$. Now $A \otimes V$ is $t$-bounded, thus
  it is a finite extension of $C_{p^r}$-equivariant $\mathbb{Z}$-modules. The
  result then follows.
\end{proof}

\begin{proposition}
  \label{prop:lim-formula-bnd-mods-tT}Let $V \in \tmop{Sp}^{\omega}$ be a
  finite spectrum, and $A \in \tmop{Sp}^{B\mathbb{T}}$ a
  $\mathbb{T}$-equivariant $\mathbb{E}_1$-ring such that $A \otimes V$ is
  $t$-bounded above, and $M \in \tmop{RMod}_A (\tmop{Sp}^{B\mathbb{T}})$ a
  $\mathbb{T}$-equivariant right $A$-module. Then the canonical map
  \[ M^{t\mathbb{T}} \longrightarrow \lim_r M^{t C_{p^r}} \]
  is exhibits the target as a $p$-completion of the source after Bousfield
  $V$-localization.
\end{proposition}

\begin{proof}
  As in the proof of \Cref{prop:Tate-orb-bnd-mods}, we may assume that $A$ is
  connective. Moreover, by \Cref{prop:Tate-orb-bnd-mods}, the limit $\lim_r 
  (M \otimes V)^{t C_{p^r}}$ is $p$-complete. Let $\mathcal{F} \subseteq
  \tmop{Sp}^{B\mathbb{T}}$ be the collection of $X \in
  \tmop{Sp}^{B\mathbb{T}}$ such that the canonical map
  \[ (M \otimes V)^{t\mathbb{T}} / p \longrightarrow \lim_r  (M \otimes V)^{t
     C_{p^r}} / p \]
  is an equivalence. Then $\mathcal{F}$ is closed under taking cofibers,
  translations, and retracts. By \Cref{lem:lim-formula-Z-mods-tT}, the
  collection $\mathcal{F}$ contains every $\mathbb{T}$-equivariant spectrum
  which admits a $\mathbb{Z}$-module structure.
  
  Now by assumption, the $\mathbb{T}$-equivariant spectrum $M \otimes A
  \otimes V$ is a finite extension of $\mathbb{T}$-equivariant
  $\mathbb{Z}$-modules, and as $M \otimes V$ is a retract of $M \otimes A
  \otimes V$, we deduce that $M \otimes V \in \mathcal{F}$.
\end{proof}

\

\end{document}